# Foliated Plateau problems and asymptotic counting of surface subgroups

# Problèmes de Plateau feuilletés et comptage asymptotique des sous-groupes de surfaces

Sébastien Alvarez[*], Ben Lowe[†], Graham Smith[‡]

**Abstract:** In 2000, Labourie initiated the study of the dynamical properties of the space of $k$-surfaces, that is, suitably complete immersed surfaces of constant extrinsic curvature in 3-dimensional manifolds, which he presented as a higher-dimensional analogue of the geodesic flow when the ambient manifold is negatively curved. In this paper, following the recent work of Calegari–Marques–Neves, we study the asymptotic counting of surface subgroups in terms of areas of $k$-surfaces. We determine a lower bound, and we prove rigidity when this bound is achieved. Our work differs from that of Calegari–Marques–Neves in two key respects. Firstly, we work with all quasi-Fuchsian subgroups as opposed to merely asymptotically Fuchsian ones. Secondly, as their proof of rigidity breaks down in the present case, we require a different approach. Following ideas recently outlined by Labourie, we prove rigidity by solving a general foliated Plateau problem in Cartan–Hadamard manifolds. To this end, we build on Labourie's theory of $k$-surface dynamics, and propose a number of new constructions, conjectures and questions.

**Résumé :** Dans les années 2000, Labourie a commencé l'étude des propriétés dynamiques de l'espace des $k$-surfaces, c'est-à-dire des surfaces complètes à courbure extrinsèque constante dans les 3-variétés à courbure négative, qu'il présente comme des analogues du flot géodésique en dimension supérieure. Dans cet article, en suivant les travaux récents de Calegari–Marques–Neves, nous étudions le comptage asymptotique de sous-groupes de surfaces selon l'aire des $k$-surfaces qui les représentent. Nous établissons une borne inférieure, et prouvons un résultat de rigidité lorsque le minimum est atteint. Notre travail se distingue de celui de Calegari–Marques–Neves en deux aspects. Premièrement, nous considérons tous les sous-groupes quasi-fuchsiens et non pas uniquement ceux qui sont asymptotiquement fuchsiens. Deuxièmement, leur preuve de la rigidité ne s'applique pas dans notre cadre et nous suivons une autre approche. Nous établissons la rigidité en résolvant un problème de Plateau feuilleté général dans les variétés de Cartan–Hadamard. Pour ce faire, nous nous basons sur la théorie dynamique des k-surfaces développée par Labourie, et proposons quelques nouvelles constructions, conjectures et questions.

## 1 - Introduction.

**1.1 - Asymptotic counting.** In his pioneering work [22], Labourie initiated a dynamical systems approach to the study of suitably complete surfaces of constant extrinsic curvature equal to $k > 0$ (henceforth referred to as $k$-surfaces). Indeed, let $X := (X, h)$ be an oriented 3-dimensional Cartan–Hadamard manifold of sectional curvature bounded above by $-1$. Suppose furthermore that $X$ is acted on cocompactly by the group $\Pi$, say, and that the quotient $X/\Pi$ is a manifold. Given $0 < k < 1$, Labourie introduced the space $\mathcal{S}_k(X/\Pi)$ of marked $k$-surfaces in $X/\Pi$. He showed how this space possesses a natural laminated structure, how it may be viewed as an extension of the geodesic lamination, and how it possesses the hyperbolic properties of the latter. In addition, he proposed the study of the asymptotic growth rate and possible rigidity properties as $E$ tends to infinity of the number $N(E, X/\Pi)$ of compact elements of $\mathcal{S}_k(X/\Pi)$ of *energy* bounded above by $E$ (see Section 1.2 of [22] and Section 3.5.1 of [23]). This hard but intriguing problem remains open. In the present paper, inspired by the recent work [8] of Calegari–Marques–Neves, we address a simpler, but related, asymptotic counting problem expressed in terms of areas of compact elements of $\mathcal{S}_k(X/\Pi)$.

---

[*] CMAT, Facultad de Ciencias, Universidad de la República, Montevideo, Uruguay
[†] Department of Mathematics, University of Chicago, Chigaco IL, USA
[‡] Pontifícia Universidade Católica do Rio de Janeiro (PUC-Rio), Rio de Janeiro, Brazil

In order to state our result, we first note that

$$\partial_\infty X = \partial_\infty \Pi = \partial_\infty \mathbb{H}^3 = \hat{\mathbb{C}}. \tag{1.1}$$

For $k > 0$, we define a *k-disk* in $X$ to be an oriented, smoothly embedded surface $D \subseteq X$, of constant extrinsic curvature equal to $k$, whose closure $\overline{D} \subseteq X \cup \hat{\mathbb{C}}$ is homeomorphic to the closed unit disk in $\mathbb{R}^2$. Note that the boundary

$$\partial_\infty D := \overline{D} \setminus D \tag{1.2}$$

of any $k$-disk is an oriented Jordan curve in $\hat{\mathbb{C}}$. Conversely, in Theorem 2.1.1, we show that, for all $0 < k < 1$, every oriented Jordan curve $c$ in $\hat{\mathbb{C}}$ defines a unique $k$-disk $D := \mathrm{D}_{k,h}(c)$ in $X$.

Recall now that a compact surface subgroup $\Gamma \subseteq \Pi$ is said to be *quasi-Fuchsian* whenever its limit set $\partial_\infty \Gamma$ is a quasicircle. Let QF denote the set of conjugacy classes of quasi-Fuchsian subgroups of $\Pi$. For all $0 < k < 1$, we define the function $\mathrm{A}_{k,h} : \mathrm{QF} \to \mathbb{R}$ by

$$\mathrm{A}_{k,h}([\Gamma]) := \mathrm{Area}(\mathrm{D}_{k,h}(\partial_\infty \Gamma)/\Gamma). \tag{1.3}$$

**Theorem 1.1.1**

*For every metric $h$ over $X/\Pi$ of sectional curvature bounded above by $-1$,*

$$\liminf_{A \to \infty} \frac{\log\left(\#\{[\Gamma] \in \mathrm{QF} \mid \mathrm{A}_{k,h}(\Gamma) \leqslant A\}\right)}{A \log(A)} \geqslant \frac{(1-k)}{2\pi}, \tag{1.4}$$

*with equality holding if and only if $h$ is hyperbolic.*

In fact, we prove a stronger, more technical, result, namely Theorem 5.1.3, which is closer in spirit to that of Calegari–Marques–Neves, and from which Theorem 1.1.1 follows immediately.

Theorem 1.1.1 is best understood within the context of the asymptotic counting result proven by Kahn–Marković in [19], namely

$$\lim_{g \to \infty} \frac{\log\left(\#\{[\Gamma] \in \mathrm{QF} \mid \mathrm{g}(\Gamma) \leqslant g\}\right)}{2g \log(g)} = 1, \tag{1.5}$$

where, for all $\Gamma \in \mathrm{QF}$, $\mathrm{g}(\Gamma)$ denotes its genus. It is worth noting that (1.5) is, in fact, the lesser of two asymptotic counting results proven by Kahn–Marković in that paper. Their second result counts, not conjugacy classes, but commensurability classes. Since commensurability classes do not distinguish covers of the same surface, whilst conjugacy classes do, a rigidity result analogous to Theorem 1.1.1 for the former would be of greater geometric interest. However, such a result lies beyond the scope of currently available techniques.

We view (1.5) as a *topological* asymptotic counting result, in the sense that it only concerns the structure of the fundamental group $\Pi$. In [8], Calegari–Marques–Neves make the key observation that minimal surfaces may be used to refine (1.5) to a *geometric* asymptotic counting result, that is, one which involves the structure of the metric $h$. In this manner, they obtain a rigidity result in the spirit of the volume entropy rigidity result [5] of Besson–Courtois–Gallot. Indeed, rigidity in Theorem 1.1 of [8] and Theorem 1.1.1 states that geometric asymptotic counting is greater than topological asymptotic counting, with equality holding if and only if the metric is hyperbolic. Note that Calegari–Marques–Neves' asymptotic counting result in fact differs from that of Kahn–Marković in that it involves a double limit – a byproduct of the various geometric challenges that arise when working with minimal surfaces. Since $k$-disks have simpler geometric properties, only a single limit is required in Theorem 1.1.1, allowing us to state a result closer in spirit to that of Kahn–Marković.

Foliated Plateau problems and asymptotic counting of surface subgroups

**1.2 - Foliated Plateau problems.** Theorem 1.1.1 differs substantially from the corresponding result of Calegari–Marques–Neves in that its proof requires an extensive study of the dynamical properties of the space of $k$-disks. In the spirit of Gromov's theory of foliated Plateau problems – laid out in [12] and [13] and applied to $k$-surfaces by Labourie in [22] and [23] – we are led to view $k$-disks "not as individuals, but as members of a community", and then to investigate the dynamical properties of such families. We find the resulting theory to be interesting in its own right, and we thus also discuss certain new constructions and open problems which lie beyond our immediate needs for the study of asymptotic counting.

We describe here two successive refinements of Labourie's laminated structure obtained upon restricting to successively smaller subsets of $\mathcal{S}_k(X/\Pi)$. For ease of presentation, we only address the case where $X/\Pi$ is compact. More general statements of the following results are given in Sections 2 and 3.

Our first refinement yields a space of marked $k$-surfaces carrying a natural *fibred* structure. For all $0 < k < 1$, we define a *marked $k$-disk* in $X$ to be a pair $(D, p)$, where $D$ is a $k$-disk in $X$ and $p$ is a point of $D$. For all $0 < k < 1$, and for all $K \geqslant 1$, let $\mathrm{MKD}_{k,h}(K)$ denote the space of marked $k$-disks $(D, p)$ in $X$ whose ideal boundary $\partial_\infty D$ is a $K$-quasicircle in $\hat{\mathbb{C}}$, and denote

$$\mathrm{MKD}_{k,h} := \bigcup_{K \geqslant 1} \mathrm{MKD}_{k,h}(K). \tag{1.6}$$

For all $K \geqslant 1$, let $\mathrm{QC}^+(K)$ denote the space of oriented $K$-quasicircles in $\hat{\mathbb{C}}$, and denote

$$\mathrm{QC}^+ := \bigcup_{K \geqslant 1} \mathrm{QC}^+(K). \tag{1.7}$$

These spaces carry natural topologies described respectively in Sections 3.2 and 3.1. In particular, with respect to these topologies, the boundary map $\partial_\infty : \mathrm{MKD}_{k,h} \to \mathrm{QC}^+$ is continuous.

**Theorem 1.2.1, Fibred Plateau problem**

*For all $0 < k < 1$, $(\mathrm{MKD}_{k,h}, \partial_\infty)$ is a topological associated $\mathrm{PSL}(2, \mathbb{R})$-bundle over $\mathrm{QC}^+$.*

**Remark 1.2.1.** In a similar manner, for all $K \geqslant 1$, we denote by $\mathrm{MMD}_k(K)$ the space of marked minimal disks in $\tilde{X}$ whose boundary curves are $K$-quasicircles in $\hat{\mathbb{C}}$. Using the work [32] of Seppi, a similar analysis to that of Section 3 shows that, for all $K$ sufficiently close to $1$, and for all $h$ sufficiently close to a hyperbolic metric, $(\mathrm{MMD}_k(K), \partial_\infty)$ likewise carries the structure of a topological associated $\mathrm{PSL}(2, \mathbb{R})$-bundle over $\mathrm{QC}^+(K)$. We aim to study this structure further in forthcoming work.

By improving upon the laminated structure of [22], this foliated structure provides a pleasing framework within which the equidistribution theory developed by Labourie in [24] can be formulated. Indeed, here his laminar measures, which generalize the invariant measures introduced by Bonahon in [6] (c.f. also [10]), are simply those measures over the bundle $\mathrm{MKD}_{k,h}$ which are factored by the Lebesgue measure of the fibre in every trivialising chart.

Our second refinement yields a space of marked $k$-surfaces carrying a natural *foliated* structure. Note first that, for all $0 < k < 1$, $\mathrm{MKD}_{k,h}(1)$ consists of those $k$-disks in $X$ whose boundary curves are round circles in $\hat{\mathbb{C}}$. Let $SX$ denote the unit sphere bundle over $X$. For all $(D, p) \in \mathrm{MKD}_{k,h}(1)$, let $\hat{D} \subseteq SX$ denote the set of unit normal vectors over $D$. Note that, for all such $(D, p)$, $(\hat{D}, p)$ is a marked embedded surface in $SX$ which we call its *Gauss lift*.

**Theorem 1.2.2, Foliated Plateau problem**

*For all $0 < k < 1$, the set of Gauss lifts of elements of $\mathrm{MKD}_{k,h}(1)$ foliates $SX$.*

Theorem 1.2.2 plays a key role in our proof of Theorem 1.1.1. This is, in fact, the most significant difference between our work and that of Calegari–Marques–Neves, namely, that the argument used by Calegari–Marques–Neves to prove rigidity in [8] fails in the present case for the following reason. The proof of rigidity in [8] involves the study of the geometry of $X$ along a smoothly embedded disk $D$ obtained through a certain limiting process. In [8], $D$ is a minimal disk, and the area functional used ensures that it is totally geodesic, and that the sectional curvature of $X$ is equal to $-1$ over its tangent bundle $TD$. The totally geodesic property together with the geodesic conjugacy rigidity result proven by Gromov in [15] implies that $TD$ projects to a dense subset of $T(X/\Pi)$, from which it readily follows that the ambient metric is hyperbolic. In the present case, $D$ is a $k$-surface, and the area functional of (1.3) no longer allows us to conclude that it is totally geodesic, or even totally umbilic. Gromov's geodesic conjugacy rigidity result no longer applies and a different approach is required in order to prove density. Without Theorem 1.2.2, we find it hard to imagine how rigidity could be proven for the area functional used here.

Our proof of Theorem 1.1.1 is in fact closer in spirit to that proposed by Labourie in [24]. We emphasize, however, that Theorems 1.2.1 and 1.2.2 hold with no other restriction on the metric than that its sectional curvature be bounded above by $-1$. This contrasts markedly with the minimal surface case, where (see [26]) foliations analogous to that of Theorem 1.2.2 are only known to exist under restrictive conditions on the ambient metric, and are even known not to exist for certain negatively curved metrics. We thus view this remarkable property as exemplary of the rich structure of the space $\mathcal{S}_k(X/\Pi)$ of $k$-surfaces in $X/\Pi$. We discuss other interesting dynamical properties of this space in Sections 2.1, 3.1 and 3.2. Overall, we are led to concur with Labourie that, in the case of negatively-curved manifolds, the space of $k$-surfaces presents a convincing higher-dimensional generalization of the geodesic flow with an intriguing structure that remains to be fully understood.

**1.3 - Structure of the paper.** This paper is structured as follows.

**A -** In Section 2, we formulate and solve various forms of the Plateau problem for $k$-disks in Cartan–Hadamard manifolds. Our starting point is a solution of the asymptotic Plateau problem, proven in [22] in the case of ambient spaces of bounded geometry, and in [34] in the general case. We elevate this result to solutions of Plateau problems for foliations in $X$ and $SX$. Gaps in the foliation of $X$ are eliminated by uniqueness together with continuity properties of families of convex sets. The foliation property in $SX$ is obtained from the combinatorial properties of the intersection set $\Gamma$ of any two $k$-disks $D$ and $D'$, which, by the unique continuation theorem [2] of Aronszajn, is a graph with vertices of valency at least 4.

**B -** In Section 3, we develop the formalism of fibred Plateau problems. We show that the set of conformal parametrizations of $k$-disks in any cocompact Cartan-Hadamard manifold has the structure of a topological principal $\mathrm{PSL}(2,\mathbb{R})$-bundle over the space of quasicircles in $\hat{\mathbb{C}}$. We then show that the space of marked $k$-disks has the structure of a topological associated $\mathrm{PSL}(2,\mathbb{R})$-bundle. This formalism provides a succinct means of encoding the conformal structure of the leaves of the space of marked $k$-disks, which proves useful in Section 4 for describing measures constructed over this space.

**C -** In Section 4, using the formalism of associated $\mathrm{PSL}(2,\mathbb{R})$-bundles, we study the properties of a certain family of measures which are defined over the space $\mathrm{QC}^+$ of oriented quasicircles in $\hat{\mathbb{C}}$ using the surfaces constructed by Kahn–Marković in [18]. Following the ideas outlined by Labourie in [24], we prove an equidistribution property for this family of measures, whereby appropriately chosen sequences weakly converge to a measure of full support on the space $\mathrm{C}^+$ of oriented circles in $\hat{\mathbb{C}}$. We view this equidistribution result as characterizing in a synthetic manner the sequences constructed in [18].

**D -** In Section 5, by combining the results of the preceding sections with the arguments developed by Calegari–Marques–Neves in [8], we prove Theorem 1.1.1.

**E -** Finally, the notations introduced in this paper are summarized in Appendix A.

**1.4 - Acknowledgements.** Sébastien Alvarez was supported by CSIC, via Grupo I+D 159 "Geometría y acciones de grupos" and the program MIA. Ben Lowe was supported by NSF grant DMS-2202830. Graham Smith was supported by CNPq Bolsa de Produtividade em Pesquisa, Nível 1D. We are grateful to François Labourie and Thibault Lefeuvre for inspiring conversations without which this paper would not have attained its current form. We are also grateful to the anonymous reviewers for helpful comments made to the earlier draft of this paper. A large part of this work was carried out whilst Graham Smith was visiting the Institut des Hautes Études Scientifiques, and he is grateful for the excellent working conditions provided during that stay.

## 2 - Foliated Plateau problems.

**2.1 - Overview.** A key step in the proof of Theorem 1.1.1 is the solution of the foliated Plateau problem for $k$-disks. We first introduce the framework for the study of this problem. Let $X := (X, h)$ be an oriented 3-dimensional Cartan–Hadamard manifold of sectional curvature pinched between $-C$ and $-1$, for some $C \geqslant 1$. Let $\partial_\infty X$ denote its ideal boundary, and let $\overline{X} := X \cup \partial_\infty X$ denote its closure, furnished with the cone topology. Recall from the introduction that, for $k > 0$, a $k$-disk in $X$ is an oriented, smoothly embedded surface $D \subseteq X$, of constant extrinsic curvature equal to $k$, whose closure in $\overline{X}$ is homeomorphic to the closed unit disk in $\mathbb{R}^2$. We will adopt the convention which orients $D$ in such a manner that its compatible unit vector field $\nu$ points outwards from the convex set bounded by $D$. For all $k > 0$, let $\mathrm{KD}_{k,h}$ denote the space of $k$-disks, furnished with the $C^\infty_{\mathrm{loc}}$ topology. Let $\mathrm{JC}^+$ denote the space of oriented Jordan curves in $\partial_\infty X$, furnished with the Hausdorff topology. Given $c \in \mathrm{JC}^+$, we say that $D \in \mathrm{KD}_{k,h}$ spans $c$ whenever

$$\partial_\infty D = c, \tag{2.1}$$

where $\partial_\infty D$ is furnished with the orientation that it inherits from $D$.

**Theorem & Definition 2.1.1, Ideal Plateau problem**

*For all $0 < k < 1$, and for all $c \in \mathrm{JC}^+$, there exists a unique $k$-disk $D := \mathrm{D}_{k,h}(c) \in \mathrm{KD}_{k,h}$ which spans $c$. Furthermore, the function $\mathrm{D}_{k,h} : \mathrm{JC}^+ \to \mathrm{KD}_{k,h}$ is continuous.*

Theorem 2.1.1 is proven in Section 2.4.

This result is the starting point of our study of foliated Plateau problems. There are two cases that will interest us. First, let $\mathbb{S}^2$ denote the unit sphere in $\mathbb{R}^3$, and let $\mathrm{F}^+$ denote its singular foliation by latitudinal circles, oriented anticlockwise about the $z$-axis.

**Theorem 2.1.2, Foliated Plateau problem I**

*For all $0 < k < 1$, and for every homeomorphism $\alpha : \mathbb{S}^2 \to \partial_\infty X$, the family $(\mathrm{D}_{k,h}(\alpha(c)))_{c \in \mathrm{F}^+}$ foliates $X$.*

Theorem 2.1.2 is proven in Section 2.5.

Now let $SX$ denote the unit sphere bundle over $X$. Note that every $k$-disk $D$ bounds a convex subset of $X$. We will say that a unit normal vector over $D$ is *outward-pointing* whenever it points outwards from this set. We define the *Gauss lift* $\hat{D}$ of $D$ to be its set of outward-pointing, unit, normal vectors, and we note that $\hat{D}$ is a properly embedded disk in $SX$. Let $\mathrm{C}^+$ denote the space of oriented round circles in $\mathbb{S}^2$.

**Theorem & Definition 2.1.3, Foliated Plateau problem II**

*For all $0 < k < 1$, and for every homeomorphism $\alpha : \mathbb{S}^2 \to \partial_\infty X$, the family $(\hat{\mathrm{D}}_{k,h}(\alpha(c)))_{c \in \mathrm{C}^+}$ foliates $SX$. We call this family the $k$-disk foliation of $SX$ associated to $\alpha$.*

Theorem 2.1.3 is proven in Section 2.6.

Theorem 2.1.3 reveals the analogy between the space of oriented $k$-disks and the space of oriented geodesics. Indeed, recall that there is a natural homeomorphism identifying the space of ordered pairs of distinct points of $\partial_\infty X$ with the space of oriented geodesics in $X$ (see, for example, [15]). Upon viewing pairs of points as 0-dimensional circles, and intervals as 1-dimensional disks, we see that Theorem 2.1.3 merely constitutes a higher-dimensional analogue of this identification. Furthermore, just as this identification has become a staple of 2-dimensional hyperbolic dynamics, we will see presently that Theorem 2.1.3 and its consequences also provide useful tools for the study of certain dynamical properties of the space of $k$-disks.

This analogy is deepened by the fact that Theorems 2.1.2 and 2.1.3 also yield natural parametrizations of $SX$ which respect the geodesic foliation. Although this has no bearing on the proof of Theorem 1.1.1, we consider it to be of independent interest. In order to construct these foliations, we introduce marked round disks in $\mathbb{S}^2$, and we describe the geometry of their moduli space. We define a *marked round disk* in $\mathbb{S}^2$ to be a pair $(D, p)$, where $D \subseteq \mathbb{S}^2$ is a round disk and $p \in D$. We denote the space of marked round disks by MD, and we furnish this space with the Hausdorff topology in the first component, and the topology of $\mathbb{S}^2$ in the second. Note that every round disk in $\mathbb{S}^2$ inherits the natural orientation of the sphere. In particular, for all $(D, p) \in \mathrm{MD}$, we furnish the circle $\partial(D, p) := \partial D$ with the orientation that it inherits from $D$, and in this manner we see that $\partial$ defines a continuous map from MD into $\mathrm{C}^+$.

The space of marked round disks carries a natural smooth manifold structure. We now describe a natural smooth foliation by curves that it carries. Given $(D, p) \in \mathrm{MD}$, let $p'$ denote the image of $p$ under the conformal reflection $R_{\partial D}$ in $\partial D$, and let $\mathcal{D}$ denote the set of all round disks which contain $p$ and which are preserved by the group of elliptic Möbius maps which fix both $p$ and $p'$. The family $(D', p)_{D' \in \mathcal{D}}$ forms a smooth curve in MD passing through $(D, p)$, and the collection of all such curves forms a smooth foliation of MD which we denote by $\mathcal{F}_1$ (see Figure 2.1.1).

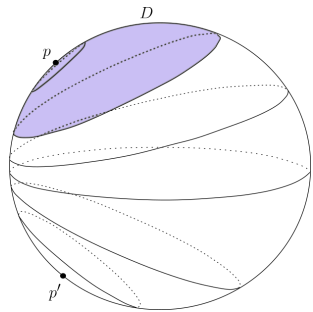

**Figure 2.1.1 - The canonical foliation by curves -** The point $p'$ is the image of $p$ under reflection in the boundary circle $\partial D$ of $D$. The family of disks preserved by the group of elliptic Möbius maps fixing these two points defines a smooth curve in the space of marked disks.

We are now ready to describe the first natural parametrization of $SX$ that Theorem 2.1.3 yields. Choose $0 < k < 1$ and a homeomorphism $\alpha : \mathbb{S}^2 \to \partial_\infty X$. We recall that the *horizon map* $\mathrm{hor} : SX \to \partial_\infty X$ is defined by

$$\mathrm{hor}(\xi) := \lim_{t \to +\infty} \gamma_\xi(t), \tag{2.2}$$

where $\gamma_\xi : \mathbb{R} \to X$ is the unique geodesic such that $\dot\gamma_\xi(0) = \xi$. For all $(D, p) \in \mathrm{MD}$, we define $\Phi_{k,h,\alpha}(D, p)$ to be the unique, unit, transverse vector $\xi$ to $\mathrm{D}_{k,h}(\alpha(\partial D))$ such that

$$\mathrm{hor}(\xi) = \alpha(p) \text{ and } \mathrm{hor}(-\xi) = (\alpha \circ R_{\partial D})(p). \tag{2.3}$$

**Theorem 2.1.4**

*For all $0 < k < 1$, and for every homeomorphism $\alpha : \mathbb{S}^2 \to \partial_\infty X$, $\Phi_{k,h,\alpha}$ defines a homeomorphism from MD into $SX$ which sends $\mathcal{F}_1$ to the geodesic foliation of $SX$.*

**Remark 2.1.1.** In the case of interest to us, where $X$ is acted on cocompactly by some hyperbolic group $\Pi$, a natural choice of homeomorphism $\alpha$ is given by (1.1), so that the homeomorphism $\Phi_{k,h,\alpha}$ really only depends on $0 < k < 1$. For this reason, we consider it to be more natural than the homeomorphism described by Gromov under the heading "Geodesic Rigidity" in [15].

Theorem 2.1.4 is proven in Section 2.6.

Interestingly, Theorem 2.1.3 also yields natural parametrizations of $SX$ which respect the $k$-surface foliation. To see this, we note first that MD carries a natural smooth foliation by surfaces. Indeed, given $(D, p) \in \mathrm{MD}$, the family $(D, q)_{q \in D}$ forms a smooth surface in MD passing through $(D, p)$, and the collection of all such surfaces yields the desired smooth foliation. We denote this foliation by $\mathcal{F}_2$, and we observe in passing that it is trivially transverse to the foliation $\mathcal{F}_1$ given above. For all $(D, p) \in \mathrm{MD}$, we now define $\Psi_{k,h,\alpha}(D, p)$ to be the unique, unit, normal vector $\xi$ to $\mathrm{D}_{k,h}(\alpha(\partial D))$ such that

$$\mathrm{hor}(\xi) = p. \tag{2.4}$$

**Theorem 2.1.5**

*For all $0 < k < 1$, and for every homeomorphism $\alpha : \mathbb{S}^2 \to \partial_\infty X$, $\Psi_{k,h,\alpha}$ defines a homeomorphism from MD into $SX$ which sends $\mathcal{F}_2$ into the $k$-surface foliation of $SX$ associated to $\alpha$.*

Theorem 2.1.5 is proven in Section 2.6.

It is natural to study how these two parametrizations are related to one-another, and how the properties of the surface foliation are related to those of the ambient manifold. A few moments' reflection leads us naturally to the following questions.

*Question 1:* Given $0 < k < 1$ and a homeomorphism $\alpha : \mathbb{S}^2 \to \partial_\infty X$, is the homeomorphism $\Phi_{k,h,\alpha}$ in some way related to the one constructed by Gromov in [15]?

*Question 2:* In the spirit of known rigidity results, given $0 < k < 1$ and a homeomorphism $\alpha : \mathbb{S}^2 \to \partial_\infty X$, if $\Phi_{k,h,\alpha}$ and $\Psi_{k,h,\alpha}$ coincide, then is it true that $X$ has constant curvature and $\alpha$ is conformal?

*Question 3:* In the spirit of [4], given $0 < k < 1$ and a homeomorphism $\alpha : \mathbb{S}^2 \to \partial_\infty X$, if the $k$-disk foliation is smooth, then is it true that $X$ has constant curvature and $\alpha$ is conformal?

*Question 4:* For all $0 < k < 1$, it is possible to construct, non-canonically, a homeomorphism $\Theta : \text{MD} \to SX$ which restricts to a smooth diffeomorphism over every leaf. Can such a homeomorphism be constructed in a canonical manner?

**2.2 - Convex subsets.** We begin by developing a formalism for the study of convex embedded surfaces in $X$. It will be convenient to work with the convex sets that such surfaces bound. Recall that a subset of $\overline{X}$ is said to be *convex* whenever it contains the unique geodesic segment joining any two of its points. Let $\text{CC} := \text{CC}(X)$ denote the space of compact, convex subsets of $\overline{X}$, furnished with the Hausdorff topology, and recall that this space is compact. For $K \in \text{CC}$, we denote

$$\begin{aligned}\partial_\infty K &:= K \cap \partial_\infty X, \\ \partial_{\text{fin}} K &:= (\partial K) \cap X, \text{ and} \\ \hat{\partial} K &:= \partial_\infty K \cup \partial_{\text{fin}} K.\end{aligned} \quad (2.5)$$

These definitions are illustrated in Figure 2.2.2.

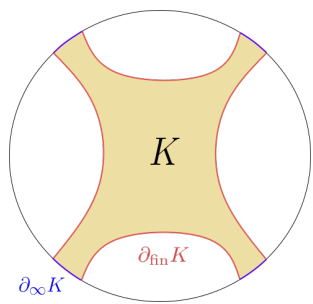

Figure 2.2.2 - Boundary components of convex sets - $K$ is a compact convex subset of the closure $\overline{X}$ of some Cartan–Hadamard manifold $X$. We denote by $\partial_\infty K$ the intersection of $K$ with the ideal boundary $\partial_\infty X$ of $X$, and we denote by $\partial_{\text{fin}} K$ the topological boundary of $K$ in $X$. The topological boundary $\hat{\partial} K$ of $K$ in $\overline{X}$ is the union of these two components.

**Lemma 2.2.1**

*If $K \subseteq \text{CC}$ has non-trivial interior, then $\hat{\partial} K$ is homeomorphic to a 2-dimensional sphere.*

**Proof:** Let $x$ be an interior point of $K$, and let $S$ denote the unit sphere in $T_x X$. By convexity, every geodesic ray leaving $x$ meets $\hat{\partial} K$ at a unique, possibly ideal, point. Since the geodesic rays leaving $x$ are determined by their initial directions, this yields a bijection $\phi : S \to \hat{\partial} K$, which is trivially continuous. Since $S$ is compact, $\phi$ is a homeomorphism, as desired. $\square$

**Lemma 2.2.2**

*If $(K_m)_{m\in\mathbb{N}}$ is a sequence in CC converging to $K_\infty \in$ CC, then $(\hat{\partial}K_m)_{m\in\mathbb{N}}$ Hausdorff converges to $\hat{\partial}K_\infty$.*

**Remark 2.2.1.** In other words $\hat{\partial}$ defines a continuous map from CC into the space of compact subsets of $\overline{X}$ furnished with the Hausdorff topology.

**Proof:** Let $(K_m)_{m\in\mathbb{N}}$ be a sequence in CC converging to $K_\infty \in$ CC. We first show that every Hausdorff accumulation point of $(\hat{\partial}K_m)_{m\in\mathbb{N}}$ is contained in $\hat{\partial}K_\infty$. It suffices to show that if $(x_m)_{m\in\mathbb{N}} \in \overline{X}$ is a sequence converging to $x_\infty \in \overline{X}$ such that $x_m \in \hat{\partial}K_m$ for all $m$, then $x_\infty \in \hat{\partial}K_\infty$. If $x_\infty \in \partial_\infty X$, then $x_\infty \in K_\infty \cap \partial_\infty X \subseteq \hat{\partial}K$, as asserted. Suppose now that $x_\infty \in X$. We may suppose that $x_m \in \partial_{\text{fin}} K_m$ for all $m$. For all $m \in \mathbb{N}\cup\{\infty\}$, let $d_m : X \to [0,\infty[$ denote the distance function in $X$ to $K_m$, and note that $(d_m)_{m\in\mathbb{N}}$ converges locally uniformly to $d_\infty$. For all $m \in \mathbb{N}\cup\{\infty\}$, let $\gamma_m : [0,\infty[\to X$ be a geodesic ray normal to $K_m$ at $x_m$. We may suppose that $(\gamma_m)_{m\in\mathbb{N}}$ converges to the geodesic ray $\gamma_\infty : [0,\infty[\to X$, say. However, for all finite $m$, and for all $t$, $(d_m \circ \gamma_m)(t) = t$. Taking limits therefore yields, for all $t$, $(d_\infty \circ \gamma_\infty)(t) = t$, so that $x_\infty = \gamma_\infty(0) \in \hat{\partial}X$, as asserted.

The reader may verify that, for any sequence $(A_m)_{m\in\mathbb{N}}$ of compact subsets of $\overline{X}$ Hausdorff converging to $A_\infty$, say, every Hausdorff accumulation point of $(\hat{\partial}A_m)_{m\in\mathbb{N}}$ contains $\hat{\partial}A_\infty$. It follows that $(\hat{\partial}K_m)_{m\in\mathbb{N}}$ only accumulates on $\hat{\partial}K_\infty$. Since the set of closed subsets of $\overline{X}$ is compact, $(\hat{\partial}K_m)_{m\in\mathbb{N}}$ Hausdorff converges towards $\hat{\partial}K_\infty$, as desired. $\square$

For every positive real number $k$, and for every open subset $\Omega$ of $X$, let $\text{CC}_{k,h}(\Omega)$ denote the subspace of CC consisting of those compact, convex subsets $K \subseteq \overline{X}$, with non-trivial interior, such that $\hat{\partial}K \cap \Omega$ is smooth and of constant extrinsic curvature equal to $k$ with respect to $h$. For all $k > 0$, we denote $\text{CC}_{k,h} := \text{CC}_{k,h}(X)$.

**Lemma 2.2.3**

*Choose $k > 0$, let $\Omega$ be an open subset of $X$, let $Y$ be a compact subset of $\Omega$, and let $\mathcal{Y}$ be a subset of $\text{CC}_{k,h}(\Omega)$ consisting of elements whose boundaries intersect $Y$ non-trivially. Then, either*

*(1) there exists $B > 0$ such that, for all $K \in \mathcal{Y}$, the norm of the shape operator of $\partial K \cap \Omega$ is bounded above by $B$ at every point of $\partial K \cap Y$; or*

*(2) the closure $\overline{\mathcal{Y}}$ of $\mathcal{Y}$ in CC contains a closed geodesic segment in $\overline{X}$ with extremities in $\overline{X} \setminus \Omega$.*

**Proof:** Let $(K_m)_{m\in\mathbb{N}}$ be a sequence in $\mathcal{Y}$ and, for all $m$, let $x_m$ be a point of $\partial K_m \cap Y$. Suppose that $(K_m)_{m\in\mathbb{N}}$ and $(x_m)_{m\in\mathbb{N}}$ converge respectively to $K_\infty$ and $x_\infty$, say. Let $\varepsilon > 0$ be such that, for all $m \in \mathbb{N}\cup\{\infty\}$, $B_\varepsilon(x_m) \subseteq \Omega$. For all finite $m$, let $\hat{S}_m \subseteq SX$ denote the set of supporting normals of $\partial K_m$ over $B_\varepsilon(x_m)$ and suppose, for convenience, that this surface consists of a single connected component. By Theorem 6.5 of [21] (see also [35]), $(\hat{S}_m)_{m\in\mathbb{N}}$ subconverges in the $C^\infty_{\text{loc}}$ sense to a smooth, embedded surface $\hat{S}_\infty$, say. In particular, $\hat{S}_\infty$ is the set of supporting normals to $K_\infty$ over $B_\varepsilon(x_\infty)$.

Now suppose that (1) does not hold. For all finite $m$, let $B_m$ denote the norm of the shape operator of $\partial K_m$ at $x_m$. We may suppose that $(B_m)_{m\in\mathbb{N}}$ tends to $+\infty$. It then follows by Theorem 6.5 of [21] (see also [35]) that $\hat{S}_\infty$ is the unit normal bundle over some open geodesic segment $\Gamma$, say, containing $x_\infty$. $K_\infty$ is thus a closed geodesic segment containing $x_\infty$ in its relative interior. It remains to show that the extremities of $K_\infty$ lie in $\overline{X}\setminus\Omega$. However, were this not the case, upon modifying $Y$ if necessary, we could choose $x_\infty$ to be an extremity of $K_\infty$, which is absurd. This completes the proof. $\square$

**Lemma 2.2.4**

*Choose $k > 0$, let $\Omega$ be an open subset of $X$, and let $(K_m)_{m\in\mathbb{N}} \in \text{CC}_{k,h}(\Omega)$ be a sequence converging towards $K_\infty \in$ CC. If $\hat{\partial}K_\infty \cap \Omega$ is non-empty, then either*

*(1) $K_\infty$ has non-trivial interior, and $(\hat{\partial}K_m \cap \Omega)_{m\in\mathbb{N}}$ converges towards $\hat{\partial}K_\infty \cap \Omega$ in the $C^\infty_{\text{loc}}$ sense; or*

*(2) $K_\infty$ is a closed geodesic segment in $\overline{X}$ with extremities in $\overline{X} \setminus \Omega$.*

**Remark 2.2.2.** In particular, the map $K \mapsto \hat{\partial}K \cap \Omega$ defines a continuous function from $\text{CC}_{k,h}(\Omega)$ into the space of properly embedded smooth hypersurfaces in $\Omega$, furnished with the $C^\infty_{\text{loc}}$ topology.

**Corollary 2.2.5**

Let $\partial\mathrm{CC}_{k,h}(\Omega)$ denote the topological boundary of $\mathrm{CC}_{k,h}(\Omega)$ in $\mathrm{CC}$. For all $K \in \partial\mathrm{CC}_{k,h}(\Omega)$, either

(1) $\hat{\partial}K \cap \Omega = \emptyset$; or

(2) $K$ is a closed geodesic segment in $\overline{X}$ with extremities in $\overline{X} \setminus \Omega$.

The case where $\Omega = X$ is worth stating separately.

**Corollary 2.2.6**

For every positive real number $k$, $\partial\mathrm{CC}_{k,h}$ consists of $X$, points of $\partial_\infty X$ and complete geodesics in $\overline{X}$.

**Proof of Lemma 2.2.4:** Suppose that $\hat{\partial}K_\infty \cap \Omega$ is non-empty and that $K_\infty$ is not a closed geodesic segment in $\overline{X}$ with extremities in $\overline{X} \setminus \Omega$. For all $m$, let $A_m$ denote the shape operator of $\hat{\partial}K_m \cap \Omega$. By Lemma 2.2.3, for every compact subset $Y$ of $\Omega$, there exists $B > 0$ such that, for all $m$, and for all $x \in \hat{\partial}K_m \cap \Omega$,

$$\|A_m(x)\| \leqslant B.$$

By classical elliptic theory, for every compact subset $Y$ of $\Omega$, and for all $k \in \mathbb{N}$, there exists $B_k > 0$ such that, for all $m$, and for all $x \in \hat{\partial}K_m \cap Y$,

$$\|\nabla^k A_m(x)\| \leqslant B_k,$$

where here $\nabla$ denotes the intrinsic Levi-Civita covariant derivative of $\hat{\partial}K_m \cap \Omega$. Since $(\hat{\partial}K_m)_{m \in \mathbb{N}}$ Hausdorff converges towards $\hat{\partial}K_\infty$, it follows that $(\hat{\partial}K_m \cap \Omega)_{m \in \mathbb{N}}$ converges in the $C^\infty_{\mathrm{loc}}$ sense towards $\hat{\partial}K_\infty \cap \Omega$, as desired. $\square$

**Proof of Corollary 2.2.5:** Indeed, by Lemma 2.2.4, if neither (1) nor (2) holds, then $K_\infty$ has non-trivial interior and $\hat{\partial}K_\infty \cap \Omega$ is smooth and of constant extrinsic curvature equal to $k$. That is $K \in \mathrm{CC}_{k,h}(\Omega)$, which is absurd, and the result follows. $\square$

**2.3 - Geometric maximum principles.** We will require the following version of the geometric maximum principle for submanifolds of Cartan–Hadamard manifolds (see Figure 2.3.3).

**Theorem 2.3.1, Global geometric maximum principle**

Let $X := (X, h)$ be a 3-dimensional Cartan–Hadamard manifold of sectional curvature bounded above by $-1$. For all $0 < k \leqslant k' < 1$, and for every open subset $\Omega$ of $X$, if $K \in \mathrm{CC}_{k,h}(\Omega)$ and $K' \in \mathrm{CC}_{k',h}(\Omega)$ satisfy

$$K \setminus \Omega \subseteq \mathrm{Int}(K') \setminus \Omega, \tag{2.6}$$

then

$$K \subseteq K'. \tag{2.7}$$

**Remark 2.3.1.** We will see presently in Corollary 2.3.3 that, for $0 < k, k' < 1$, neither $K$ nor $K'$ can be wholly contained in $\Omega$, so that Condition (2.6) is always non-trivial.

Our starting point will be the following more classical result, the proof of which we include for the reader's convenience.

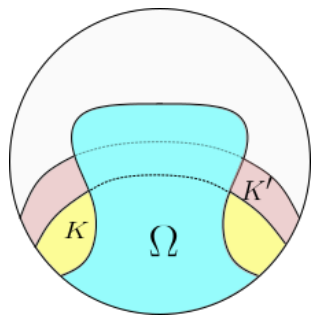

Figure 2.3.3 - **The geometric maximum principle** - If $K$ is contained within the interior of $K'$ outside $\Omega$, then $K$ is also contained within $K'$ inside $\Omega$.

**Lemma 2.3.2, Local geometric maximum principle**

Let $X := (X, h)$ be a riemannian manifold. Let $K$ and $K'$ be closed subsets of $X$ with non-trivial interiors such that $K \subseteq K'$. If $x \in \partial K \cap \partial K'$, if $\partial K$ and $\partial K'$ are both smooth near this point, and if $\mathrm{II}_x$ and $\mathrm{II}'_x$ denote their respective shape operators with respect to their outward pointing unit normals at this point, then

$$\mathrm{II}_x \geqslant \mathrm{II}'_x, \qquad (2.8)$$

in the sense that $\mathrm{II}_x - \mathrm{II}'_x$ is non-negative semi-definite.

**Corollary 2.3.3**

Let $X := (X, h)$ be a Cartan–Hadamard manifold of sectional curvature bounded above by $-1$, and let $\Omega \subseteq X$ be an open subset. If $0 < k < 1$, and if $K \in \mathrm{CC}_{k,h}(\Omega)$, then

$$K \setminus \Omega \neq \emptyset. \qquad (2.9)$$

**Proof:** Indeed, suppose the contrary, so that $K \subseteq \Omega$. Choose $x \in X$. Since $K$ is compact, and since $K \subseteq \Omega \subseteq X$, the distance in $X$ to this point attains its maximum over $K$ at some point $y$, say. Since $y \in \Omega$, $\partial K$ is smooth at this point with extrinsic curvature equal to $k < 1$. Now let $K'$ denote the closed ball of radius $r := d(x, y)$ about $x$. Since $X$ has sectional curvature bounded above by $-1$, $\partial K'$ is smooth and locally strictly convex with extrinsic curvature at every point strictly greater than 1. This is absurd by Lemma 2.3.2, and it follows that $K$ is not wholly contained in $\Omega$, as desired. $\square$

Lemma 2.3.2 is a consequence of the following simple, but very useful, relation, for which we know of no reference in the literature.

**Lemma 2.3.4**

Let $X := (X, h)$ be a riemannian manifold, let $f : X \to \mathbb{R}$ be a smooth function, and let $Y \subseteq X$ be a smooth hypersurface. The hessians of $f$ and $f|_Y$ are related by

$$\mathrm{Hess}(f|_Y) = \mathrm{Hess}(f)|_{TY} - df(\nu)\mathrm{II}, \qquad (2.10)$$

where $\nu$ denotes the unit normal vector field over $Y$, and $\mathrm{II}$ denotes its shape operator.

**Remark 2.3.2.** Note that reversing the sign of $\nu$ also reverses the sign of $\mathrm{II}$, so that (2.10) is independent of any choice of orientation of $X$ and $Y$.

**Proof:** Let $\nabla$ and $\overline{\nabla}$ denote the respective Levi-Civita covariant derivatives of $Y$ and $X$. For all vector fields $\xi$ and $\mu$ tangent to $Y$,

$$\mathrm{Hess}(f|_Y)(\xi, \mu) = D_\xi D_\mu f - df(\nabla_\xi \mu) = D_\xi D_\mu f - df(\overline{\nabla}_\xi \mu) - df(\nu)\mathrm{II}(\xi, \mu) = \mathrm{Hess}(f)(\xi, \mu) - df(\nu)\mathrm{II}(\xi, \mu),$$

as desired. $\square$

**Proof of Lemma 2.3.2:** Let $d, d' : X \to \mathbb{R}$ denote respectively the signed distances in $X$ to $\partial K$ and $\partial K'$, taken to have positive signs over the respective complements of $K$ and $K'$. Upon restricting to a neighbourhood $U$ of $x$, we may assume that these functions are both smooth and that $d \geqslant d'$. Since $d(x) = d'(x) = 0$,

$$\nabla d(x) = \nabla d'(x), \text{ and}$$
$$\mathrm{Hess}(d)(x) \geqslant \mathrm{Hess}(d')(x).$$

Note that $\nabla d(x)$ and $\nabla d'(x)$ are the respective outward-pointing unit normal vectors of $\partial K$ and $\partial K'$ at $x$. In particular,

$$T_x \partial K = T_x \partial K'.$$

Since $d$ and $d'$ are constant over $\partial K$ and $\partial K'$ respectively, (2.10) yields

$$\mathrm{II}_x = \mathrm{Hess}(d)(x)|_{T_x \partial K} \geqslant \mathrm{Hess}(d')(x)|_{T_x \partial K'} = \mathrm{II}'_x,$$

as desired. □

We now return to the case where $X$ is a 3-dimensional Cartan–Hadamard manifold of sectional curvature bounded above by $-1$. Theorem 2.3.1 will follow from the study of the curvatures of surfaces which are equidistant to certain other geometric objects. We first address surfaces equidistant to other surfaces. Let $S$ be a convex, embedded surface in $X$ and let $\nu$ denote its outward-pointing, unit, normal vector field. For all $t \geqslant 0$, denote

$$e_t : S \to X; x \mapsto \mathrm{Exp}(t\nu(x)), \tag{2.11}$$

where Exp denotes the exponential map of $X$. We call $S_t := (S, e_t)$ the *equidistant surface* of $S$ at distance $t$. Note that $S_t$ is embedded for all $t$.

**Lemma 2.3.5**

Let $X := (X, h)$ be a 3-dimensional Cartan–Hadamard manifold of sectional curvature bounded above by $-1$. Choose $0 < k_0 < 1$ and let $S$ be a convex, embedded surface in $X$ of extrinsic curvature everywhere at least $k_0$. For all $t > 0$, the extrinsic curvature of $S_t$ is everywhere at least

$$k_t := \tanh^2(t_0 + t), \tag{2.12}$$

where

$$\tanh^2(t_0) := k_0. \tag{2.13}$$

In particular, for all $t > 0$, the extrinsic curvature of $S_t$ is everywhere strictly greater than $k_0$.

**Proof:** For all $t$, let $A_t$ denote the shape operator of $e_t$ and let $K_t := \mathrm{Det}(A_t)$ denote its extrinsic curvature. Recall (see, for example, the Tube Formula (∗∗) of [14]) that $A_t$ satisfies the ordinary differential equation

$$\partial_t A_t = W_t - A_t^2,$$

where $W_t \in \mathrm{End}(TS)$ is defined by

$$\langle W_t \xi, \xi \rangle := \langle R(\nu_t, \xi)\nu_t, \xi \rangle,$$

$R$ denotes the Riemann curvature tensor of $X$, and $\nu_t$ denotes the unit normal vector field over $e_t$. Since $X$ has sectional curvature bounded above by $-1$, this yields

$$\partial_t A_t - \mathrm{Id} + A_t^2 \geqslant 0,$$

in the sense that the matrix-valued function on the left-hand side is everywhere non-negative semi-definite. As long as $A_t$ is positive-definite, taking determinants yields

$$\partial_t K_t = K_t \mathrm{Tr}(A_t^{-1} \partial_t A_t) \geqslant K_t \mathrm{Tr}(A_t^{-1} - A_t) \geqslant \mathrm{Tr}(A_t)(1 - K_t) \geqslant 2\sqrt{K_t}(1 - K_t).$$

The result follows upon solving this ordinary differential inequality. □

We now address surfaces in $X$ equidistant to curves. Let $\Gamma$ be a smoothly embedded curve in $X$, let $T : \Gamma \to SX$ denote its unit tangent vector field, let $\kappa : \Gamma \to \mathbb{R}$ denote its geodesic curvature, and let $N_1\Gamma$ denote its unit normal bundle. We say that $\Gamma$ is *regular* whenever $\kappa$ never vanishes. In this case, its unit *normal* and *binormal* vector fields are respectively defined by

$$N := \frac{1}{|\kappa|}\nabla_T T, \text{ and} \tag{2.14}$$
$$B := T \wedge N.$$

We say that a vector $\nu_x \in N_1\Gamma$ is a *convex direction* of $\Gamma$ whenever

$$\langle \nu_x, N(x) \rangle < 0. \tag{2.15}$$

We denote the set of convex directions of $\Gamma$ by $C_1\Gamma$. Note that, for any unit vector $\nu_x \in N_1\Gamma$, there trivially exists an embedded surface $S$, normal to $\nu_x$, containing an open segment of $\Gamma$ about $x$. The convex directions of $\Gamma$ at $x$ are precisely those unit normal vectors $\nu_x$ for which these surfaces can be chosen to be convex. For all $t$, we denote

$$e_t : C_1\Gamma \to X; x \mapsto \mathrm{Exp}(tx), \tag{2.16}$$

where Exp again denotes the exponential map of $X$. We call $\Gamma_t := (C_1\Gamma, e_t)$ the *convex equidistant surface* of $\Gamma$ at distance $t$.

**Lemma 2.3.6**

Let $X := (X, h)$ be a 3-dimensional Cartan–Hadamard manifold of sectional curvature bounded above by $-1$. Let $\Gamma$ be a smoothly embedded regular curve in $X$. For all $t > 0$, the extrinsic curvature of $\Gamma_t$ is everywhere at least $1$.

**Proof:** Choose $\nu_x \in C_1\Gamma$. Let $S$ be an embedded surface in $X$ which contains $\Gamma$ and is normal to $\nu_x$ at $x$. Let $d_\Gamma$ denote distance to $\Gamma$ in $X$, and let $d_S$ denote signed distance to $S$ in $X$, with sign chosen such that $\nabla d_S(x) = \nu_x$. For $A > 0$, consider the function
$$f_A(y) := d_S(y) + A d_\Gamma(y)^2.$$
We claim that, for sufficiently large $A$, the restriction of $\mathrm{Hess}(f_A(x))$ to the plane orthogonal to $\nu_x$ has arbitrarily large, positive determinant. Indeed, note first that $\nabla d_\Gamma^2$ vanishes along $\Gamma$ so that, for any two vector fields $\xi_1$ and $\xi_2$, with $\xi_1$ tangent to $\Gamma$,
$$\mathrm{Hess}(d_\Gamma^2)(\xi_1, \xi_2)(x) = \left(D_{\xi_1}\langle \nabla d_\Gamma^2, \xi_2\rangle\right)(x) - \left(\langle \nabla d_\Gamma^2, \nabla_{\xi_1}\xi_2\rangle\right)(x) = 0.$$
Next if $\xi$ is normal to $\Gamma$, then
$$\mathrm{Hess}(d_\Gamma^2)(\xi, \xi)(x) = \frac{\partial^2}{\partial t^2} d_\Gamma^2(\mathrm{Exp}(t\xi(x)))\bigg|_{t=0} = 2.$$
Now let $N_S$ and $\mathrm{II}_S$ denote respectively the unit normal vector field and second fundamental form of $S$, and let $\xi$ be a unit-length vector field tangent to $\Gamma$. By Lemma 2.3.4,
$$\mathrm{Hess}(d_S)(\xi, \xi)(x) = \mathrm{II}_S(\xi, \xi)(x) = \langle \nabla_\xi N_S, \xi\rangle(x) = -\langle (\nabla_\xi \xi)(x), \nu_x\rangle.$$
Since $\nu_x$ is a convex direction of $\Gamma$,
$$\delta := \mathrm{Hess}(d_S)(\xi, \xi)(x) > 0.$$
Upon combining these relations, we see that
$$\mathrm{Det}\big(\mathrm{Hess}(f_A)(x)|_{T_x S}\big) = 2A\delta + \mathrm{Det}\big(\mathrm{Hess}(d_S)(x)|_{T_x S}\big).$$
Upon choosing $A$ sufficiently large and positive, this quantity can also be made arbitrarily large and positive, and the assertion follows.

For $A > 0$, and for $\delta > 0$ sufficiently small, denote
$$S'_{A,\delta} := \{y \in X \mid d(y,x) < \delta \,,\, f_{A,\delta}(y) = 0\},$$
and, for all $t > 0$, let $S'_{A,\delta,t}$ denote the equidistant surface of $S'_{A,\delta}$ at distance $t$. By Lemma 2.3.4 and the above discussion, upon choosing $A$ sufficiently large, we may suppose that $S'_{A,\delta}$ has extrinsic curvature at least $1$ at every point. By Lemma 2.3.5, we may suppose that the same holds for $S'_{A,\delta,t}$ for all $t$. Since $\Gamma_t$ is an interior tangent to $S_t$ at $\mathrm{Exp}(t\nu_x)$, it follows by Lemma 2.3.2 that $\Gamma_t$ also has extrinsic curvature at least $1$ at this point. Since $\nu_x \in C_1\Gamma$ is arbitrary, the result follows. $\square$

We are now ready to prove Theorem 2.3.1.

**Proof of Theorem 2.3.1:** Since $\partial_\infty K$ is contained in the interior of $K'$, we may remove the closure of a neighbourhood of $\partial_\infty X$ from $\Omega$ without changing the hypotheses of the theorem. We may therefore suppose that $\Omega$ is bounded. Denote $K'' := K \cap K'$, and note that $K'' \setminus \Omega = K \setminus \Omega$. It trivially suffices to show that $K$ is contained in $K''$. However, suppose the contrary. By hypothesis, $K \setminus K''$ is contained in $\Omega$. Since $K \cap \overline{\Omega}$ is compact, it contains a point $x$, say, maximising distance to $K''$. Trivially, $x$ is an element of $\Omega$. Denote $r := d(x, K'')$. Then $K$ is an interior tangent to $\partial B_r(K'')$ at this point. This is illustrated in Figure 2.3.4.

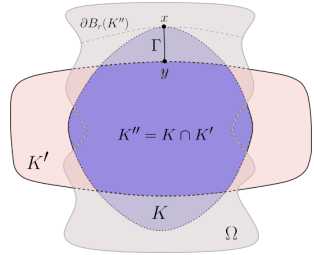

**Figure 2.3.4 - Proof of the geometric maximum principle I -** We suppose the contrary, so that some portion of $K$ lies outside $K'$ in $\Omega$. We denote $K'' := K \cap K'$ and we choose a point $x \in K$ maximising distance to $K''$.

Following [7], we say that a $C^1$ embedded surface $S$ has extrinsic curvature greater than $k > 0$ in the weak sense at some point $x$ whenever there exists a smooth, convex surface $S'$, which is an exterior tangent to $S$ at $x$, whose extrinsic curvature at this point is greater than $k$. Note that, by Lemma 2.3.2, when $S$ is smooth, this property implies that $S$ has extrinsic curvature greater than $k$ in the classical sense at this point.

We claim that $\partial B_r(K'')$ has extrinsic curvature strictly greater than $k$ in the weak sense at $x$. Indeed, let $y \in \partial K''$ denote the closest point in $K''$ to $x$, let $\Gamma$ denote the geodesic segment from $x$ to $y$, and note that $\Gamma$ is an external normal to $\partial K''$ at $y$. There are three cases to consider. If $y \in \partial K \setminus K'$ then, by convexity, $\Gamma$ can only meet $K$ at $y$, which is absurd. If $y \in \partial K' \setminus K$ then, in particular, $y \in \Omega$ so that $\partial K'$ is smooth and of constant extrinsic curvature equal to $k'$ at this point. By Lemma 2.3.5, near $x$, $\partial B_r(K'') = \partial B_r(K')$ has extrinsic curvature strictly greater than $k$, as asserted.

We now address the case where $y \in \partial K \cap \partial K'$. As before, $y \in \Omega$, so that $\partial K$ and $\partial K'$ are smooth with unique normals at this point. We denote these normals respectively by $\nu_y$ and $\nu_y'$. Note that $\nu_y$ and $\nu_y'$ cannot be antipodal, for otherwise, by strict convexity, $K''$ would consist of a single point, which is absurd. Let $S_y$ denote the unit sphere in $T_y X$ and note (see, for example, Theorem 5.20 of [36]) that the set of supporting normals to $K''$ at $y$ is the shorter geodesic segment $\Delta_y$ in $S_y$ joining $\nu_y$ and $\nu_y'$. Let $\nu_y'' \in \Delta_y$ denote the tangent to $\Gamma$ at $y$.

There are now three cases to consider. As before, $\nu_y'' \neq \nu_y$, for otherwise $\Gamma$ would meet $K$ only at $y$, which is absurd. If $\nu_y'' = \nu_y'$, then, near $x$, $\partial B_r(K'')$ is an interior tangent to $\partial B_r(K')$ at $x$ and thus, as before, has extrinsic curvature strictly greater than $k$ in the weak sense at this point, as asserted. Finally (see Figure 2.3.5), if $\nu_y''$ lies in the relative interior of $\Delta_y$ then, in particular, $\Delta_y$ is non-trivial, so that $\partial K$ and $\partial K'$ meet transversally at $y$. Their intersection is then a smooth curve $c$, say, of which $\nu_y''$ is a convex direction. Since $\partial B_r(K'')$ coincides with $\partial B_r(c)$ near $x$, it follows by Lemma 2.3.6 that $\partial B_r(K'')$ has extrinsic curvature strictly greater than $k$ at this point, as asserted.

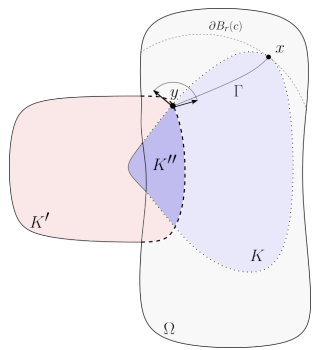

**Figure 2.3.5 - Proof of the geometric maximum principle II -** It is possible for $\Gamma$ to meet $K''$ along the intersection of $\partial K$ and $\partial K'$.

We conclude that, in all cases, $\partial B_r(K'')$ has extrinsic curvature strictly greater than $k$ in the weak sense at $x$. Since $\partial K$ is an interior tangent to $\partial B_r(K'')$ at $x$, it also has extrinsic curvature strictly greater than $k$ at this point. This is absurd, and it follows that $K$ is contained in $K''$, as desired. $\square$

**2.4 - Ideal Plateau problems.** We now recall the main results of [22] and [34] as they pertain to the present paper. Given an oriented Jordan curve $c$ in an oriented topological sphere $\mathbb{S}^2$, we denote its interior and exterior in $\mathbb{S}^2$ respectively by $\text{Int}(c)$ and $\text{Ext}(c)$, and we denote the closures of these sets respectively by $\overline{\text{Int}}(c)$ and $\overline{\text{Ext}}(c)$. In the present framework, the main theorem of [34] immediately yields the following result.

**Theorem & Definition 2.4.1**

*Let $X := (X, h)$ be an oriented, 3-dimensional Cartan–Hadamard manifold of sectional curvature pinched between $-C$ and $-1$, for some $C \geqslant 1$. For all $0 < k < 1$, and for every oriented Jordan curve $c \in \text{JC}^+$, there exists a unique element $K := K_{k,h}(c) \in \text{CC}_{k,h}$ such that*

$$\partial_\infty K = \overline{\text{Ext}}(c). \tag{2.17}$$

*We denote*
$$\mathrm{D}_{k,h}(c) := \partial_{\mathrm{fin}} \mathrm{K}_{k,h}(c). \tag{2.18}$$

*In particular, $\mathrm{D}_{k,h}(c)$ is the unique k-disk which spans c.*

**Remark 2.4.1.** We underline that, in the case that will concern us in the sequel, namely where $X$ is acted upon cocompactly by some group $\Pi$, this result also follows from Theorems $A$ and $E$ of [22].

We also require the following non-existence result, proven in [34] in the general case, and proven in Theorem $B$ of [22] in the cocompact case.

**Theorem 2.4.2**

*Let $X := (X, h)$ be an oriented, 3-dimensional Cartan–Hadamard manifold of sectional curvature pinched between $-C$ and $-1$, for some $C \geqslant 1$. For all $0 < k < 1$, there exists no element $K \in \mathrm{CC}_{k,h}$ such that $\partial_\infty K$ consists of a single point.*

In this context, the geometric maximum principle of Theorem 2.3.1 yields the following result.

**Lemma 2.4.3**

*For all $0 < k \leqslant k' < 1$, and for every pair $c, c' \in \mathrm{JC}^+$, if $c \subseteq \mathrm{Ext}(c')$, then*
$$\mathrm{K}_{k,h}(c) \subseteq \mathrm{K}_{k',h}(c'). \tag{2.19}$$

In order to prove Theorem 2.1.1, it remains only to prove continuity.

**Lemma 2.4.4**

*The function $\mathrm{K}_{k,h} : \mathrm{JC}^+ \to \mathrm{CC}_{k,h}$ is continuous.*

**Proof:** Let $(c_m)_{m \in \mathbb{N}}$ be a sequence in $\mathrm{JC}^+$ converging to some element $c_\infty$, say. For all $m$, denote $K_m := \mathrm{K}_{k,h}(c_m)$. Let $K_\infty$ be an accumulation point in CC of the sequence $(K_m)_{m \in \mathbb{N}}$. Let $c_\infty^\pm \in \mathrm{JC}^+$ be such that
$$c_\infty^- \subseteq \mathrm{Ext}(c_\infty) \text{ and } c_\infty \subseteq \mathrm{Ext}(c_\infty^+).$$

For sufficiently large $m$,
$$c_\infty^- \subseteq \mathrm{Ext}(c_m) \text{ and } c_m \subseteq \mathrm{Ext}(c_\infty^+).$$

By Lemma 2.4.3 that, for all such $m$,
$$\mathrm{K}_{k,h}(c_\infty^-) \subseteq K_m \subseteq \mathrm{K}_{k,h}(c_\infty^+),$$

so that, upon taking limits,
$$\mathrm{K}_{k,h}(c_\infty^-) \subseteq K_\infty \subseteq \mathrm{K}_{k,h}(c_\infty^+).$$

Hence
$$\overline{\mathrm{Ext}}(c_\infty^-) = \partial_\infty \mathrm{K}_{k,h}(c_\infty^-) \subseteq \partial_\infty K_\infty \subseteq \partial_\infty \mathrm{K}_{k,h}(c_\infty^+) = \overline{\mathrm{Ext}}(c_\infty^+).$$

Since $c_\infty^\pm$ are arbitrary, it follows that
$$\partial_\infty K_\infty = \overline{\mathrm{Ext}}(c_\infty^+).$$

In particular, by Corollary 2.2.6, $K_\infty \subseteq \mathrm{CC}_{k,h}$ so that, by uniqueness
$$K_\infty = \mathrm{K}_{k,h}(c_\infty).$$

The sequence $(K_m)_{m \in \mathbb{N}}$ thus only accumulates on $\mathrm{K}_{k,h}(c_\infty)$. Since CC is compact, it follows that $(K_m)_{m \in \mathbb{N}}$ converges towards $\mathrm{K}_{k,h}(c_\infty)$, as desired. $\square$

**Proof of Theorem 2.1.1:** Existence and uniqueness follow from Theorem 2.4.1. To prove continuity, note first that, by (2.18), $\mathrm{D}_{k,h} = \partial_{\mathrm{fin}} \circ \mathrm{K}_{k,h}$. By Lemma 2.4.4, $\mathrm{K}_{k,h}$ is continuous, whilst, by Lemma 2.2.4, $\partial_{\mathrm{fin}}$ is continuous, so that $\mathrm{D}_{k,h}$ is indeed continuous, as desired. $\square$

**2.5 - The first foliated Plateau problem.** Let $X := (X, h)$ be a Cartan–Hadamard manifold of sectional curvature pinched between $-C$ and $-1$, for some $C \geqslant 1$. In order to address the first foliated Plateau problem, it is necessary to understand the action of $\mathrm{K}_{k,h}$ on divergent sequences of curves.

**Lemma 2.5.1**

Let $(c_m)_{m \in \mathbb{N}}$ be a sequence of curves in $\mathrm{JC}^+$ and, for all $m$, denote $K_m := \mathrm{K}_{k,h}(c_m)$ and $D_m := \mathrm{D}_{k,h}(c_m)$. If $(\overline{\mathrm{Int}}(c_m))_{m \in \mathbb{N}}$ converges to the singleton $\{x_\infty\}$, then $(K_m)_{m \in \mathbb{N}}$ converges to $\overline{X}$ and $(\overline{D}_m)_{m \in \mathbb{N}}$ Hausdorff converges to $\{x_\infty\}$.

**Proof:** We first show that $(K_m)_{m \in \mathbb{N}}$ converges to $\overline{X}$. Indeed, let $K_\infty \in \mathrm{CC}$ be an accumulation point of $(K_m)_{m \in \mathbb{N}}$. Since $\mathrm{Ext}(c_m) \subseteq K_m$ for all $m$, it follows upon taking limits that $\partial_\infty X \subseteq K_\infty$ so that, by convexity, $\overline{X} \subseteq K_\infty$. $(K_m)_{m \in \mathbb{N}}$ thus only accumulates on $\overline{X}$. Since $\mathrm{CC}$ is compact, it follows that $(K_m)_{m \in \mathbb{N}}$ converges to $\overline{X}$, as desired.

We now show that $(\overline{D}_m)_{m \in \mathbb{N}}$ Hausdorff converges to $\{x_\infty\}$. Indeed, let $\overline{D}_\infty$ be a Hausdorff accumulation point of this sequence, and let $c \in \mathrm{JC}^+$ be such that $x_\infty \in \mathrm{Int}(c)$. Let $M$ be such that, for $m \geqslant M$, $\overline{\mathrm{Int}}(c_m) \subseteq \mathrm{Int}(c)$, so that $c \subseteq \mathrm{Ext}(c_m)$. By Lemma 2.4.3, for all $m \geqslant M$,

$$\mathrm{K}_{k,h}(c) \subseteq K_m.$$

Let $x$ be an interior point of $\mathrm{K}_{k,h}(c)$, let $S_x$ denote the unit sphere in $T_x X$, and let $\pi : \overline{X} \setminus \{x\} \to S_x$ denote the radial projection. By Lemma 2.2.1, for all $m \geqslant M$, the restriction of $\pi$ to $\hat{\partial} K_m$ is a homeomorphism, so that

$$\pi(\overline{D}_m) \cap \pi(\mathrm{Ext}(c)) \subseteq \pi(\overline{D}_m) \cap \pi(\mathrm{Ext}(c_m)) = \emptyset.$$

Upon taking limits, it follows that

$$\pi(\overline{D}_\infty) \cap \pi(\mathrm{Ext}(c)) = \emptyset.$$

However, by Lemma 2.2.2, $\overline{D}_\infty \subseteq \partial_\infty \overline{X}$, so that

$$\overline{D}_\infty \subseteq \overline{\mathrm{Int}}(c).$$

Since $c$ is arbitrary, it follows that $\overline{D}_\infty = \{x_\infty\}$. The sequence $(\overline{D}_m)_{m \in \mathbb{N}}$ thus Hausdorff accumulates only on $\{x_\infty\}$. Since the set of closed subsets of $\overline{X}$ is compact, it follows that $(\overline{D}_m)_{m \in \mathbb{N}}$ Hausdorff converges to $\{x_\infty\}$, as desired. $\square$

**Lemma 2.5.2**

Let $(c_m)_{m \in \mathbb{N}}$ be a sequence of curves in $\mathrm{JC}^+$, and for all $m$, denote $K_m := \mathrm{K}_{k,h}(c_m)$ and $D_m := \mathrm{D}_{k,h}(c_m)$. If $(\overline{\mathrm{Ext}}(c_m))_{m \in \mathbb{N}}$ converges to the singleton $\{x_\infty\}$, then $(K_m)_{m \in \mathbb{N}}$ converges to $\{x_\infty\}$ and $(D_m)_{m \in \mathbb{N}}$ Hausdorff converges to $\{x_\infty\}$.

**Proof:** It trivially suffices to show that $(K_m)_{m \in \mathbb{N}}$ converges to $\{x_\infty\}$. Let $K_\infty$ be an accumulation point of $(K_m)_{m \in \mathbb{N}}$. Let $c \in \mathrm{JC}^+$ be an oriented Jordan curve in $\partial_\infty X$ such that $x_\infty \in \mathrm{Ext}(c)$. For sufficiently large $m$, $c_m \subseteq \mathrm{Ext}(c)$ so that, by Lemma 2.4.3,

$$K_m \subseteq \mathrm{K}_{k,h}(c).$$

Upon taking limits it follows that

$$K_\infty \subseteq \mathrm{K}_{k,h}(c),$$

so that

$$\partial_\infty K_\infty \subseteq \partial_\infty \mathrm{K}_{k,h}(c) = \overline{\mathrm{Ext}}(c).$$

Since $c$ is arbitrary, it follows that

$$\partial_\infty K_\infty = \{x_\infty\}.$$

In particular, by Theorem 2.4.2, $K_\infty \notin \mathrm{CC}_{k,h}$ so that, by Corollary 2.2.6, $K_\infty = \{x_\infty\}$. The sequence $(K_m)_{m \in \mathbb{N}}$ thus only accumulates on the singleton $\{x_\infty\}$. Since $\mathrm{CC}$ is compact, it follows that $(K_m)_{m \in \mathbb{N}}$ converges towards $\{x_\infty\}$, as desired. $\square$

In order to complete the proof of Theorem 2.1.2, we require a version of the strong geometric maximum principle. For later applications, we will take a non-standard approach using the following unique continuation result of Aroszajn (see [2] and Theorem 2.3.4 of [29]).

**Theorem 2.5.3, Aronszajn (1957)**

Let $\Omega \subseteq \mathbb{R}^m$ be a connected open set. Let $u : \Omega \to \mathbb{R}$ be a smooth function, and suppose that there exists $C > 0$ such that, for all $x \in \Omega$,
$$|\Delta u(x)| \leqslant C\big(|u(x)| + \|Du(x)\|\big). \tag{2.20}$$
If the Taylor series of $u$ vanishes at some point, then $u$ vanishes identically.

In order to apply this result to our setting, we recall the concept of Fermi coordinates about an embedded surface in $X$. Let $S$ be an embedded surface in $X$, let $\nu : S \to SX$ denote its unit normal vector field compatible with the orientation and, for $\varepsilon > 0$, define $E : S \times \,]{-}\varepsilon, \varepsilon[\, \to X$ by
$$E(x, t) := \operatorname{Exp}(t\nu(x)). \tag{2.21}$$

Upon reducing $S$ and $\varepsilon$ if necessary, we may suppose that $E$ is a diffeomorphism onto its image. We call the resulting parametrization *Fermi coordinates* of $X$ about $S$. Given any function $f : S \to \,]{-}\varepsilon, \varepsilon[$, we call the image under $E$ of its graph the *Fermi graph* of $f$. Let $x_0$ be a base point of $S$. Upon composing with the exponential map of $S$, we may also identify this surface with a neighbourhood $U$ of the origin of its tangent space at that point. We may trivially suppose that $U$ is connected.

**Lemma 2.5.4**

Choose $k > 0$ and suppose that $S$ has constant extrinsic curvature equal to $k$. Let $f : U \to \,]{-}\varepsilon, \varepsilon[$ be a smooth function whose Fermi graph $S^f$ also has constant extrinsic curvature equal to $k$. If $f$ is non-vanishing, then its Taylor series is non-trivial. Furthermore, the lowest order term of this Taylor series is a polynomial $P : T_{x_0} S \to \mathbb{R}$ satisfying
$$\operatorname{Tr}((\mathrm{II}(x_0) + \mathrm{II}^f(x_0))^{-1} \operatorname{Hess}(P)) = 0, \tag{2.22}$$
where $\mathrm{II}$ and $\mathrm{II}^f$ here denote respectively the second fundamental forms of $S$ and $S^f$.

**Proof:** Indeed, consider a smooth function $f : U \to \,]{-}\varepsilon, \varepsilon[$ and let $\kappa_f$ denote the extrinsic curvature of its Fermi graph. Note that
$$\kappa_f(x) = F(x, \mathrm{J}^2 f(x)),$$
for some smooth function $F$, where $\mathrm{J}^2 f$ here denotes the 2-jet of $f$. There therefore exist smooth functions $a$, $b$ and $c$ such that
$$a(x, \mathrm{J}^2 f(x))^{ij} \operatorname{Hess}(f)(x)_{ij} + b(x, \mathrm{J}^2 f(x))^i Df(x)_i + c(x, \mathrm{J}^2 f(x)) f(x) = \kappa_0(x) - \kappa_f(x). \tag{2.23}$$

We now claim that
$$a(x, \mathrm{J}^2 f(x))^{ij} = \frac{1}{2} \operatorname{Adj}(\mathrm{II}(x) + \mathrm{II}^f(x))^{ij} + \operatorname{O}(\mathrm{J}^1 f(x)), \tag{2.24}$$

where Adj here denotes the adjugate operator of matrices, and $\operatorname{O}(\mathrm{J}^1 f(x))$ represents terms which vanish when $\mathrm{J}^1 f(x)$ vanishes. Indeed, choose $x \in S$, and suppose that $\mathrm{J}^1 f(x) = 0$. Then, by (2.10),

$$\mathrm{II}^f(x) = \operatorname{Hess}(t - (f \circ \pi))(x)|_{TS} = \operatorname{Hess}(t)(x)|_{TS} - \operatorname{Hess}(f \circ \pi)(x)|_{TS} = \mathrm{II}(x) - \operatorname{Hess}(f)(x).$$

Observe now that, for any pair $(A, B)$ of 2-dimensional matrices,
$$\operatorname{Det}(A) - \operatorname{Det}(B) = \frac{1}{2} \operatorname{Tr}(\operatorname{Adj}(A + B)(A - B)).$$

It follows that
$$\kappa_0(x) - \kappa_f(x) = \operatorname{Det}(\mathrm{II}(x)) - \operatorname{Det}(\mathrm{II}^f(x))$$
$$= \frac{1}{2} \operatorname{Tr}((\operatorname{Adj}(\mathrm{II}^f(x) + \mathrm{II}(x)) \operatorname{Hess}(f)(x))$$
$$= \frac{1}{2} \operatorname{Adj}(\mathrm{II}^f(x) + \mathrm{II}(x))^{ij} \operatorname{Hess}(f)(x)_{ij},$$

and (2.24) follows, as asserted.

Suppose now that $S^f$ has constant extrinsic curvature equal to $k$ and choose $x_0 \in \Omega$. We may suppose that $f(x_0), Df(x_0) = 0$, for otherwise the result holds trivially. By (2.23), near $x_0$,

$$a^{ij}\text{Hess}(f)_{ij} + b^i Df_i + cf = 0, \qquad (2.25)$$

for suitable $a$, $b$ and $c$. Since $\text{II}(x_0)$ and $\text{II}^f(x_0)$ are positive-definite, by (2.24), so too is $a^{ij}$. By Theorem 2.5.3, $f$ either has non-trivial Taylor series at $x_0$ or vanishes identically in a neighbourhood of this point. Finally, upon determining the lowest-order term of (2.25) using (2.24), we obtain (2.22), and this completes the proof. $\square$

**Corollary 2.5.5**

With the notation of Lemma 2.5.4, if $f(0) = Df(0) = 0$, then either $f$ vanishes identically or the critical point of $f$ at $0$ is isolated.

**Proof:** Indeed, by (2.22), up to precomposition with a linear map, the lowest order term of the Taylor series of $f$ at $0$ is a non-linear, homogeneous, harmonic polynomial. Since such functions have isolated critical points at the origin, the result follows. $\square$

In particular, Lemma 2.5.4 yields an alternative proof of the strong maximum principle for $k$-disks.

**Lemma 2.5.6**

Let $X$ be a riemannian manifold. Let $K$ and $K'$ be closed subsets of $X$ with non-trivial interiors such that $K \subseteq K'$. If $x \in \partial K \cap \partial K'$, and if $\partial K$ and $\partial K'$ are both smooth and of constant extrinsic curvature equal to $k$ near this point, then there exists a neighbourhood $\Omega$ of $x$ such that

$$\partial K \cap \Omega = \partial K' \cap \Omega. \qquad (2.26)$$

**Proof:** Let $\Omega$ be a neighbourhood of $x$ such that $\partial K \cap \Omega$ is a smooth surface $S$, say. Upon reducing $\Omega$ if necessary, we may suppose that $\partial K' \cap \Omega$ is the Fermi graph of some smooth function $f : S \to \mathbb{R}$. By hypothesis, $f \geqslant 0$. Suppose now that $f$ does not vanish near $x$. By Lemma 2.5.4, it has non-trivial Taylor series at this point, with lowest-order term $P$ satisfying

$$\text{Tr}((\text{II} + \text{II}^f)(x)^{-1}\text{Hess}(P)(x)) = 0.$$

Furthermore, since $f(0) = 0$, $P$ has degree at least $1$. However, harmonic polynomials of non-trivial degree assume positive and negative values in every neighbourhood of the origin. This is absurd, and it follows that $f$ vanishes near $0$, as desired. $\square$

We now prove Theorem 2.1.2.

**Proof of Theorem 2.1.2:** Let $S, N \subseteq \mathbb{S}^2$ denote respectively the south and north poles of $\mathbb{S}^2$. By Lemma 2.4.3, $(\text{K}_{k,h}(\alpha(c)))_{c \in \text{F}^+}$ is a nested family. We first show that this family covers $X$. Indeed, choose $x \in X$. By Lemma 2.5.1,

$$x \in \overline{X} = \lim_{c \to \{N\}} \text{K}_{k,h}(\alpha(c)) = \left(\bigcup_{c \in \text{F}^+} \text{K}_{k,h}(\alpha(c))\right) \cup \{\alpha(N)\},$$

and by Lemma 2.5.2,

$$x \notin \{S\} = \lim_{c \to \{S\}} \text{K}_{k,h}(\alpha(c)) = \bigcap_{c \in \text{F}^+} \text{K}_{k,h}(\alpha(c)).$$

It follows that

$$x \in \text{K}_{k,h}(\alpha(c_1)) \text{ and } x \notin \text{K}_{k,h}(\alpha(c_2)),$$

for some $c_1, c_2 \in \text{F}^+$. We now denote

$$I := \left\{c \in \text{F}^+ \mid x \in \text{K}_{k,h}(\alpha(c))\right\}.$$

Since $I$ is non-empty with non-empty complement, its boundary $\partial I$ is non-trivial. Let $c \in I$ be a boundary point. By continuity, $x \in \mathrm{K}_{k,h}(\alpha(c))$. We now claim that

$$x \in \partial_{\mathrm{fin}} \mathrm{K}_{k,h}(\alpha(c)) = \mathrm{D}_{k,h}(\alpha(c)).$$

Indeed, suppose the contrary, and let $r > 0$ be such that

$$B_{2r}(x) \cap \hat{\partial} \mathrm{K}_{k,h}(\alpha(c)) = \emptyset.$$

Since $c$ is a boundary point of $I$, there exists a sequence $(c_m)_{m \in \mathbb{N}} \in \mathrm{F}^+ \setminus I$ converging to $c$. By Lemma 2.2.2, there exists $M > 0$ such that, for all $m \geqslant M$,

$$B_r(x) \cap \hat{\partial} \mathrm{K}_{k,h}(\alpha(c_m)) = \emptyset.$$

In particular, since $x \notin K_{k,h}(\alpha(c_m))$, for all $m$,

$$B_r(x) \cap \mathrm{K}_{k,h}(\alpha(c_m)) = \emptyset,$$

so that taking limits yields

$$B_r(x) \cap \mathrm{K}_{k,h}(\alpha(c)) = \emptyset.$$

This is absurd, and it follows that $x \in \partial_{\mathrm{fin}} \mathrm{K}_{k,h}(\alpha(c))$, as asserted.

We have thus shown that the family $(\mathrm{D}_{k,h}(\alpha(c)))_{x \in \mathrm{F}^+}$ covers the whole of $X$. It remains only to show that it forms a foliation. Suppose, however, that there exists $c \neq c' \in \mathrm{F}^+$ such that

$$Y := \mathrm{D}_{k,h}(\alpha(c)) \cap \mathrm{D}_{k,h}(\alpha(c')) \neq \emptyset.$$

By Lemma 2.5.6, $Y$ is a relatively open subset of both $\mathrm{D}_{k,h}(\alpha(c))$ and $\mathrm{D}_{k,h}(\alpha(c'))$. Since it is trivially also closed, it follows by connectedness that $\mathrm{D}_{k,h}(\alpha(c)) = \mathrm{D}_{k,h}(\alpha(c'))$, so that

$$c = \partial_\infty \mathrm{D}_{k,h}(\alpha(c)) = \partial_\infty \mathrm{D}_{k,h}(\alpha(c')) = c'.$$

This is absurd, and it follows that $(\mathrm{D}_{k,h}(\alpha(c)))_{x \in \mathrm{F}^+}$ foliates $X$, as desired. $\square$

**2.6 - The second foliated Plateau problem.** The proof of Theorem 2.1.3 now follows by a standard argument from the theory of minimal surfaces (c.f. Chapter 6.2 of [9]).

**Lemma 2.6.1**

Let $X := (X, h)$ be a Cartan–Hadamard manifold. Choose $k > 0$, and let $D, D' \subseteq X$ be complete, properly embedded surfaces in $X$ of constant extrinsic curvature equal to $k$. If $D$ and $D'$ share the same normal at some point, then their intersection is a graph whose edges are smoothly embedded curves and whose vertices are isolated points with even valency at least $4$.

**Proof:** Denote $\Gamma := D \cap D'$. Near every point where $D$ meets $D'$ transversally, $\Gamma$ is a smooth curve. It thus suffices to determine the geometry of $\Gamma$ at points where these surfaces share a common tangent. Let $x_0 \in D \cap D'$ be such a point. We first claim that $D$ and $D'$ share a common normal at this point. Indeed, otherwise, their normals would point in opposite directions so that, by strict convexity, they would intersect at no other point. In particular, they would share no common normal at any point, which is absurd, and the claim follows. Now let $\Omega$ be a neighbourhood of $x_0$ in $X$ such that $D' \cap \Omega$ is the Fermi graph of some function $f$ over $D$. By Corollary 2.5.5, the critical point of $f$ at $x_0$ is isolated. For all sufficiently small $r > 0$, let $C_r$ denote the geodesic circle in $D$ of radius $r$ about $x_0$. By (2.22), the lowest order term of $f$ is, up to precomposition with a linear map, a harmonic polynomial of degree $m$, say, for some $m \geqslant 2$. It follows that, for all sufficiently small $r$, the set $f^{-1}(\{0\})$ intersects $C_r$ transversally at $2m$ distinct points. There therefore exists $\delta > 0$ and finitely many distinct, smooth curves $\gamma_1, \cdots, \gamma_{2m} : ]0, \delta[ \to D$ such that, for all $r \in ]0, \delta[$,

$$f^{-1}(\{0\}) \cap C_r = \{\gamma_1(r), \cdots, \gamma_{2m}(r)\}.$$

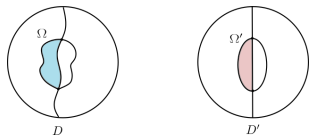

Figure 2.6.6 - **The intersection locus** - The intersection of $D$ and $D'$ is a graph $\Gamma$ contained in each of these surfaces all of whose vertices have valency at least 4.

The point $x_0$ is thus a vertex of $\Gamma$ of finite valency equal to $2m$, and this completes the proof. □

**Proof of Theorem 2.1.3:** We first show that the surfaces in the family $(\hat{D}_{k,h}(\alpha(c)))_{c \in \mathrm{C}^+}$ are pairwise disjoint. Indeed, suppose the contrary. Let $c, c' \in \mathrm{C}^+$ be such that $\hat{D} := \hat{D}_{k,h}(\alpha(c))$ and $\hat{D}' := \hat{D}_{k,h}(\alpha(c'))$ intersect non-trivially, and let $D, D' \subseteq X$ denote their respective images under the canonical projection. Upon applying a rotation of $\mathbb{S}^2$ if necessary, we may suppose that $c$ is an element of $\mathrm{F}^+$. By Lemma 2.6.1, $\Gamma := D \cap D'$ is a graph whose edges are smoothly embedded curves and whose vertices are isolated points of even valency at least 4 (see Figure 2.6.6).

Now let $x_0$ be a base point of $D$ and, for all $r$, let $B_r$ denote the geodesic ball of radius $r$ about $x_0$ in $D$. Choose $r$ such that $\partial B_r$ is transverse to $\Gamma$. We claim that if $B_r$ contains a single vertex of $\Gamma$, then $B_r \setminus \Gamma$ consists of at least 3 connected components. Note first that if $\Gamma$ does not meet $\partial B_r$, then we can join $\Gamma$ to $\partial B_r$ by adding an edge without changing the number of connected components of $B_r \setminus \Gamma$. Likewise, if $\Gamma \cap \overline{B}_r$ is not connected, we can also add edges to make it connected without changing the number of connected components of $B_r \setminus \Gamma$. Now let $V$ and $E$ denote respectively the sets of vertices and edges of $\Gamma$, and consider the graph $\Gamma_r$ whose vertex and edge sets are

$$V_r := \{v \in V \mid v \in B_r\} \cup \{e \cap \partial B_r \mid e \in E\}, \text{ and}$$
$$E_r := \{e \cap B_r \mid e \in E\} \cup \{\text{components of } \partial B_r \setminus \Gamma\}.$$

Let $F_r$ denote the set of connected components of $B_r \setminus \Gamma_r$. By Euler's formula

$$\#F_r - \#E_r + \#V_r = 1.$$

Since vertices lying on $\partial B_r$ are 3-valent and vertices lying on $B_r$ are at least 4-valent,

$$3\#(V_r \cap \partial B_r) + 4\#(V_r \cap B_r) \leqslant 2\#E_r.$$

Hence

$$\#F_r \geqslant 1 + \frac{1}{2}\#(V_r \cap \partial B_r) + \#(V_r \cap B_r) \geqslant \frac{5}{2},$$

so that $\#F_r \geqslant 3$, as asserted.

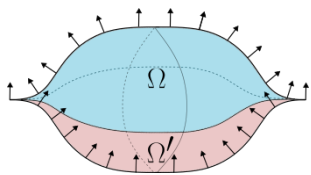

Figure 2.6.7 - **Towards a contradiction I** - We obtain a ball in $X$ whose boundary consists of a portion of $D$ joined to a portion of $D'$ along a Jordan curve in $\Gamma$.

Since $\partial_\infty D \cap \partial_\infty D' = c \cap c'$ consists of at most 2 points, $D \setminus \Gamma$ contains a connected component $\Omega$ which meets $\partial_\infty D$ at at most these two points. Note that $\partial \Omega$ is then a Jordan curve in $\overline{D}$ and also in $\overline{D}'$. Let $\Omega'$ denote the open subset of $D'$ that it bounds. Note that $\overline{\Omega} \cup \overline{\Omega}'$ is a Jordan sphere in $\overline{X}$ which therefore bounds an open ball $B$, say. Furthermore, since $D$ and $D'$ share a normal at at least one point of the common boundary of $\Omega$ and $\Omega'$, upon exchanging $D$ and $D'$ if necessary, we may suppose that the outward pointing normal over $\Omega'$ points into $B$, as in Figure 2.6.7. In particular, $\Omega'$ is contained in the interior of $K_{k,h}(\alpha(c))$. Since $\Omega'$ is relatively compact with boundary in $\overline{D}_{k,h}(\alpha(c))$, and since $(D_{k,h}(\alpha(c'')))_{c''\in F^+}$ foliates $X$, there exists $c'' \neq c \in F^+$ such that $K_{k,h}(\alpha(c''))$ intersects $\overline{B}$ only at points of $\Omega'$, and that $D_{k,h}(\alpha(c'')) = \partial_{\mathrm{fin}} K_{k,h}(\alpha(c''))$ and $\Omega'$ share the same outward-pointing normal at these points, as in Figure 2.6.8. It then follows by Lemma 2.5.6 that $D_{k,h}(\alpha(c')) = D_{k,h}(\alpha(c''))$. This is absurd, and it follows that $(\hat{D}_{k,h}(\alpha(c)))_{c\in C^+}$ consists of disjoint surfaces in $SX$. It remains only to show that this family covers $SX$. However, by conservation of the domain, the set that it covers is open whilst, by a straightforward compactness argument, it is also closed. It follows by connectedness that this family indeed covers $SX$, and this completes the proof. $\square$

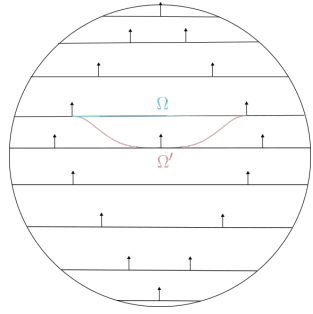

**Figure 2.6.8 - Towards a contradiction II -** The strong maximum principle together with the foliated Plateau problem addressed in the preceding section allows us to conclude that $D'$ coincides with another leaf of the foliation containing $D$, which is absurd.

We conclude this section by proving Theorems 2.1.4 and 2.1.5.

**Proof of Theorem 2.1.4:** We define the map $\phi_{k,h,\alpha} : SX \to \mathrm{MD}$ as follows. Choose $\xi \in SX$ with base point $x \in X$. Denote $p := (\alpha^{-1} \circ \mathrm{hor})(\xi)$ and $q := (\alpha^{-1} \circ \mathrm{hor})(-\xi)$. Let $\mathcal{D}$ denote the set of all disks which contain $p$ and which are preserved by the group of elliptic Möbius maps fixing both $p$ and $q$. By Theorem 2.1.2, there exists a unique disk $D \in \mathcal{D}$ such that $x \in D_{k,h}(\alpha(\partial D))$. We denote $\phi_{k,h,\alpha}(\xi) := (D,p)$. This function is trivially continuous, and we now show that it inverts $\Phi_{k,h,\alpha}$. First choose $(D,p) \in \mathrm{MD}$. Denote $q := R_{\partial D} p$. Let $\gamma : \mathbb{R} \to X$ denote the unique unit-speed geodesic such that $\gamma(+\infty) = \alpha(p)$, $\gamma(-\infty) = \alpha(q)$, and $\gamma(0) \in D_{k,h}(\alpha(\partial D))$, and denote $\xi := \dot\gamma(0)$. By definition,

$$\Phi_{k,h,\alpha}(D,p) = \xi.$$

However, $p = (\alpha^{-1} \circ \mathrm{hor})(\xi)$, $q = (\alpha^{-1} \circ \mathrm{hor})(-\xi)$, and $D$ is the unique disk preserved by the group of elliptic Möbius maps fixing both $p$ and $q$ such that $x := \gamma(0) \in D_{k,a}(\alpha(\partial D))$. Consequently,

$$\phi_{k,h,\alpha}(\xi) = (D,p).$$

Since $(D,p)$ is arbitrary, it follows that

$$\phi_{k,h,\alpha} \circ \Phi_{k,h,\alpha} = \mathrm{Id}.$$

Upon reversing this argument, we see that

$$\Phi_{k,h,\alpha} \circ \phi_{k,h,\alpha} = \mathrm{Id},$$

as desired. $\square$

**Proof of Theorem 2.1.5:** We define the map $\psi_{k,h,\alpha} : SX \to \mathrm{MD}$ as follows. Choose $\xi \in SX$. Denote $p := (\alpha^{-1} \circ \mathrm{hor})(\xi)$. By Theorem 2.1.3, there exists a unique $c \in C^+$ such that $\xi \in \hat{D}_{k,h}(\alpha(c))$. We define $\psi_{k,h,\alpha}(\xi) := (\mathrm{Int}(c), p)$. This function is trivially continuous, and we now show that it inverts $\Psi_{k,h,\alpha}$. First choose $(D,p) \in \mathrm{MD}$. Let $x \in D_{k,h}(\alpha(\partial D))$ denote the orthogonal projection of $\alpha(p)$, and let $\xi$ denote the outward-pointing unit normal to $D_{k,h}(\alpha(\partial D))$ at this point. By definition,

$$\Psi_{k,h,\alpha}(D,p) = \xi.$$

However, $p = (\alpha^{-1} \circ \text{hor})(\xi)$, and $\partial D$ is trivially the unique element of $\text{C}^+$ such that $\xi \in \hat{\text{D}}_{k,h}(\alpha(\partial D))$, so that
$$\psi_{k,h,\alpha}(\xi) = (D, p).$$

Since $(D, p)$ is arbitrary, it follows that
$$\psi_{k,h,\alpha} \circ \Psi_{k,h,\alpha} = \text{Id}.$$

Upon reversing this argument, we see also that
$$\Psi_{k,h,\alpha} \circ \psi_{k,h,\alpha} = \text{Id},$$

and the result follows. $\square$

## 3 - Fibred Plateau problems.

**3.1 - Principal $\text{PSL}(2, \mathbb{R})$-bundles.** The second foliated Plateau problem of Theorem 2.1.3 lies at the centre of a family of fibred Plateau problems which, as we will see in Section 4, provides the natural framework for developing the equidistribution theory described by Labourie in [24]. Let $X := (X, h)$ be an oriented, 3-dimensional Cartan–Hadamard manifold of sectional curvature bounded above by $-1$, acted on cocompactly by the hyperbolic group $\Pi$. Recall that this yields the following natural identifications.
$$\partial_\infty X = \partial_\infty \Pi = \partial_\infty \mathbb{H}^3 = \hat{\mathbb{C}}. \tag{3.1}$$

For all $K$, let $\text{QC}^+(K)$ denote the space of oriented $K$-quasicircles in $\hat{\mathbb{C}}$ furnished with the Hausdorff topology. Note that this is a nested family with
$$\text{C}^+ := \text{QC}^+(1) = \bigcap_{K \geqslant 1} \text{QC}^+(K). \tag{3.2}$$

We denote
$$\text{QC}^+ := \text{QC}^+(\infty) := \bigcup_{K \geqslant 1} \text{QC}^+(K), \tag{3.3}$$

and we furnish this space with the colimit topology. Recall that this is the finest topology over $\text{QC}^+$ with respect to which the canonical inclusion $\text{QC}^+(K) \subseteq \text{QC}^+$ is continuous for all $K$. More concretely, a sequence $(c_m)_{m \in \mathbb{N}}$ in $\text{QC}^+$ converges to $c_\infty$ in the colimit topology if and only if there exists $K \geqslant 1$ such that $c_m \in \text{QC}^+(K)$ for all $m \in \mathbb{N} \cup \{\infty\}$, and $(c_m)_{m \in \mathbb{N}}$ converges to $c_\infty$ in this subspace.

For all $c \in \text{QC}^+$, we denote its quasisymmetry constant by $\text{QS}(c)$, so that
$$\text{QS}(c) := \text{Inf}\{K \mid c \in \text{QC}^+(K)\}. \tag{3.4}$$

We will see presently, in Corollary 3.3.5, that $\text{QS}$ defines a lower semicontinuous function from $\text{QC}^+$ into $[1, \infty[$. That is, for every sequence $(c_m)_{m \in \mathbb{N}}$ in $\text{QC}^+$ converging to $c_\infty$,
$$\text{QS}(c_\infty) \leqslant \liminf_{m \to \infty} \text{QS}(c_m). \tag{3.5}$$

We will make considerable use of this function presently.

We will express our results in the language of topological bundles. Since these are less standard than smooth bundles, we recall their definition. Given two topological spaces $B$ and $F$, an *F-bundle* over $B$ is defined to be a pair $(T, \pi)$, where $T$ is a topological space and $\pi : T \to B$ is a continuous function having the property that, for every point $b \in B$, there exists a neighbourhood $\Omega$ of $b$ in $B$ and a function $f : \Omega \times F \to T$ such that

(1) $f$ is a homeomorphism onto its image; and

(2) $(\pi \circ f)$ coincides with projection onto the first factor.

We call $B$ the *base*, $T$ the *total space*, and $(\Omega, f)$ a *trivialization* of $(T, \pi)$ about $b$. Given a Lie group $G$, we define a *principal $G$-bundle* over $B$ to be a $G$-bundle $(P, \pi)$ over $B$ furnished with a continuous right $G$-action $\alpha : P \times G \to P$ which preserves every fibre and which acts transitively without fixed points over every fibre.

We will be interested in principal $\mathrm{PSL}(2, \mathbb{R})$-bundles over $\mathrm{QC}^+$. Let $\mathbb{H}$ denote the upper half-space in $\mathbb{C}$. For all $K \geqslant 1$, let $\mathrm{UQD}(K)$ denote the space of univalent maps $\alpha : \mathbb{H} \to \hat{\mathbb{C}}$ whose image is a $K$-quasidisk, furnished with the compact-open topology. We denote

$$\mathrm{UQD} := \mathrm{UQD}(\infty) := \bigcup_{K \geqslant 1} \mathrm{UQD}(K), \tag{3.6}$$

and we furnish this space with the colimit topology. We define $\partial : \mathrm{UQD} \to \mathrm{QC}^+$ by

$$\partial \alpha := \partial \mathrm{Im}(\alpha), \tag{3.7}$$

where this curve is oriented so that $\mathrm{Im}(\alpha)$ lies on its inside. Note that precomposition defines a continuous right-action of $\mathrm{PSL}(2, \mathbb{R})$ on $\mathrm{UQD}$.

**Theorem 3.1.1**

*The operator $\partial$ is continuous, and, with the action of $\mathrm{PSL}(2, \mathbb{R})$ on $\mathrm{UQD}$ by precomposition, $(\mathrm{UQD}, \partial)$ is a principal $\mathrm{PSL}(2, \mathbb{R})$-bundle over $\mathrm{QC}^+$.*

Theorem 3.1.1 is proven in Section 3.4.

For all $0 < k < 1$, and for all $K \geqslant 1$, let $\mathrm{KD}_{k,h}(K)$ denote the space of $k$-disks in $X$ spanning $K$-quasicircles in $\partial_\infty X$, furnished with the $C^\infty_{\mathrm{loc}}$ topology, and let $\mathrm{UKD}_{k,h}(K)$ denote the space of conformal parametrizations $\alpha : \mathbb{H} \to D$, for some $D \in \mathrm{KD}_{k,h}(K)$, furnished with the compact-open topology. For all $0 < k < 1$, we denote

$$\mathrm{UKD}_{k,h} := \mathrm{UKD}_{k,h}(\infty) := \bigcup_{K \geqslant 1} \mathrm{UKD}_{k,h}(K), \tag{3.8}$$

and we furnish this space with the colimit topology. We define $\partial_\infty : \mathrm{UKD}_{k,h} \to \mathrm{QC}^+$ by

$$\partial_\infty \alpha := \partial_\infty \alpha(\mathbb{H}), \tag{3.9}$$

and we furnish this curve with the orientation that it inherits from $\alpha(\mathbb{H})$. As before, for all $0 < k < 1$, precomposition defines a continuous right-action of $\mathrm{PSL}(2, \mathbb{R})$ on $\mathrm{UKD}_{k,h}$.

**Theorem 3.1.2**

*For all $0 < k < 1$, the operator $\partial_\infty$ is continuous, and, with the action of $\mathrm{PSL}(2, \mathbb{R})$ on $\mathrm{UKD}_{k,h}$ by precomposition, $(\mathrm{UKD}_{k,h}, \partial_\infty)$ is a principal $\mathrm{PSL}(2, \mathbb{R})$-bundle over $\mathrm{QC}^+$.*

Theorem 3.1.2 is proven in Section 3.5.

**3.2 - Associated $\mathrm{PSL}(2, \mathbb{R})$-bundles.** Let $G$ be a Lie group and let $(P, \pi)$ be a principal $G$-bundle with base $B$. Given another topological space $F$ over which $G$ acts continuously, we define the left $G$-action on the cartesian product $P \times F$ such that, for all $\gamma \in G$, and for all $(p, f) \in P \times F$,

$$\gamma \cdot (p, f) = (p \cdot \gamma^{-1}, \gamma \cdot f). \tag{3.10}$$

We denote the quotient of $P \times G$ under this action by $P \otimes_G F$. Note that $\pi$ descends to a continuous projection from this space to $B$ which makes $(P \otimes_G F, \pi)$ into an $F$-bundle over $B$. We call any bundle constructed in this manner an *associated $G$-bundle* over $B$. Given an open subset $\Omega \subseteq B$ and a continuous section $\sigma : \Omega \to P$, we define the trivialization $\alpha_\sigma : \Omega \times F \to P \otimes_G F$ by

$$\alpha_\sigma(b, f) := [\sigma(b), f]. \tag{3.11}$$

We call any trivialization constructed in this manner a *G-trivialization* of $P \otimes_G F$. Given two *G*-trivializations $\alpha, \alpha' : \Omega \times F \to P \otimes_G F$, the transition map $\tau := \alpha^{-1} \circ \alpha'$ takes the form

$$\tau(x, f) := (x, \beta(x) \cdot f), \tag{3.12}$$

for some continuous function $\beta : \Omega \to G$. The reader may verify that an associated *G*-bundle structure over a given *F*-bundle $(T, \pi)$ is equivalent to an atlas of trivializations, all of whose transition maps have this form.

We now study bundles over $\mathrm{QC}^+$ associated to the principal bundles constructed in the previous section. For all $K \geqslant 1$, let $\mathrm{QD}(K)$ denote the space of *K*-quasidisks in $\hat{\mathbb{C}}$, furnished with the Hausdorff topology. For all such *K*, let $\mathrm{MQD}(K)$ denote the space of marked *K*-quasidisks in $\hat{\mathbb{C}}$, that is

$$\mathrm{MQD}(K) := \{(D, z) \mid D \in \mathrm{QD}(K) \text{ and } z \in D\}, \tag{3.13}$$

furnished with the topology that it inherits as a subset of $\mathrm{QD}(K) \times \hat{\mathbb{C}}$. We denote

$$\mathrm{MQD} := \mathrm{MQD}(\infty) := \bigcup_{K \geqslant 1} \mathrm{MQD}(K), \tag{3.14}$$

and we furnish this space with the colimit topology. We define $\partial : \mathrm{MQD} \to \mathrm{QC}^+$ by

$$\partial(D, z) := \partial D. \tag{3.15}$$

Finally, we define $\Phi : \mathrm{UQD} \otimes_{\mathrm{PSL}(2,\mathbb{R})} \mathbb{H} \to \mathrm{MQD}$ by

$$\Phi([\alpha, z]) := (\mathrm{Im}(\alpha), \alpha(z)). \tag{3.16}$$

**Theorem 3.2.1**

$\Phi : \mathrm{UQD} \otimes_{\mathrm{PSL}(2,\mathbb{R})} \mathbb{H} \to \mathrm{MQD}$ *is a bundle homeomorphism. In particular, the operator $\partial$ is continuous and* $(\mathrm{MQD}, \partial)$ *is an $\mathbb{H}$-bundle over $\mathrm{QC}^+$.*

Theorem 3.2.1 is proven in Section 3.4.

For all $0 < k < 1$, and for all $K \geqslant 1$, recall that $\mathrm{KD}_{k,h}(K)$ denotes the space of *k*-disks in $(X, h)$ spanning *K*-quasicircles in $\hat{\mathbb{C}}$, furnished with the $C^\infty_{\mathrm{loc}}$ topology. For all such *k* and *K*, let $\mathrm{MKD}_{k,h}(K)$ denote the space of marked *k*-disks in $(X, h)$ spanning *K*-quasicircles in $\hat{\mathbb{C}}$, that is

$$\mathrm{MKD}_{k,h}(K) := \{(D, p) \mid D \in \mathrm{KD}_{k,h}(K) \text{ and } p \in D\}, \tag{3.17}$$

furnished with the topology that it inherits as a subset of $\mathrm{KD}_{k,h}(K) \times X$. For all $0 < k < 1$, we denote

$$\mathrm{MKD}_{k,h} := \mathrm{MKD}_{k,h}(\infty) := \bigcup_{K \geqslant 1} \mathrm{MKD}_{k,h}(K), \tag{3.18}$$

and we furnish this space with the colimit topology. We define $\partial_\infty : \mathrm{MKD}_{k,h} \to \mathrm{QC}^+$ by

$$\partial_\infty(D, z) := \partial_\infty D. \tag{3.19}$$

Finally, we define $\Phi_{k,h} : \mathrm{UKD}_{k,h} \otimes_{\mathrm{PSL}(2,\mathbb{R})} \mathbb{H} \to \mathrm{MKD}_{k,h}$ by

$$\Phi_{k,h}([\alpha, z]) := (\mathrm{Im}(\alpha), \alpha(z)). \tag{3.20}$$

**Theorem 3.2.2**

$\Phi_{k,h} : \mathrm{UKD}_{k,h} \otimes_{\mathrm{PSL}(2,\mathbb{R})} \mathbb{H} \to \mathrm{MKD}_{k,h}$ is a bundle homeomorphism. In particular, the operator $\partial_\infty$ is continuous and $\mathrm{MKD}_{k,h}$ is an $\mathbb{H}$-bundle over $\mathrm{QC}^+$.

Theorem 3.2.2 is proven in Section 3.5.

**3.3 - Quasiconformal homeomorphisms and quasicircles.** We now review the structures of the space of quasiconformal homeomorphisms and quasicircles. For completeness, we study this in more depth than is necessary for our present purposes. For all $K \geqslant 1$, let $\mathrm{QH}^+(K)$ denote the space of orientation-preserving, $K$-quasiconformal homeomorphisms of $\hat{\mathbb{C}}$ furnished with the compact-open topology. Note that this is a nested family, with

$$\mathrm{PSL}(2,\mathbb{C}) = \mathrm{QH}^+(1) = \bigcap_{K \geqslant 1} \mathrm{QH}^+(K). \tag{3.21}$$

More generally, this family is semicontinuous in the sense that, for all $K$,

$$\mathrm{QH}^+(K) = \bigcap_{K' \geqslant K} \mathrm{QH}^+(K'). \tag{3.22}$$

We define

$$\mathrm{QH}^+ := \mathrm{QH}^+(\infty) := \bigcup_{K \geqslant 1} \mathrm{QH}^+(K), \tag{3.23}$$

and we furnish this space with the colimit topology. Note that, for all $\alpha \in \mathrm{QH}^+(K)$ and $\beta \in \mathrm{QH}^+(K')$,

$$\alpha \circ \beta^{-1} \in \mathrm{QH}^+(KK'), \tag{3.24}$$

so that $\mathrm{QC}^+$ is a topological group.

We will require the following compactness criterion for closed subsets of $\mathrm{QH}^+$ (see Theorem 5.1 of [25]). Let $\Delta \subseteq \hat{\mathbb{C}}^3$ denote the diagonal, that is

$$\Delta := \{(z_1, z_2, z_3) \mid \#\{z_1, z_2, z_3\} \leqslant 2\}. \tag{3.25}$$

Given a Möbius map $\mu$, we define $\mathrm{T}_\mu : \mathrm{QH}^+ \to \hat{\mathbb{C}}^3 \setminus \Delta$ by

$$\mathrm{T}_\mu(\alpha) := ((\alpha \circ \mu)(0), (\alpha \circ \mu)(1), (\alpha \circ \mu)(\infty)). \tag{3.26}$$

**Theorem 3.3.1**

Let $X$ be a closed subset of $\mathrm{QH}^+(K)$, for some $K \geqslant 1$. The following affirmations are equivalent.

(1) $X$ is compact.

(2) $\mathrm{T}_\mu(X)$ does not accumulate on $\Delta$ for some $\mu \in \mathrm{PSL}(2,\mathbb{C})$.

(3) $\mathrm{T}_\mu(X)$ does not accumulate on $\Delta$ for all $\mu \in \mathrm{PSL}(2,\mathbb{C})$.

We now realize $\mathrm{QC}^+$ as a homogeneous space of $\mathrm{QH}^+$. To this end, for all $K \geqslant 1$, let $\mathrm{QH}^+(K, \hat{\mathbb{R}})$ denote the space of orientation-preserving $K$-quasiconformal homeomorphisms of $\hat{\mathbb{C}}$ which preserve the oriented extended real line $\hat{\mathbb{R}}$. We define

$$\mathrm{QH}^+(\hat{\mathbb{R}}) := \mathrm{QH}^+(\infty, \hat{\mathbb{R}}) := \bigcup_{K \geqslant 1} \mathrm{QH}^+(K, \hat{\mathbb{R}}), \tag{3.27}$$

and we furnish this space with the colimit topology. Note that $\mathrm{QH}^+(\hat{\mathbb{R}})$ is a closed subgroup of $\mathrm{QH}^+$.

Define $\tilde{c} : \mathrm{QH}^+ \to \mathrm{QC}^+$ by

$$\tilde{c}(\alpha) := \alpha(\hat{\mathbb{R}}), \tag{3.28}$$

where this curve is furnished with the natural orientation inherited from $\hat{\mathbb{R}}$. Trivially, this function is continuous and, for all $\alpha, \beta \in \mathrm{QH}^+$,

$$\tilde{c}(\alpha) = \tilde{c}(\beta) \iff \alpha^{-1} \circ \beta \in \mathrm{QH}^+(\hat{\mathbb{R}}). \tag{3.29}$$

so that $\tilde{c}$ descends to a continuous bijection

$$c : \mathrm{QH}^+ / \mathrm{QH}^+(\hat{\mathbb{R}}) \to \mathrm{QC}^+; [\alpha] \mapsto \tilde{c}(\alpha). \tag{3.30}$$

**Theorem 3.3.2**

c *defines a homeomorphism from* $\mathrm{QH}^+/\mathrm{QH}^+(\hat{\mathbb{R}})$ *into* $\mathrm{QC}^+$.

This result is a consequence of the following useful lemma.

**Lemma 3.3.3**

*Let* $\pi : \mathrm{QH}^+ \to \mathrm{QH}^+/\mathrm{PSL}(2,\mathbb{R})$ *denote the canonical projection. If $A$ is a subset of $\mathrm{QH}^+(K)$, for some $K \geqslant 1$, then $\pi(A)$ is precompact in $\mathrm{QH}^+/\mathrm{PSL}(2,\mathbb{R})$ if and only if $\mathrm{c}(A)$ does not accumulate on any singleton.*

**Proof:** Suppose first that $\mathrm{c}(A)$ does not accumulate on any singleton. Let $(\alpha_m)_{m\in\mathbb{N}}$ be a sequence in $A$ and, for all $m$, denote $c_m := \tilde{\mathrm{c}}(\alpha_m)$. We may suppose that $(c_m)_{m\in\mathbb{N}}$ Hausdorff converges to some non-trivial, closed subset $c_\infty$, say, of $\hat{\mathbb{C}}$ which, by hypothesis, contains at least two distinct points $z_1$ and $z_2$, say. By connectedness, $c_\infty$ must also contain a third point $z_3$, say. There therefore exists a sequence $(\beta_m)_{m\in\mathbb{N}}$ of elements of $\mathrm{PSL}(2,\mathbb{R})$ such that
$$\lim_{m\to\infty} \mathrm{T}_{\mathrm{Id}}(\alpha_m \circ \beta_m) = (z_1, z_2, z_3).$$
In particular, $(\mathrm{T}_{\mathrm{Id}}(\alpha_m \circ \beta_m))_{m\in\mathbb{N}}$ does not accumulate on the diagonal. It follows by Theorem 3.3.1 that $(\alpha_m \circ \beta_m)_{m\in\mathbb{N}}$ is precompact in $\mathrm{QC}^+$, so that $(\pi(\alpha_m))_{m\in\mathbb{N}}$ is precompact in $\mathrm{QC}^+/\mathrm{PSL}(2,\mathbb{R})$. Since $(\alpha_m)_{m\in\mathbb{N}}$ in $A$ is arbitrary, it follows that $\pi(A)$ is precompact, as desired. The converse is trivial, and this completes the proof. $\square$

**Proof of Theorem 3.3.2:** It remains only to show that, given a sequence $(\alpha_m)_{m\in\mathbb{N}}$ in $\mathrm{QH}^+$, if $(\tilde{\mathrm{c}}(\alpha_m))_{m\in\mathbb{N}}$ converges in $\mathrm{QC}^+$, then the sequence $([\alpha_m])_{m\in\mathbb{N}}$ of equivalence classes converges in $\mathrm{QH}^+/\mathrm{QH}^+(\hat{\mathbb{R}})$. Suppose, however, that $(\tilde{\mathrm{c}}(\alpha_m))_{m\in\mathbb{N}}$ converges to $c_\infty$, say. By definition of the colimit topology, we may suppose that there exists $K > 0$ such that, for all $m$, $c_m \in \mathrm{QC}^+(K)$, and thus also that $\alpha_m \in \mathrm{QH}^+(K)$. By Lemma 3.3.3, the sequence $(\pi(\alpha_m))_{m\in\mathbb{N}}$ is precompact $\mathrm{QH}^+/\mathrm{PSL}(2,\mathbb{R})$, from which it follows that the sequence $([\alpha_m])_{m\in\mathbb{N}}$ is also precompact in $\mathrm{QH}^+/\mathrm{QH}^+(\hat{\mathbb{R}})$. By continuity, every accumulation point $[\alpha_\infty]$ of this sequence satisfies $\mathrm{c}([\alpha_\infty]) = c_\infty$, so that, by (3.29), $([\alpha_m])_{m\in\mathbb{N}}$ has a unique accumulation point, and therefore converges, as desired. $\square$

Theorem 3.3.2 and Lemma 3.3.3 together yield the well-known precompactness criterion for quasicircles.

**Theorem 3.3.4**

*Let $C$ be a closed subset of $\mathrm{QC}^+(K)$ for some $K \geqslant 1$. This set is compact if and only if it does not accumulate on any singleton.*

**Corollary 3.3.5**

QS *is lower semicontinuous over* $\mathrm{QC}^+$.

**Proof:** Indeed, by Theorem 3.3.4, for all $t \geqslant 1$, $\mathrm{QC}^+(t)$ is closed. However, by (3.22) and Lemma 3.3.3, for all $t \geqslant 1$, $\mathrm{QS}^{-1}(]-\infty,t]) = \mathrm{QC}^+(t)$, so that QS is indeed lower semicontinuous, as desired. $\square$

**3.4 - Principal $\mathrm{PSL}(2,\mathbb{R})$ bundles I.** We now show that $\partial : \mathrm{UQD} \to \mathrm{QC}^+$ is continous and that $(\mathrm{UQD}, \partial)$ is a principal $\mathrm{PSL}(2,\mathbb{R})$-bundle over $\mathrm{QC}^+$. First recall from Theorem 14.19 of [31] and the subsequent remark that every $\alpha \in \mathrm{UQD}$ extends to a homeomorphism $\overline{\alpha} : \overline{\mathbb{H}} \to \overline{\alpha(\mathbb{H})}$. In fact, since $\partial \alpha$ is a quasicircle, we have the following stronger result.

**Lemma 3.4.1**

*For all $K \geqslant 1$, every $\alpha \in \mathrm{UQD}(K)$ extends to a $K^2$-quasiconformal homeomorphism $\tilde{\alpha} : \hat{\mathbb{C}} \to \hat{\mathbb{C}}$.*

**Proof:** Let $\beta$ be a $K$-quasiconformal homeomorphism such that $\tilde{\mathrm{c}}(\beta) = \partial \alpha$, and let $\mu$ denote its Beltrami differential. Define the Beltrami differential $\tilde{\mu}$ by
$$\tilde{\mu}(z) := \begin{cases} \mu(z) & \text{if } \mathrm{Im}(z) \geqslant 0, \text{ and} \\ \overline{\mu}(\overline{z}) & \text{otherwise.} \end{cases}$$
Let $\gamma$ be a $K$-quasiconformal homeomorphism with Beltrami differential $\tilde{\mu}$ fixing the points $0$, $1$ and $\infty$, and observe that $\gamma$ preserves $\hat{\mathbb{R}}$. Denote $\alpha' := \beta \circ \gamma^{-1}$ and observe that its restriction to $\mathbb{H}$ is conformal. That is, $\alpha'$ restricts to a uniformising map of $\alpha(\mathbb{H})$ and is thus equal to $\alpha \circ \delta$, for some Möbius map $\delta$, say. In particular, $\alpha$ extends to the $K^2$-quasiconformal homeomorphism $\alpha' \circ \delta^{-1}$, and this completes the proof. $\square$

**Lemma 3.4.2**

Let $(\alpha_m)_{m\in\mathbb{N}}$ be a sequence in $\mathrm{UQD}(K)$, for some $K \geqslant 1$. If $(\alpha_m)_{m\in\mathbb{N}}$ converges to $\alpha_\infty$, say, then $(\overline{\alpha}_m)_{m\in\mathbb{N}}$ converges to $\overline{\alpha}_\infty$ uniformly over $\overline{\mathbb{H}}$.

**Proof:** For all $m$, let $\tilde{\alpha}_m$ be a $K^2$-quasiconformal extension of $\alpha_m$ as in Lemma 3.4.1. Let $\mu$ be a Möbius map such that $\mu(0), \mu(1), \mu(\infty) \in \mathbb{H}$. Since $(\alpha_m)_{m\in\mathbb{N}}$ converges to $\alpha_\infty$, $(\mathrm{T}_\mu(\tilde{\alpha}_m))_{m\in\mathbb{N}}$ does not accumulate on $\Delta$. It follows by Theorem 3.3.1 that $(\tilde{\alpha}_m)_{m\in\mathbb{N}}$ is precompact in $\mathrm{QH}^+(K^2)$, and the sequence $(\overline{\alpha}_m)_{m\in\mathbb{N}}$ is thus precompact in the uniform topology over $\overline{\mathbb{H}}$. Since $\overline{\alpha}_\infty$ is trivially the unique accumulation point of this sequence, it follows that $(\overline{\alpha}_m)_{m\in\mathbb{N}}$ converges uniformly to $\overline{\alpha}_\infty$, as desired. □

**Lemma 3.4.3**

$\partial : \mathrm{UQD} \to \mathrm{QC}^+$ is continuous.

**Proof:** Let $(\alpha_m)_{m\in\mathbb{N}}$ be a sequence in $\mathrm{UQD}$ converging to $\alpha_\infty$, say. For all $m \in \mathbb{N} \cup \{\infty\}$, let $\overline{\alpha}_m$ denote the continuous extension of $\alpha_m$. By Lemma 3.4.2, $(\overline{\alpha}_m)_{m\in\mathbb{N}}$ converges uniformly over $\overline{\mathbb{H}}$ to $\overline{\alpha}_\infty$. Consequently,

$$\partial \alpha_\infty = \overline{\alpha}_\infty(\hat{\mathbb{R}}) = \lim_{m\to\infty} \overline{\alpha}_m(\hat{\mathbb{R}}) = \lim_{m\to\infty} \partial \alpha_m,$$

as desired. □

**Lemma 3.4.4**

Let $(\alpha_m)_{m\in\mathbb{N}}$ be a sequence in $\mathrm{UQD}(K)$, for some $K \geqslant 1$ and, for all $m$, denote $c_m := \partial \alpha_m$. If $(c_m)_{m\in\mathbb{N}}$ converges to $c_\infty$, say, and if $(\alpha_m(i))_{m\in\mathbb{N}}$ does not accumulate on $c_\infty$, then $(\alpha_m)_{m\in\mathbb{N}}$ subconverges to some element $\alpha_\infty$, say.

**Proof:** Indeed, for all $m$, let $\tilde{\alpha}_m : \hat{\mathbb{C}} \to \hat{\mathbb{C}}$ be a $K^2$-quasiconformal extension of $\alpha_m$, as in Lemma 3.4.1. For all $m$, $\tilde{c}(\tilde{\alpha}_m) = c_m$ so that, by Lemma 3.3.3, we may suppose that there exists a sequence $(\beta_m)_{m\in\mathbb{N}}$ in $\mathrm{PSL}(2,\mathbb{R})$ such that $(\tilde{\alpha}_m \circ \beta_m)_{m\in\mathbb{N}}$ converges uniformly to some limit $\tilde{\alpha}'_\infty$, say. Since $(\tilde{\alpha}_m(i))_{m\in\mathbb{N}}$ does not accumulate on $c_\infty$, we may suppose furthermore that $(\beta_m)_{m\in\mathbb{N}}$ converges to $\beta_\infty$, say, so that $(\tilde{\alpha}_m)_{m\in\mathbb{N}}$ converges uniformly to $\tilde{\alpha}_\infty := \tilde{\alpha}'_\infty \circ \beta_\infty^{-1}$. In particular, $(\alpha_m)_{m\in\mathbb{N}}$ converges to $\alpha_\infty$, as desired. □

We now prove Theorems 3.1.1 and 3.2.1.

**Proof of Theorem 3.1.1:** Since continuity is proven in Lemma 3.4.3, it remains only to verify the local trivialization property. Choose $c_0 \in \mathrm{QC}^+$, $z \in \mathrm{Int}(c_0)$ and $\xi \in T_z\hat{\mathbb{C}}$. Let $U$ be a neighbourhood of $c_0$ in $\mathrm{QC}^+$ such that $z \in \mathrm{Int}(c)$ for all $c \in U$. For all $c \in U$, let $\alpha_c$ denote the unique univalent map $\alpha_c : \mathbb{H} \to \mathrm{Int}(c)$ such that $\alpha_c(i) = z$ and $D\alpha_c(i) \cdot \partial_x \propto \xi$. By Lemma 3.4.4, $\alpha_c$ depends continuously on $c$. We now define $A : U \times \mathrm{PSL}(2,\mathbb{R}) \to \mathrm{UQD}$ by

$$A(c, \beta) := \alpha_c \circ \beta.$$

Trivially, this function is continuous and maps $U \times \mathrm{PSL}(2,\mathbb{R})$ bijectively onto $\partial^{-1}(U)$. By Lemma 3.4.3, its inverse is also continuous. It is thus a local trivialization of $(\mathrm{UQD}, \partial)$ over $c_0$, as desired. □

**Proof of Theorem 3.2.1:** Since $\Phi$ is trivially bijective, it remains only to show that it is a homeomorphism. We first prove continuity. Let $(\alpha_m, z_m)_{m\in\mathbb{N}}$ be a sequence in $\mathrm{UQD} \times \mathbb{H}$ converging to some limit $(\alpha_\infty, z_\infty)$, say. By Lemma 3.4.3, the sequence $(\Phi([\alpha_m, z_m]))_{m\in\mathbb{N}} = ((\mathrm{Int} \circ \partial)(\alpha_m), \alpha_m(z_m))_{m\in\mathbb{N}}$ converges to $\Phi([\alpha_\infty, z_\infty]) = ((\mathrm{Int} \circ \partial)(\alpha_\infty), \alpha_\infty(z_\infty))$, and continuity follows.

We now prove continuity of its inverse. Let $(D_m, z_m)_{m\in\mathbb{N}}$ be a sequence in $\mathrm{MQD}$ converging to $(D_\infty, z_\infty)$, say, and, for all $m$, denote $c_m := \partial D_m$. By definition of the colimit topology, there exists $K \geqslant 1$ such that $D_m \in \mathrm{QD}(K)$ and $c_m \in \mathrm{QC}^+(K)$ for all $m$. Note now that the restriction of $\partial$ to $\mathrm{QD}(K)$ is the inverse of the restriction of $\mathrm{Int}$ to $\mathrm{QC}^+(K)$. We now apply a compactness argument to show that the former is continuous, and, in particular, that $(c_m)_{m\in\mathbb{N}}$ converges to $c_\infty$. To see this, suppose the contrary. By Theorem 3.3.4, upon extracting a subsequence, we may suppose that $(c_m)_{m\in\mathbb{N}}$ converges to $c'_\infty$, say, which is either a $K$-quasicircle not equal to $c_\infty$, or a single point, $z_\infty$, say. It then follows by continuity of the operator $\overline{\mathrm{Int}}$ that $(\overline{D}_m)_{m\in\mathbb{N}}$ converges in the Hausdorff sense to one of $\overline{\mathrm{Int}}(c'_\infty)$, $\{z_\infty\}$, or $\hat{\mathbb{C}}$. Since none of these is equal to $\overline{D}_\infty$, this is absurd, and it follows that $(c_m)_{m\in\mathbb{N}}$ indeed converges to $c_\infty$, as asserted.

For all $m$, let $\xi_m$ be a non-trivial tangent vector to $\hat{\mathbb{C}}$ at $z_m$ and suppose that $(\xi_m)_{m\in\mathbb{N}}$ converges to $\xi_\infty$. For all $m$, let $\alpha_m \in \mathrm{UQD}(K)$ denote the unique element such that $\alpha_m(\mathbb{H}) = D_m$, $\alpha_m(i) = z_m$ and $D\alpha_m(i) \cdot \partial_x \propto \xi_m$. For all $m$,

$$\partial \alpha_m = \partial \alpha_m(\mathbb{H}) = \partial D_m = c_m.$$

Since $(c_m)_{m\in\mathbb{N}}$ converges to $c_\infty$, and since $(z_m)_{m\in\mathbb{N}}$ does not accumulate on this curve, it follows by Lemma 3.4.4 that $(\alpha_m)_{m\in\mathbb{N}}$ is relatively compact. Since this sequence accumulates only on $\alpha_\infty$, it therefore converges towards this element, so that

$$\lim_{m\to\infty} \Phi^{-1}(D_m, z_m) = \lim_{m\to\infty} [\alpha_m, i] = [\alpha_\infty, i] = \Phi^{-1}(D_\infty, z_\infty).$$

Continuity of $\Phi^{-1}$ follows, and this completes the proof. $\square$

**3.5 - Principal $\mathrm{PSL}(2,\mathbb{R})$-bundles II.** We now show that, for all $0 < k < 1$, $\partial_\infty : \mathrm{UKD}_{k,h} \to \mathrm{QC}^+$ is continuous and that $(\mathrm{UKD}_{k,h}, \partial_\infty)$ is a principal $\mathrm{PSL}(2,\mathbb{R})$-bundle over $\mathrm{QC}^+$.

**Lemma 3.5.1**

*For all $0 < k < 1$, $\partial_\infty : \mathrm{UKD}_{k,h} \to \mathrm{QC}^+$ is continuous.*

**Proof:** Let $(\alpha_m)_{m\in\mathbb{N}}$ be a sequence in $\mathrm{UKD}_{k,h}$ converging to $\alpha_\infty$, say. By definition of the colimit topology, there exists $K \geqslant 1$ such that, for all $m \in \mathbb{N}$, $\alpha_m \in \mathrm{UKD}_{k,h}(K)$. For all $m \in \mathbb{N}\cup\{\infty\}$, denote $c_m := \partial_\infty \alpha_m$ and $K_m := \mathrm{K}_{k,h}(c_m)$, and note that

$$\alpha_m(\mathbb{H}) = D_m := \partial_{\mathrm{fin}} K_m.$$

In particular, for all such $m$,

$$\alpha_m(i) \in D_m.$$

Since $(\alpha_m(i))_{m\in\mathbb{N}}$ converges to $\alpha_\infty(i)$, the sequence $(D_m)_{m\in\mathbb{N}}$ does not accumulate on a point of $\partial_\infty X$. It follows by Lemmas 2.5.1 and 2.5.2 that the sequence $(c_m)_{m\in\mathbb{N}}$ also does not accumulate on a point. By Theorem 3.3.4, it is thus precompact in $\mathrm{QC}^+$. Let $c'_\infty$ be an accumulation point of this sequence. We claim that $c'_\infty = c_\infty$. Indeed, by continuity, $(K_m)_{m\in\mathbb{N}}$ subconverges towards $\mathrm{K}_{k,h}(c'_\infty)$. By Lemma 2.2.2, $(\hat{\partial} K_m)_{m\in\mathbb{N}}$ subconverges towards $\hat{\partial}\mathrm{K}_{k,h}(c'_\infty)$. Since $(\alpha_m)_{m\in\mathbb{N}}$ converges to $\alpha_\infty$, upon extracting a subsequence, we obtain

$$\partial_{\mathrm{fin}} K_\infty = \alpha_\infty(\mathbb{H}) \subseteq \lim_{m\to\infty} \hat{\partial} K_m \cap X = \hat{\partial}\mathrm{K}_{k,h}(c'_\infty) \cap X = \partial_{\mathrm{fin}}\mathrm{K}_{k,h}(c'_\infty),$$

from which it follows that

$$K_\infty = \mathrm{K}_{k,h}(c'_\infty),$$

and

$$c'_\infty = \partial_\infty \alpha_\infty = c_\infty,$$

as asserted. The sequence $(c_m)_{m\in\mathbb{N}}$ is thus precompact in $\mathrm{QC}^+$ and accumulates only on $c_\infty$. It thus converges towards $c_\infty$, and continuity of $\partial_\infty$ follows. $\square$

**Lemma 3.5.2**

*Let $(\alpha_m)_{m\in\mathbb{N}}$ be a sequence in $\mathrm{UKD}_{k,h}(K)$ for some $K \geqslant 1$ and, for all $m$, denote $c_m := \partial_\infty \alpha_m$. If $(c_m)_{m\in\mathbb{N}}$ converges to $c_\infty$, say, and if $(\alpha_m(i))_{m\in\mathbb{N}}$ does not accumulate on $\partial_\infty X$, then $(\alpha_m)_{m\in\mathbb{N}}$ subconverges in the $C^\infty_{\mathrm{loc}}$ sense to some element $\alpha_\infty$, say, such that $c_\infty = \partial_\infty \alpha_\infty$.*

**Proof:** Upon extracting a subsequence, we may suppose that $(\alpha_m(i))_{m\in\mathbb{N}}$ converges to some point $p_\infty$, say, in $X$. Note now that, by uniqueness in Theorem 2.4.1, $\Pi$ acts cocompactly on $\mathrm{UKD}_{k,h}(K)$. Thus, by Theorem 2.6 of [3], $(\alpha_m)_{m\in\mathbb{N}}$ subconverges in the $C^\infty_{\mathrm{loc}}$ sense to a conformal parametrization $\alpha_\infty$, say, of $\mathrm{D}_{k,h}(c_\infty)$, and the result follows. $\square$

We are now ready to prove Theorems 3.1.2 and 3.2.2.

**Proof of Theorem 3.1.2:** Since continuity is proven in Lemma 3.5.1, it remains only to prove the local trivialization property. Choose $c_0 \in \mathrm{QC}^+$, $x \in \mathrm{D}_{k,h}(c_0)$, and a non-trivial tangent vector $\xi_x \in \mathrm{T}_x \mathrm{D}_{k,h}(c_0)$.

Let $U$ be a neighbourhood of $c_0$ in $QC^+$ such that, for all $c \in U$, the closest point projection onto $D_{k,h}(c)$ is smooth near $x$ and the image of $\xi_x$ under its derivative is non-trivial. For all such $c$, let $\xi_c$ denote this image, and observe that this vector varies continuously with $c$. For all $c \in U$, let $\alpha_c$ denote the unique conformal parametrization of $D_{k,h}(c)$ such that $\alpha_c(i) = x$ and $D\alpha_c(i) \cdot \partial_x \propto \xi_c$. By Lemma 3.5.2, $\alpha_c$ depends continuously on $c$. We now define $A : U \times PSL(2,\mathbb{R}) \to UKD_{k,h}$ by

$$A(c, \beta) := \alpha_c \circ \beta.$$

Trivially, this function is continuous and maps $U \times PSL(2,\mathbb{R})$ bijectively onto $\partial_\infty^{-1}(U)$. By Lemma 3.5.1, its inverse is also continuous. We have thus constructed a local trivialization of $(UKD_{k,h}, \partial_\infty)$ near $c_0$, and this completes the proof. $\square$

**Proof of Theorem 3.2.2:** Since $\Phi_{k,h}$ is trivially bijective, it remains only to prove that it is a homeomorphism. We first prove continuity. Let $(\alpha_m, z_m)_{m \in \mathbb{N}}$ be a sequence in $UKD_{k,h} \times \mathbb{H}$ converging to some limit $(\alpha_\infty, z_\infty)$, say. By Theorem 2.1.1 and Lemma 3.5.1, $(\Phi_{k,h}([\alpha_m, z_m]))_{m \in \mathbb{N}} = ((D_{k,h} \circ \partial_\infty)(\alpha_m), \alpha_m(z_m))_{m \in \mathbb{N}}$ converges to $\Phi_{k,h}([\alpha_\infty, z_\infty]) = ((D_{k,h} \circ \partial_\infty)(\alpha_\infty), \alpha_\infty(z_\infty))$, and continuity follows.

We now prove continuity of its inverse. Let $(D_m, p_m)_{m \in \mathbb{N}}$ be a sequence in $MKD_{k,h}$ converging to $(D_\infty, p_\infty)$, say, and, for all $m$, denote $c_m := \partial_\infty D_m$. By definition of the colimit topology, there exists $K \geqslant 1$ such that $D_m \in KD_{k,h}(K)$ and $c_m \in QC^+(K)$ for all $m$. Note now that the restriction of $\partial_\infty$ to $KD_{k,h}(K)$ is the inverse of the restriction of $D_{k,h}$ to $QC^+(K)$. We now apply a compactness argument to show that the former is continuous, and, in particular, that $(c_m)_{m \in \mathbb{N}}$ converges to $c_\infty$. Indeed, suppose that contrary. By Theorem 3.3.4, upon extracting a subsequence, we may suppose that $(c_m)_{m \in \mathbb{N}}$ converges to $c'_\infty$, say, which is either a $K$-quasicircle not equal to $c_\infty$, or a single point, $z_\infty$, say. It then follows by Theorem 2.1.1 and Lemmas 2.5.1 and 2.5.2 that either $(D_m)_{m \in \mathbb{N}}$ converges in the $C^\infty_{\text{loc}}$ sense to $D_{k,h}(c'_\infty)$ or $(\overline{D}_m)_{m \in \mathbb{N}}$ converges in the Hausdorff sense to $\{z_\infty\}$. Since none of these is equal to $D_\infty$, this is absurd, and it follows that $(c_m)_{m \in \mathbb{N}}$ indeed converges to $c_\infty$, as asserted.

Let $\xi_\infty$ be a non-vanishing tangent vector to $D_\infty$ at $p_\infty$. There exists $M \in \mathbb{N}$ such that, for all $m \geqslant M$, the closest point projection onto $D_m$ is well-defined and smooth near $p_\infty$. For all $m \geqslant M$, let $p'_m$ and $\xi'_m$ denote the respective images of $p_\infty$ and $\xi_\infty$ under this projection, and let $\xi_m$ denote the parallel transport of $\xi'_m$ along the unique geodesic in $D_m$ from $p'_m$ to $p_m$. We may suppose that, for all $m \geqslant M$, $\xi_m \neq 0$, and, for all such $m$, we denote by $\alpha_m : \mathbb{H} \to D_m$ the unique conformal parametrization such that $\alpha_m(i) = p_m$ and $D\alpha_m(i) \cdot \partial_x \propto \xi_m$. Since $(c_m)_{m \in \mathbb{N}}$ converges to $c_\infty$, and since $(p_m)_{m \in \mathbb{N}}$ does not accumulate on $\partial_\infty X$, by Lemma 3.5.2, the sequence $(\alpha_m)_{m \in \mathbb{N}}$ is relatively compact in $UKD_{k,h}(K)$. Since this sequence only accumulates on $\alpha_\infty$, it therefore converges towards this element, so that

$$\lim_{m \to \infty} \Phi_{k,h}^{-1}(D_m, p_m) = \lim_{m \to \infty}[\alpha_m, i] = [\alpha_\infty, i] = \Phi_{k,h}^{-1}(D_\infty, p_\infty).$$

Continuity of $\Phi_{k,h}^{-1}$ follows, and this completes the proof. $\square$

## 4 - Equidistribution.

**4.1 - Overview.** In Section 5 of [24], Labourie presents a reformulation of the work [8] of Calegari–Marques–Neves in terms of an equidistribution result for a certain class of measures over $QC^+$. He achieves this by adapting to the 2-dimensional setting the theory of invariant measures and conformal currents developed by Bonahon in [6] (c.f. also [10]). This equidistribution result can in fact be viewed as a synthetic characterization of the sequences constructed by Kahn–Marković in [18] which encapsulates many of their essential properties. In particular, as we will see Section 5, with Kahn–Marković's result expressed in this manner, it is straightforward to derive interesting corollaries, including, but not limited to, the main result of this paper, namely Theorem 1.1.1.

We now show how the framework of fibred Plateau problems developed in the preceding sections provides the natural context for the development of Labourie's ideas. Let $X := (X, h)$ be a Cartan–Hadamard manifold of sectional curvature bounded above by $-1$ acted upon cocompactly by the group $\Pi$. By hyperbolization, $X$ also carries a $\Pi$-invariant hyperbolic metric which we denote by $hyp$.

For every positive integer $g$, and for all $K \geqslant 1$, let $\mathrm{QC}_g^+(K)$ denote the space of oriented $K$-quasicircles $c$ whose stabilizer $\Gamma(c)$ in $\Pi$ is a compact surface subgroup of genus at most $g$. For all $K \geqslant 1$, we denote

$$\mathrm{QC}_*^+(K) := \bigcup_{g \geqslant 1} \mathrm{QC}_g^+(K). \tag{4.1}$$

We denote

$$\mathrm{QC}_*^+ := \mathrm{QC}_*^+(\infty) := \bigcup_{K \geqslant 1} \mathrm{QC}_*^+(K). \tag{4.2}$$

For any $c \in \mathrm{QC}^+$, let $\delta(c)$ denote the Dirac measure over $\mathrm{QC}^+$ supported on $c$, that is, for every Borel subset $A \subseteq \mathrm{QC}^+$,

$$\delta(c)(A) := \begin{cases} 1 & \text{if } c \in A, \text{ and} \\ 0 & \text{otherwise.} \end{cases} \tag{4.3}$$

For all $c \in \mathrm{QC}_*^+$, we define the $\Pi$-invariant, Borel regular measure $\mu(c)$ over $\mathrm{QC}^+$ by

$$\mu(c) := \frac{1}{2\pi |\chi(\Gamma(c))|} \sum_{\alpha \in \Pi/\Gamma(c)} \delta(\sigma(\alpha) \cdot c), \tag{4.4}$$

where $\chi(\Gamma(c))$ denotes the Euler characteristic of $\Gamma(c)$, and $\sigma : \Pi/\Gamma(c) \to \Pi$ is any section.

Finally, recall that, given any measure $m$, its *support* $\mathrm{Supp}(m)$ is defined to be the complement of the union of all open sets of measure zero. Note, in particular, that this set is closed. Recall from Section 3.1 that, for all $c \in \mathrm{QC}^+$, $\mathrm{QS}(c)$ denotes the quasisymmetry constant of $c$. The desired equidistribution result is now stated as follows.

**Theorem 4.1.1, Equidistribution**

There exists a sequence $(c_m)_{m \in \mathbb{N}}$ in $\mathrm{QC}_*^+$ such that

(1) $(\mathrm{QS}(c_m))_{m \in \mathbb{N}}$ tends to $1$; and

(2) $(\mu(c_m))_{m \in \mathbb{N}}$ weakly converges to a Borel regular measure $\mu_\infty$ over $\mathrm{QC}^+$ satisfying

$$\mathrm{Supp}(\mu_\infty) = \mathrm{QC}^+(1) = \mathrm{C}^+. \tag{4.5}$$

In particular, $(\chi(\Gamma(c_m)))_{m \in \mathbb{N}}$ converges to $-\infty$.

Theorem 4.1.1 will be proven in Sections 4.2 and 4.4.

**4.2 - Invariant measures over bundles.** Theorem 4.1.1 will be proven with the help of a number of auxiliary structures. In particular, since $\mathrm{C}^+/\Pi$ is liable to be non-Hausdorff, it will be helpful to work with invariant measures defined over associated $\mathrm{PSL}(2,\mathbb{R})$-bundles over $\mathrm{QC}^+$.

Let $G$ be a unimodular Lie group, let $F$ be the homogeneous space of a unimodular subgroup $H$ of $G$, and let $\mu_F$ denote its Haar measure. Although the following construction is quite general, we will only be interested in the present paper in the case where $G = \mathrm{PSL}(2,\mathbb{R})$, $F = \mathbb{H}^2$, and $\mu_F$ is the area measure of the hyperbolic metric. Given any measurable space $A$, we say that a measure over $A \times F$ is *G-invariant* whenever it is preserved by the action of $G$ on the second factor. Trivially, $G$-invariant measures over $A \times F$ are precisely those measures $\mu$ of the form

$$\mu := \mu_A \times \mu_F, \tag{4.6}$$

for some measure $\mu_A$ over $A$ which we call the *projection* of $\mu$. Now let $(B, \pi)$ be an associated $G$-bundle over $A$ with fibre $F$. We say that a measure over $B$ is *G-invariant* whenever it has this property over every $G$-trivialization. Trivially, $G$-invariant measures over $(B, \pi)$ are precisely those measures $\mu$ which in every $G$-trivialization $\alpha : U \times F \to B$ have the form

$$\alpha^* \mu = (\mu_A|_U) \times \mu_F, \tag{4.7}$$

for some measure $\mu_A$ over $A$ which we likewise call the *projection* of $\mu$.

**Lemma 4.2.1**

*If $A$ is separable and locally compact, then, for every associated $G$-bundle $(B,\pi)$ over $A$, and for every Borel regular measure $\mu$ over $A$, there exists a unique $G$-invariant Borel regular measure $\hat{\mu}$ over $B$ with projection $\mu$. Furthermore, $\hat{\mu}$ depends continuously on $\mu$ with respect to the weak topology.*

**Remark 4.2.1.** Since the converse is trivial, this result yields a homeomorphism between the space of $G$-invariant Borel regular measures over $B$ and the space of Borel regular measures over $A$.

**Remark 4.2.2.** In the present case, where $G = \mathrm{PSL}(2,\mathbb{R})$, $F = \mathbb{H}^2$, and $A = \mathrm{QC}^+$, we will use the fact that $\mathrm{QC}^+(K)$, furnished with the Hausdorff metric, satisfies the hypotheses of Lemma 4.2.1 for all $K$.

**Remark 4.2.3.** $G$-invariant measures over $(B,\pi)$ correspond to the laminar measures introduced by Labourie in [24], which in turn correspond to the invariant measures introduced by Bonahon in [6]. Likewise, measures over $A$ correspond to Labourie's conformal currents, which in turn correspond to Bonahon's geodesic currents.

**Proof:** Let $(U_i)_{i\in I}$ be a covering of $A$ by open sets trivializing $B$. Let $(\phi_i)_{i\in I}$ be a locally finite partition of unity of $A$ subordinate to $(U_i)_{i\in I}$. For each $i$, we identify $\pi^{-1}(U_i)$ with $U_i \times F$ and define $\hat{\mu}_i$ over this set by

$$\hat{\mu}_i := \left(\phi_i \mu|_{U_i}\right) \times \mu_F.$$

Since the measure

$$\hat{\mu} := \sum_{i\in I} \hat{\mu}_i$$

trivially has the desired properties, this proves existence. Uniqueness and continuity follow trivially, and this completes the proof. $\square$

For $0 < k < 1$, consider now the associated $\mathrm{PSL}(2,\mathbb{R})$-bundle $\mathrm{MKD}_{k,h}$ over $\mathrm{QC}^+$ introduced in Section 3.2. For all $c \in \mathrm{QC}_*^+$, let $\hat{\delta}_{k,h}(c)$ and $\hat{\mu}_{k,h}(c)$ denote respectively the unique $\mathrm{PSL}(2,\mathbb{R})$-invariant measures over $\mathrm{MKD}_{k,h}$ with respective projections $\delta(c)$ and $\mu(c)$. In particular, by (4.4),

$$\hat{\mu}_{k,h}(c) := \frac{1}{2\pi\,|\chi(\Gamma(c))|} \sum_{\alpha\in\Pi/\Gamma(c)} \hat{\delta}_{k,h}(\sigma(\alpha)\cdot c), \tag{4.8}$$

where $\sigma : \Pi/\Gamma(c) \to \Pi$ is any section.

At this stage, it is worth noting that there are, in fact, two natural ways to lift each $\delta(c)$ to a measure over $\mathrm{MKD}_{k,h}$. The first, which we have already introduced, weights the fibre over $c$ uniformly with respect to its conformal structure, whilst the second weights this fibre uniformly with respect to its riemannian structure. More precisely, recalling that $\mathrm{MKD}_{k,h}$ is a subset of $\mathrm{KD}_{k,h} \times X$, we denote by $\pi : \mathrm{MKD}_{k,h} \to X$ the projection onto the second factor, and we note that $\pi$ restricts to a smooth embedding of every fibre into $X$. For every such fibre $D$, we denote by $\mathrm{Area}_{k,h}$ the area measure that it inherits from $h$ through $\pi$, and we define the measure $\delta(D)$ over $\mathrm{MKD}_{k,h}$ such that, for every Borel subset $A$ of $\mathrm{MKD}_{k,h}$,

$$\delta(D)(A) = \mathrm{Area}_{D,h}(A \cap D). \tag{4.9}$$

We call $\delta(D)$ the *riemannian Dirac measure* of $\mathrm{MKD}_{k,h}$ supported on $D$, and we see that the second natural lift of $\delta(c)$ is given by $\delta(\mathrm{D}_{k,h}(c))$. In what follows, it will be useful to compare the above two lifts of $\delta(c)$.

**Lemma 4.2.2**

*For $C \geqslant 1$, if the sectional curvature of $h$ is pinched between $-C$ and $-1$, then, for all $0 < k < 1$, and for all $c \in \mathrm{QC}_*^+$,*

$$(1-k)\delta(\mathrm{D}_{k,h}(c)) \leqslant \hat{\delta}_{k,c}(c) \leqslant (C-k)\delta(\mathrm{D}_{k,h}(c)). \tag{4.10}$$

*In particular, for all such $k$ and $c$,*

$$\frac{(1-k)}{2\pi\,|\chi(\Gamma(c))|} \sum_{\alpha\in\Pi/\Gamma(c)} \delta(\mathrm{D}_{k,h}(\sigma(\alpha)\cdot c)) \leqslant \hat{\mu}_{k,h}(c) \leqslant \frac{(C-k)}{2\pi\,|\chi(\Gamma(c))|} \sum_{\alpha\in\Pi/\Gamma(c)} \delta(\mathrm{D}_{k,h}(\sigma(\alpha)\cdot c)), \tag{4.11}$$

where $\chi(\Gamma(c))$ denotes the Euler characteristic of $\Gamma(c)$ and $\sigma : \Pi/\Gamma \to \Pi$ is any section. In particular, when $h$ is hyperbolic, equality holds in (4.11).

**Proof:** Consider first a fibre $D := D_{k,h}(c)$ of $MKD_{k,h}$, and let $h_\pi := \pi^* h$ denote the metric that it inherits from $h$ by the projection $\pi : D \to X$. By Gauss' formula, $h_\pi$ has sectional curvature pinched between $-(C-k)$ and $-(1-k)$. Let $h_{-1}$ denote the unique complete hyperbolic metric in the conformal class of $h_\pi$. It follows from the curvature bounds on $h_\pi$ (see, for example, Theorem $A$ of [1] and Theorem 2 of [37]) that

$$(1-k)h_\pi \leqslant h_{-1} \leqslant (C-k)h_\pi.$$

Denoting by $\text{Area}_\pi$ and $\text{Area}_{-1}$ the respective area measures of $h_\pi$ and $h_{-1}$, it follows that

$$(1-k)\text{Area}_\pi \leqslant \text{Area}_{-1} \leqslant (C_k)\text{Area}_\pi.$$

Consequently,

$$(1-k)\delta(D) \leqslant \hat{\delta}_{k,c}(c) \leqslant (C-k)\delta(D),$$

as desired. Finally (4.11) follows by (4.4), and this completes the proof. $\square$

With $\pi : MKD_{k,h} \to X$ as above, for all $0 < k < 1$, and for all $c \in QC_*^+$, we define the $\Pi$-invariant measure $\hat{\nu}_{k,h}(c)$ over $X$ by

$$\hat{\nu}_{k,h}(c) := \pi_* \hat{\mu}_{k,h}(c). \tag{4.12}$$

Note that $\hat{\nu}_{k,h}(c)$ projects to a measure over the quotient $X/\Pi$ which we denote by $\nu_{k,h}(c)$, and which is none other than the conformal Dirac measure of the compact surface $\pi(D_{k,h}(c))$ – that is, its Dirac measure with respect to its Poincaré metric – normalized so as to have unit total mass. The different measures constructed so far are summarized in Figure 4.2.9. Note, in particular, that the operation sending measures over $QC^+$ to (possibly infinite) measures over $X$, indicated in Figure 4.2.9 by the broken arrow, is a well-defined function which is continuous with respect to the weak topology.

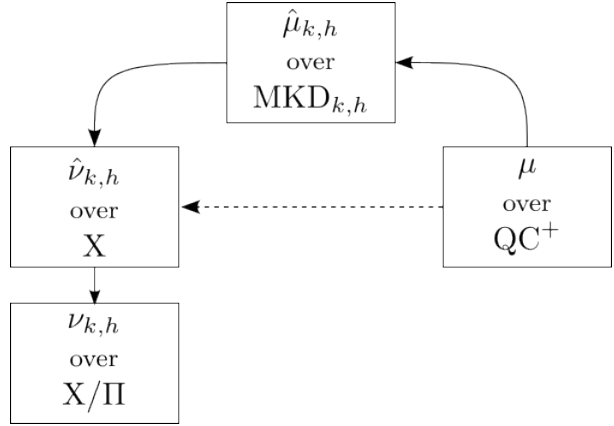

**Figure 4.2.9 - Invariant measures -** We aim to construct invariant measures over different spaces. The measure $\hat{\mu}_{k,h}$, defined over $MKD_{k,h}$, will be $(PSL(2,\mathbb{R}),\Pi)$-biinvariant, whilst the measures $\hat{\nu}_{k,h}$ and $\mu$, defined respectively over $X$ and $QC^+$, will both be $\Pi$-invariant.

**Lemma 4.2.3**

For $0 < k < 1$, and for any sequence $(c_m)_{m \in \mathbb{N}} \in QC_*^+$ such that $(QS(c_m))_{m \in \mathbb{N}}$ converges to $1$, the sequence $(\hat{\mu}_{k,h}(c_m))_{m \in \mathbb{N}}$ subconverges in the weak sense over $MKD_{k,h}$ to a locally finite, Borel regular measure $\hat{\mu}_\infty$ such that

$$\text{Supp}(\hat{\mu}_\infty) \subseteq MKD_{k,h}(1). \tag{4.13}$$

**Proof:** Let $U \subseteq X$ be an open set such that, for all $\gamma \in \Pi$,

$$(\gamma \cdot U) \cap U = \emptyset.$$

By the preceding discussion, for all $m$,

$$\hat{\mu}_{k,h}(c_m)(\pi^{-1}(U)) = \hat{\nu}_{k,h}(c_m)(U) \leqslant 1.$$

The sequence $(\hat{\mu}_{k,h}(c_m))_{m\in\mathbb{N}}$ is thus locally uniformly bounded. Now choose $K > 1$. Since the quasisymmetry constants of $(c_m)_{m\in\mathbb{N}}$ converge to 1, for sufficiently large $m$,

$$\mathrm{Supp}(\hat{\mu}_{k,h}(c_m)) \subseteq \mathrm{MKD}_{k,h}(K).$$

Since $\mathrm{MKD}_{k,h}(K)$ is metrizable and separable, it follows that $(\hat{\mu}_{k,h}(c_m))_{m\in\mathbb{N}}$ weakly subconverges to some Borel regular measure $\hat{\mu}_\infty$, say, such that

$$\mathrm{Supp}(\hat{\mu}_\infty) \subseteq \bigcap_{K\geqslant 1} \mathrm{MKD}_{k,h}(K) = \mathrm{MKD}_{k,h}(1),$$

as desired. $\square$

This yields the first part of Theorem 4.1.1.

**Lemma 4.2.4**

*For all $0 < k < 1$, and for any sequence $(c_m)_{m\in\mathbb{N}} \in \mathrm{QC}^+_*$ such that $(\mathrm{QS}(c_m))_{m\in\mathbb{N}}$ converges to 1, the sequence $(\mu(c_m))_{m\in\mathbb{N}}$ subconverges in the weak sense over $\mathrm{QC}^+$ to a locally finite, $\Pi$-invariant, Borel regular measure $\mu_\infty$ such that*

$$\mathrm{Supp}(\mu_\infty) \subseteq \mathrm{QC}^+(1) = \mathrm{C}^+. \tag{4.14}$$

**Proof:** Let $\hat{\mu}_\infty$ be as in Lemma 4.2.3. Since $\hat{\mu}_\infty$ is locally finite and $(\mathrm{PSL}(2,\mathbb{R}), \Pi)$-biinvariant, it projects to a locally finite, $\Pi$-invariant, Borel regular measure $\mu_\infty$ supported over $\mathrm{C}^+$. By continuity of the projection, $(\mu(c_m))_{m\in\mathbb{N}}$ weakly subconverges to $\mu_\infty$, and the result follows. $\square$

**4.3 - Kahn–Marković sequences.** In [18], Kahn–Marković develop a method of constructing families of curves in $\mathrm{QC}^+_*$ with various useful properties. We will use this method to construct a sequence satisfying the conclusion of Theorem 4.1.1. The reader is doubtless aware of the technical complexity of Kahn–Marković's construction, and we recommend the excellent presentations of this challenging work given by Hamenstädt in [17] and Kassel in [20]. Since Kahn–Marković's construction is purely hyperbolic, we will only be concerned in this and the next section with the metric hyp. In what follows, $\lambda$ will denote volume measure of $X/\Pi$ with respect to this metric, normalized so as to have unit total mass. Likewise, $\hat{\lambda}$ will denote the volume measure of the Sasaki metric that hyp induces over $SX/\Pi$, also normalized so as to have unit total mass. Finally, we will denote by $Y$ the union of all cocompact totally geodesic surfaces in $X$. Note that, since there are only countably many such surfaces, this set has vanishing Lebesgue measure.

Let $(c_m)_{m\in\mathbb{N}}$ be a sequence in $\mathrm{QC}^+_*$. For all $m$, let $K_m := \mathrm{QS}(c_m)$ denote the quasisymmetry constant of $c_m$, and let $\Gamma_m := \Gamma(c_m)$ denote the stabilizer of $c_m$ in $\Pi$. We now describe the many properties that this sequence will be required to possess. We will require that

(1) $(K_m)_{m\in\mathbb{N}}$ tends to 1.

We will require that, for all $m$, there exists a piecewise smooth embedded disk $S_m \subseteq (X, \mathrm{hyp})$ spanning $c_m$ such that

(2) for all $m$, $S_m$ is $\Gamma_m$-invariant.

We will require that

(3) $(\mathrm{Area}(S_m/\Gamma_m)/|\chi(\Gamma_m)|)_{m\in\mathbb{N}}$ converges to $2\pi$.

Given a piecewise smooth, properly embedded surface $S \in (X, \mathrm{hyp})$, let $d(S)$ denote the riemannian Dirac measure over $X$ supported on $S$, that is, for every Borel subset $A \subseteq X$,

$$d(S)(A) := \mathrm{Area}(S \cap A), \tag{4.15}$$

where $S$ is furnished with the area form induced by hyp. We will require that, for every sequence $(\gamma_m)_{m\in\mathbb{N}} \in \Pi$ such that $(\gamma_m \cdot S_m)_{m\in\mathbb{N}}$ meets some fixed bounded set,

(4) $(\gamma_m \cdot S_m)_{m\in\mathbb{N}}$ Hausdorff subconverges over $\overline{X}$ to a complete, totally geodesic surface; and

(5) $(d(\gamma_m \cdot S_m))_{m\in\mathbb{N}}$ weakly subconverges over $X$ to the riemannian Dirac measure of the same complete, totally geodesic surface.

Finally, for all $m$, let $\hat{\nu}'_m$ denote the Borel regular measure over $X$ given by

$$\hat{\nu}'_m := \frac{1}{\text{Area}(S_m/\Gamma_m)} \sum_{\alpha \in \Pi/\Gamma_m} d(\sigma(\alpha) \cdot S_m), \tag{4.16}$$

where $\sigma : \Pi/\Gamma_m \to \Pi$ is any section, and note that this measure projects to a Borel regular probability measure $\nu'_m$ over the quotient space $X/\Pi$. We will require that

(6) $(\hat{\nu}'_m)_{m\in\mathbb{N}}$ weakly subconverges to a Borel regular measure $\hat{\nu}'_\infty$ not wholly supported over $Y$, that is, whose restriction to $X \setminus Y$ is non-trivial.

**Definition 4.3.1**

*Let $(c_m)_{m\in\mathbb{N}}$ be a sequence in $\text{QC}^+_*$. We say that $(c_m)_{m\in\mathbb{N}}$ is a Kahn–Marković sequence whenever it satisfies the Conditions (1) to (6) given above.*

The remainder of this section will be devoted to proving the existence of such sequences. We will then show in the following section that every Kahn-Marković sequence satisfies the conclusion of Theorem 4.1.1.

**Theorem 4.3.2**

*There exists a Kahn–Marković sequence $(c_m)_{m\in\mathbb{N}} \in \text{QC}^+_*$.*

Theorem 4.3.2 follows from the proof of Theorem 4.2 of [8] together with Proposition 6.1 of [27]. However, given the key role that it plays in the present paper, we provide a complete proof for the reader's convenience.

We first recall the relevant elements of Kahn–Marković's construction as presented by Hamenstädt in [17]. Given $\varepsilon > 0$, Hamenstädt constructs a family

$$\mathcal{P}^\varepsilon := (P_i^\varepsilon)_{i \in I} \tag{4.17}$$

of oriented pants in $(X/\Pi, \text{hyp})$, every one of which has geodesic boundary and principal curvatures no greater than $\varepsilon$. Each pant $P_i^\varepsilon$ in fact consists of two hexagons $H_{i,1}^\varepsilon$ and $H_{i,2}^\varepsilon$, centred respectively on the points $c_{i,1}^\varepsilon$ and $c_{i,2}^\varepsilon$, with respective normals $n_{i,1}^\varepsilon$ and $n_{i,2}^\varepsilon$ at these points. Furthermore, every such hexagon contains a disk of radius $r_0$ about its centre, for some $r_0$ independent of $\varepsilon$. Let $\sigma^\varepsilon$ denote the counting measure of the set $\{n_{i,1}^\varepsilon, n_{i,2}^\varepsilon \mid i \in I\}$ over $SX/\Pi$, normalized so as to have unit total mass. That is, for every Borel subset $B \subseteq SX/\Pi$,

$$\sigma^\varepsilon(B) := \sum_{j=1}^{2} \#\left\{ i \mid n_{i,j}^\varepsilon \in B \right\} / 2\#\mathcal{P}^\varepsilon. \tag{4.18}$$

This measure will be of use to us presently. Hamenstädt joins pairs of pants in $\mathcal{P}^\varepsilon$, using every pant in this set exactly once, to yield a family $\Sigma_1^\varepsilon, \cdots, \Sigma_{m_\varepsilon}^\varepsilon$ of compact, piecewise smooth immersed surfaces in $X/\Pi$ with respective positive integer multiplicities $a_1^\varepsilon, \cdots, a_{m_\varepsilon}^\varepsilon$. The key property of this joining construction is that, for each $i$, $\Sigma_i^\varepsilon$ lifts to a piecewise smooth, properly embedded disk in $X$ spanning a $(1+\varepsilon)$-quasicircle.

Let $(i(k))_{k\in\mathbb{N}}$ be a sequence such that, for all $k$, $1 \leqslant i(k) \leqslant m_{1/k}$. For all $k$, let $S_k$ be a lift of $\Sigma_{i(k)}^{1/k}$ and denote $c_k := \partial_\infty S_k$. By the preceding discussion, the sequence $(c_k)_{k\in\mathbb{N}}$ satisfies Conditions (1) and (2) of Definition 4.3.1. Hamenstädt likewise shows that it satisfies Conditions (4) and (5).

Foliated Plateau problems and asymptotic counting of surface subgroups**Lemma 4.3.3**

*The sequence $(c_k)_{k\in\mathbb{N}}$ satisfies Condition (3) of Definition 4.3.1.*

**Proof:** Indeed, by Gauss' theorem, for all $k$, every pant $P$ used in the construction of $S_k$ has intrinsic curvature contained in the interval $[-1-1/k^2, -1+1/k^2]$, so that, by the Gauss–Bonnet theorem,

$$\frac{2\pi}{1+1/k^2} \leqslant \mathrm{Area}(P) \leqslant \frac{2\pi}{1-1/k^2}. \tag{4.19}$$

Since the cardinality of any pair of pants decomposition of $S_k/\Gamma_k$ is equal to $\chi(\Gamma_k)$, this proves (3). $\square$

It remains only to address Condition (6) of Definition 4.3.1. For all $k$, and for all $1 \leqslant i \leqslant m_{1/k}$, let $\omega_i^{1/k}$ denote the area of $\Sigma_i^{1/k}$, and define the probability measure $\nu_{k,i}''$ over $X/\Pi$ by

$$\nu_{k,i}'' := d(\Sigma_i^{1/k})/\omega_i^{1/k}. \tag{4.20}$$

For all $k$, denote

$$\nu_k'' := \sum_{i=1}^{m_{1/k}} a_i^{1/k}\omega_i^{1/k}\nu_{k,i}'' \bigg/ \sum_{i=1}^{m_{1/k}} a_i^{1/k}\omega_i^{1/k}. \tag{4.21}$$

**Lemma 4.3.4**

*The sequence $(\nu_k'')_{k\in\mathbb{N}}$ weakly subconverges to a Borel regular probability measure $\nu_\infty''$ over $X/\Pi$ satisfying*

$$\nu_\infty'' \geqslant b\lambda, \tag{4.22}$$

*for some $0 < b \leqslant 1$.*

**Proof:** Indeed, choose $x \in X$. For all $0 < r < r_0/2$, let $\Omega_{x,r} \subseteq SX/\Pi$ denote the set of all $n := n_y \in SX/\Pi$ such that the totally geodesic plane $P_n$ normal to $n$ at $y$ meets $B_{r/2}(x) \cap B_{r_0/2}(y)$. For sufficiently small $r$ and sufficiently large $k$, if $n_{i,j}^{1/k} \in \Omega_{x,r}$, then

$$\mathrm{Area}(H_{i,j}^{1/k} \cap B_r(x)) \geqslant \frac{1}{8}\pi r^2,$$

where, we recall, $H_{i,j}^{1/k}$ denotes the hexagon centred on the point $c_{i,j}^{1/k}$ with normal $n_{i,j}^{1/k}$ at this point. It follows that, for all such $r$ and $k$,

$$\sum_{i=1}^{m_{1/k}} a_i^{1/k}\omega_i^{1/k}\nu_{k,i}''(B_r(x)) = \sum_{i=1}^{m_{1/k}} a_i^{1/k}d(\Sigma_i^{1/k})(B_r(x))$$

$$= \sum_{P\in\mathcal{P}^{1/k}} \mathrm{Area}(P \cap B_r(x))$$

$$\geqslant \sum_{n_{i,1}^{1/k}\in\Omega_{x,r}} \frac{1}{8}\pi r^2 + \sum_{n_{i,2}^{1/k}\in\Omega_{x,r}} \frac{1}{8}\pi r^2$$

$$= \frac{1}{4}\pi r^2 \sigma^{1/k}(\Omega_{x,r})(\#\mathcal{P}^{1/k}),$$

where here $\sigma^{1/k}$ is as in (4.18). However, Hamenstädt shows that $(\sigma^{1/k})_{k\in\mathbb{N}}$ weakly converges to $\hat{\lambda}$. Bearing in mind (4.19), we therefore obtain, for all sufficiently small $r$,

$$\liminf_{k\to\infty} \frac{1}{r^3}\nu_k''(B_r(x)) \geqslant \frac{1}{16r}\hat{\lambda}(\Omega_{x,r}).$$

Finally, an exercise of classical calculus shows that, for all sufficiently small $r$,

$$\hat{\lambda}(\Omega_{x,r}) \geqslant br,$$

for some $b > 0$, independent of $x$, so that, for all sufficiently small $r$,

$$\liminf_{k \to \infty} \frac{1}{r^3} \nu_k''(B_r(x)) \geqslant \frac{b}{16},$$

and the result follows. $\square$

We underline that, even though $Y$ has zero Lebesgue measure, Lemma 4.3.4 is not in itself sufficient to verify Condition (6) of Definition 4.3.1, since, for all $k$, the measure $\nu_k''$ is constructed using a possibly disconnected family of surfaces, whilst the desired sequence consists of measures, each of which is constructed using a single surface. However, the desired refinement is a straightforward consequence of the following technical result of classical measure theory. Let $\mathcal{P}(X/\Pi)$ denote the space of Borel regular probability measures over $X/\Pi$, furnished with the weak topology. Recall (see, for example, Theorem 4.4 of [33]) that this space is compact. Note also that every convex combination of elements of $\mathcal{P}(X/\Pi)$ is also an element of $\mathcal{P}(X/\Pi)$. We define a *convex polytope* in $\mathcal{P}(X/\Pi)$ to be any closed set which is the convex hull of a finite set of points in $\mathcal{P}(X/\Pi)$. Every convex polytope is generated by a unique set of minimal cardinality whose elements we call its *vertices*.

**Lemma 4.3.5**

*Let $(P_m)_{m \in \mathbb{N}}$ be a sequence of convex polytopes in $\mathcal{P}(X/\Pi)$. For all $m$, let $p_m$ be a point of $P_m$, and suppose that $(p_m)_{m \in \mathbb{N}}$ weakly converges towards $p_\infty \in \mathcal{P}(X/\Pi)$. For every measurable subset $A$ of $X/\Pi$ such that $p_\infty(A) < 1$, there exists a sequence $(p_m')_{m \in \mathbb{N}} \in \mathcal{P}(X/\Pi)$ such that*

*(1) for all $m$, $p_m'$ is a vertex of $P_m$; and*

*(2) $(p_m')_{m \in \mathbb{N}}$ weakly subconverges to a measure $p_\infty'$, say, satisfying $p_\infty'(A) < 1$.*

**Proof:** Choose $\varepsilon > 0$ such that $p_\infty(A) < 1 - 3\varepsilon$. By Borel regularity, there exists a neighbourhood $U$ of $A$ such that $p_\infty(U) < 1 - 2\varepsilon$. Since $X/\Pi$ is separable, there exists a sequence $(B_n)_{n \in \mathbb{N}}$ of open balls in $X/\Pi$, whose closures are contained in $U$, and whose union is $U$. For all $n$, denote $U_n := B_1 \cup \cdots \cup B_n$, let $V_n$ denote the complement of the closure of this set, and note that

$$p_\infty(V_n) = 1 - p_\infty(\overline{U}_n) > 1 - p_\infty(U) > 2\varepsilon.$$

Now fix $n$. By a suitable variant of Fatou's lemma (see, for example, Theorem 1 of Section 1.9 of [11]), for all sufficiently large $m$,

$$p_m(V_n) > \varepsilon,$$

so that

$$p_m(U_n) < 1 - \varepsilon.$$

For all sufficiently large $m$, since $p_m$ is a weighted mean of the vertices of $P_m$, there exists a vertex $p_{m,n}'$ of $P_m$ such that

$$p_{m,n}'(U_n) < 1 - \varepsilon.$$

Since $(U_n)_{n \in \mathbb{N}}$ is an increasing nested family, upon applying a diagonal argument, we obtain a sequence $(p_m')_{m \in \mathbb{N}}$ of vertices such that, for all $n$, and for all sufficiently large $m$,

$$p_m'(U_n) < 1 - \varepsilon.$$

Upon extracting a subsequence, we may suppose that $(p_m')_{m \in \mathbb{N}}$ weakly converges towards some element $p_\infty'$, say, in $\mathcal{P}(X/\Pi)$. By Fatou's lemma again, for all $n$,

$$p_\infty'(U_n) < 1 - \varepsilon.$$

It follows that

$$p_\infty'(A) < p_\infty'(U) = \lim_{n \to \infty} p_\infty'(U_n) \leqslant 1 - \varepsilon,$$

as desired. $\square$

We are now ready to verify Condition (6) of Definition 4.3.1.

**Lemma 4.3.6**

*The sequence $(i(k))_{k\in\mathbb{N}}$ can be chosen in such a manner that the sequence $(\hat{\nu}'_k)_{k\in\mathbb{N}} := (\hat{\nu}''_{k,i(k)})_{k\in\mathbb{N}}$ weakly subconverges to a Borel regular measure $\hat{\nu}'_\infty$ not wholly supported over $Y$.*

**Proof:** Indeed, note that $Y$ is $\Pi$-invariant, and that $Y/\Pi$ also has vanishing Lebesgue measure. Since $\nu''_\infty$ is a probability measure, by (4.22),
$$\nu''_\infty(Y/\Pi) \leqslant 1 - b\lambda(X/\Pi),$$
for some $b > 0$, and hence
$$\nu''_\infty(Y/\Pi) < 1.$$
The result now follows by Lemma 4.3.5. $\square$

We now complete the proof of Theorem 4.3.2.

**Proof of Theorem 4.3.2:** Indeed, as discussed above, Conditions (1), (2), (4) and (5) of Definition 4.3.1 follow readily from Hamenstädt's presentation of Kahn–Marković's construction. Conditions (3) and (6) follow respectively from Lemmas 4.3.3 and 4.3.6. This completes the proof. $\square$

**4.4 - Equidistribution.** We continue to use the notation of the preceding section. The equidistribution property (4.5) will be a consequence of the following result.

**Theorem 4.4.1**

*For every Kahn–Marković sequence $(c_m)_{m\in\mathbb{N}} \in \mathrm{QC}_*^+$, the sequence $(\hat{\nu}_{k,\mathrm{hyp}}(c_m))_{m\in\mathbb{N}}$ subconverges weakly to a Borel regular measure $\hat{\nu}_\infty$ of full support over $X$ which projects to a probability measure over $X/\Pi$.*

We will require the following consequence of the Ratner-Shah Theorem, proven in Theorem 2.6 of [8].

**Theorem 4.4.2**

*If $L \subseteq \mathrm{C}^+$ is a closed $\Pi$-invariant subset, then either $L$ is discrete or $L = \mathrm{C}^+$. Furthermore, if $L$ is discrete, then the stabilizer in $\Pi$ of every point of $L$ is a compact surface subgroup.*

**Sketch of proof:** Note first that, although Theorem 2.5 of [8] is stated for the space C of *unoriented* circles in $\hat{\mathbb{C}}$, this result, and therefore also Theorem 2.6 of [8], also holds for the space $\mathrm{C}^+$ of oriented circles in $\hat{\mathbb{C}}$. Now, if $L \neq \mathrm{C}^+$ then, in particular, no element of $L$ has dense $\Pi$-orbit in $\mathrm{C}^+$. It follows from Theorem 2.6 of [8] that $L$ is discrete, and the stabilizer in $\Pi$ of every point of $L$ is a compact surface subgroup, as desired. $\square$

We prove Theorem 4.4.1 by contradiction. Let $(c_m)_{m\in\mathbb{N}}$ be a Kahn–Marković sequence in $\mathrm{QC}_*^+$ and, for all $m$, denote $\Gamma_m := \Gamma(c_m)$, $D_m := \mathrm{D}_{k,\mathrm{hyp}}(c_m)$ and $\hat{\nu}_m := \hat{\nu}_{k,\mathrm{hyp}}(c_m)$. Note now that, by definition, for all $m$, $D_m$ is $\Gamma_m$-invariant. Note also that, for all $m$, since $D_m$ has constant intrinsic curvature equal to $(k-1)$, by the Gauss–Bonnet theorem,
$$\frac{\mathrm{Area}(D_m/\Gamma_m)}{|\chi(\Gamma_m)|} = \frac{2\pi}{(1-k)}. \tag{4.23}$$

Suppose now that $(\hat{\nu}_m)_{m\in\mathbb{N}}$ converges to the Borel regular measure $\hat{\nu}_\infty$ over $X$. We will show that $\hat{\nu}_\infty$ has full support. To this end, suppose the contrary, and let $\phi \in C_0^0(X)$ be a non-trivial, non-negative valued function such that
$$\mathrm{Supp}(\phi) \cap \mathrm{Supp}(\hat{\nu}_\infty) = \emptyset. \tag{4.24}$$
In particular,
$$\lim_{m\to\infty} \int_X \phi\, d\hat{\nu}_m = \int_X \phi\, d\hat{\nu}_\infty = 0. \tag{4.25}$$
Define $\hat{\phi} \in C^0(X)$ by
$$\hat{\phi}(y) := \sum_{\gamma\in\Pi} (\phi \circ \gamma)(y). \tag{4.26}$$
Note that, since $\phi$ has compact support, this sum is everywhere locally finite.

For all $x \in X$, and for all $\varepsilon, r > 0$, we denote

$$\hat{P}_m(x, \varepsilon, r) := \left\{ \gamma \in \Pi \;\middle|\; \int_{(\gamma \cdot D_m) \cap B_r(x)} \hat{\phi} dA \geqslant \varepsilon \right\}, \text{ and}$$
$$\hat{Q}_m(x, r) := \left\{ \gamma \in \Pi \;\middle|\; (\gamma \cdot S_m) \cap B_r(x) \neq \emptyset \right\},$$
(4.27)

where, for all $m$, $S_m$ is the piecewise smooth embedded disk described in Section 4.3. For all such $x$, $\varepsilon$ and $r$, we denote

$$P_m(x, \varepsilon, r) := \hat{P}_m(x, \varepsilon, r)/\Gamma_m, \text{ and}$$
$$Q_m(x, r) := \hat{Q}_m(x, r)/\Gamma_m,$$
(4.28)

and we denote by $M_m(x, \varepsilon, r)$ and $N_m(x, r)$ the respective cardinalities of these sets.

Heuristically, $M_m(x, \varepsilon, r)$ is the number of copies of $D_m$ in which the restriction of $\hat{\phi}$ to some fixed ball about $x$ has uniformly non-zero mass. In accordance with (4.25), we expect this quantity to grow "slowly" with $m$. Similarly, $N_m(x, r)$ is the number of copies of $S_m$ having non-trivial intersection with some fixed ball about $x$. For suitable balls $B_r(x)$, chosen in accordance with Condition (6) of Definition 4.3.1, we expect this quantity to grow "quickly" with $m$. These assertions, which will be formalized in Lemmas 4.4.3 and 4.4.4, together imply that, for suitable $(x', r')$ and $(x, \varepsilon, r)$, and for sufficiently large $m$,

$$\hat{Q}_m(x', r') \setminus \hat{P}_m(x, \varepsilon, r) \neq \emptyset.$$
(4.29)

from which Theorem 4.4.1 will follow.

We will say that a fundamental domain $\Delta$ of $\Pi$ is *generic* whenever it is polyhedral, and every one of its faces meets $(\gamma \cdot D_m)$ and $(\gamma \cdot S_m)$ transversally for every $\gamma \in \Pi$ and for all $m$. By standard results of transversality theory (see, for example, [16]), generic fundamental domains always exist. Given a generic fundamental domain $\Delta$, for all $m$, and for any section $\sigma : \Pi/\Gamma_m \to \Pi$, we denote

$$\Delta_{m,\sigma} := \bigcup_{\alpha \in \Pi/\Gamma_m} \sigma(\alpha)^{-1}(\Delta).$$
(4.30)

This is a fundamental domain of the action of $\Gamma_m$ over $X$ and, by genericity, $\Delta_{m,\sigma} \cap D_m$ and $\Delta_{m,\sigma} \cap S_m$ are respectively fundamental domains of the action of $\Gamma_m$ over $D_m$ and $S_m$.

**Lemma 4.4.3**

For all $x \in X$, and for all $\varepsilon, r > 0$,

$$\limsup_{m \to \infty} \frac{M_m(x, \varepsilon, r)}{|\chi(\Gamma_m)|} = 0.$$
(4.31)

**Proof:** Let $\Delta$ be a generic fundamental domain of $\Pi$, choose $m \in \mathbb{N}$, and let $\sigma : \Pi/\Gamma_m \to \Pi$ be a section. For all $\gamma \in \Pi$, define the section $\gamma^{-1}\sigma$ such that, for all $\alpha$,

$$(\gamma^{-1}\sigma)(\alpha) := \gamma^{-1} \cdot \sigma(\alpha).$$

For all $\gamma$, since $D_m \cap \Delta_{m,\gamma^{-1}\sigma}$ is a fundamental domain of $\Gamma_m$, by (4.23),

$$\int_X \phi d\hat{\nu}_m = \frac{(1-k)}{2\pi |\chi(\Gamma_m)|} \int_{D_m \cap \Delta_{m,\gamma^{-1}\sigma}} \hat{\phi} dA.$$

Now let $N$ be such that

$$B_r(x) \subseteq \bigcup_{|\gamma| \leqslant N} \gamma(\Delta),$$

where $|\cdot|$ here denotes the word metric of $\Pi$ with respect to some finite generating set, and denote

$$A_1 := \#\{\gamma \in \Pi \mid |\gamma| \leqslant N\}.$$

Then

$$\begin{aligned}
\frac{2\pi |\chi(\Gamma_m)| A_1}{(1-k)} \int_X \phi d\hat{\nu}_m &= \sum_{|\gamma| \leqslant N} \int_{D_m \cap \Delta_{m,\gamma^{-1}\sigma}} \hat{\phi} dA \\
&= \sum_{|\gamma| \leqslant N} \sum_{\alpha \in \Pi/\Gamma_m} \int_{(\sigma(\alpha) \cdot D_m) \cap \gamma(\Delta)} \hat{\phi} dA \\
&= \sum_{\alpha \in \Pi/\Gamma_m} \sum_{|\gamma| \leqslant N} \int_{(\sigma(\alpha) \cdot D_m) \cap \gamma(\Delta)} \hat{\phi} dA \\
&\geqslant \sum_{\alpha \in \Pi/\Gamma_m} \int_{(\sigma(\alpha) \cdot D_m) \cap B_r(x)} \hat{\phi} dA \\
&\geqslant \varepsilon M_m(x, \varepsilon, r).
\end{aligned}$$

It follows by (4.25) that

$$\limsup_{m \to \infty} \frac{M_m(x, \varepsilon, r)}{|\chi(\Gamma_m)|} \leqslant \frac{2\pi A_1}{(1-k)\varepsilon} \limsup_{m \to \infty} \int_X \phi d\hat{\nu}_m = 0,$$

as desired. □

**Lemma 4.4.4**

For all $x \in X$, and for all $r > 0$, such that $\hat{\nu}'_\infty(B_r(x)) > 0$,

$$\liminf_{m \to \infty} \frac{N_m(x, r)}{|\chi(\Gamma_m)|} > 0. \tag{4.32}$$

**Proof:** We may suppose that $r$ is small enough that $B_r(x)$ is contained in some generic fundamental domain $\Delta$ of $\Pi$. For all $m$, let $\sigma_m : \Pi/\Gamma_m \to \Pi$ be a section. Note that, by Condition (5) of Definition 4.3.1, for sufficiently large $m$, and for all $\gamma$,

$$d(\gamma \cdot S_m)(B_r(x)) \leqslant 4\pi(\cosh(r) - 1) =: A_1.$$

For all $m$,

$$\begin{aligned}
\hat{\nu}'_m(B_r(x)) &= \frac{1}{\text{Area}(S_m/\Gamma_m)} \sum_{\alpha \in \Pi/\Gamma_m} d(\sigma_m(\alpha) \cdot S_m)(B_r(x)) \\
&= \frac{1}{\text{Area}(S_m/\Gamma_m)} \sum_{\alpha \in Q_m(x,r)} d(\sigma_m(\alpha) \cdot S_m)(B_r(x)) \\
&\leqslant \frac{A_1 N_m(x, r)}{\text{Area}(S_m/\Gamma_m)}.
\end{aligned}$$

However, by hypotheses,

$$\liminf_{m \to \infty} \hat{\nu}'_m(B_r(x)) \geqslant \hat{\nu}'_\infty(B_r(x)) > 0.$$

so that, by Condition (3) of Definition 4.3.1,

$$\liminf_{m \to \infty} \frac{N_m(x, r)}{|\chi_m(\Gamma_m)|} \geqslant \frac{2\pi \hat{\nu}'_\infty(B_r(x))}{A_1} > 0,$$

as desired. □

**Proof of Theorem 4.4.1.** Choose $x_0 \in X$, and let $C_0^+ \subseteq C^+$ denote the set of all oriented circles $c$ having the property that there exists a sequence $(\gamma_m)_{m \in \mathbb{N}}$ of elements of $\Pi$ and sequences $(\varepsilon_m)_{m \in \mathbb{N}}$ and $(R_m)_{m \in \mathbb{N}}$ converging to $0$ and $\infty$ respectively such that

(1) $(\gamma_m \cdot c_m)_{m \in \mathbb{N}}$ subconverges to $c$; and

(2) for all $m$,
$$\gamma_m \notin \hat{P}_m(x_0, \varepsilon_m, R_m).$$

We claim that $C_0^+ = C^+$. Indeed, suppose the contrary. Since this set is closed and $\Pi$-invariant, by Theorem 4.4.2, $C_0^+$ is discrete, and every element $c$ of $C_0^+$ bounds a cocompact totally geodesic disk in $X$. Now let $Y' \subseteq X$ denote the union of all totally geodesic disks bounded by elements of $C^+$. Note that $Y'$ is a subset of $Y$ and that its complement is open. By Condition (6) of Definition 4.3.1, $\hat{\nu}'_\infty$ has non-trivial mass over $X \setminus Y'$, and there therefore exists $x_1 \in X$ and $r > 0$ such that $\overline{B}_1(x_1)$ lies in the complement of $Y'$ and $\hat{\nu}'_\infty(B_r(x_1)) > 0$. By (4.31) and (4.32), for all $\varepsilon, R > 0$, and for all sufficiently large $m$, there exists an element $\gamma_m^{\varepsilon,R} \in \Pi$ such that
$$\gamma_m^{\varepsilon,R} \notin \hat{P}_m(x_0, \varepsilon, R),$$
and
$$(\gamma_m^{\varepsilon,R} \cdot S_m) \cap B_r(x_1) \neq \emptyset.$$

A diagonal argument then yields sequences $(\varepsilon_m)_{m \in \mathbb{N}}$ and $(R_m)_{m \in \mathbb{N}}$, converging respectively to $0$ and $\infty$, and a sequence $(\gamma_m)_{m \in \mathbb{N}} \in \Pi$ such that, for all sufficiently large $m$,
$$\gamma_m \notin \hat{P}_m(x_0, \varepsilon_m, R_m),$$
and
$$(\gamma_m \cdot S_m) \cap B_r(x_1) \neq \emptyset.$$

By Condition (4) of Definition 4.3.1, we may now suppose that $(\gamma_m \cdot c_m)_{m \in \mathbb{N}}$ converges to an element $c_\infty \in C_0^+$ such that
$$D_{0,\text{hyp}}(c_\infty) \cap \overline{B}_r(x_1) \neq \emptyset.$$

This is absurd, and it follows that $C_0^+ = C^+$, as asserted.

Now choose $c \in C_0^+ = C^+$, let $(\gamma_m)_{m \in \mathbb{N}}$, $(\varepsilon_m)_{m \in \mathbb{N}}$, and $(R_m)_{m \in \mathbb{N}}$ be as above, and denote $D := D_{k,\text{hyp}}(c)$. Choose $\delta, r > 0$. Since $\gamma_m \notin P_m(x_0, \varepsilon_m, R_m)$, for all sufficiently large $m$,
$$\int_{(\gamma_m \cdot D_m) \cap B_r(x_0)} \hat{\phi} \, dA < \delta.$$

Upon taking limits, it follows that
$$\int_{D \cap B_r(x_0)} \hat{\phi} \, dA < \delta.$$

Since $\delta, r > 0$ are arbitrary, it follows that
$$\int_D \hat{\phi} \, dA = 0.$$

Since $\hat{\phi}$ is non-negative, it therefore vanishes over $D$. Since $c \in C^+$ is arbitrary, it follows by Theorem 2.1.3 that $\hat{\phi}$, and therefore also $\phi$, vanishes over $X$. This is absurd, and the result follows. □

The proof of Theorem 4.1.1 now follows.

**Proof of Theorem 4.1.1:** Condition (1) follows by definition of Kahn–Marković sequences. By Lemma Lemma 4.2.4, in order to prove (2), it suffices to verify (4.5). However, suppose the contrary. By Theorem 4.4.2, $\text{Supp}(\hat{\mu}_\infty)$ is discrete, and $\hat{\nu}_\infty$ therefore does not have full support over $X$. This is absurd by Theorem 4.4.1, and (2) follows.

Finally, we show that $(\chi(\Gamma(c_m)))_{m\in\mathbb{N}}$ converges to minus infinity. Indeed, suppose the contrary. Upon extracting a subsequence, we may suppose that $(|\chi(\Gamma(c_m))|)_{m\in\mathbb{N}}$ remains bounded. For all $m$, denote $D_m := D_{k,\mathrm{hyp}}(c_m)$, and denote

$$Z := \bigcup_{m\in\mathbb{N}} D_m.$$

By hypothesis, the surfaces $(D_m/\Gamma(c_m))_{m\in\mathbb{N}}$ have uniformly bounded genus. Furthermore, since $\Pi$ acts cocompactly on $X$, and since $(\mathrm{QS}(c_m))_{m\in\mathbb{N}}$ converges to 1, they also have uniformly bounded shape operators. It follows that this sequence is finite, and that $Z$ is the complement of a dense, open subset of $X/\Pi$. However, since $\nu_m$ is supported over $Z$ for all $m$, so too is $\nu_\infty$. This is absurd, and it follows that $(\chi(\Gamma(c_m)))_{m\in\mathbb{N}}$ converges to minus infinity, as desired. $\square$

## 5 - Asymptotic counting.

**5.1 - Proof of main result.** We conclude this paper by applying the techniques developed in the preceding sections to the study of the marked area spectrum of compact $k$-surfaces in compact, negatively-curved 3-manifolds. Let $(X,h)$ be a Cartan–Hadamard manifold of sectional curvature bounded above by $-1$, acted on cocompactly by the group $\Pi$. We first note that Lemma 4.2.1 yields the following useful reformulation of Theorem 4.1.1.

**Theorem 5.1.1**

There exists a sequence $(c_m)_{m\in\mathbb{N}}$ in $\mathrm{QC}^+$ such that

(1) $(\mathrm{K}(c_m))_{m\in\mathbb{N}}$ converges to 1; and

(2) for all $0 < k < 1$, $(\hat{\mu}_{k,h}(c_m))_{m\in\mathbb{N}}$ weakly subconverges to a $\mathrm{PSL}(2,\mathbb{R})$-invariant Borel regular measure $\hat{\mu}_\infty$ over $\mathrm{MKD}_{k,h}$ satisfying

$$\mathrm{Supp}(\hat{\mu}_\infty) = \mathrm{MKD}_{k,h}(1). \tag{5.1}$$

Following [8], we refine our asymptotic counting problem as follows. Let $\mathrm{QF} := \mathrm{QF}(\Pi)$ denote the set of conjugacy classes of quasi-Fuchsian surface subgroups of $\Pi$. Recall that, for all $[\Gamma] \in \mathrm{QF}$, the limit set $\partial_\infty \Gamma$ is a quasicircle. For all $[\Gamma] \in \mathrm{QF}$, let $\mathrm{g}([\Gamma])$ denote the genus of $\Gamma$, let $\mathrm{QS}([\Gamma])$ denote the quasisymmetry constant of $\partial_\infty \Gamma$, and let $\mathrm{A}_{k,h}([\Gamma])$ denote the area of the surface $\mathrm{D}_{k,h}(\partial_\infty \Gamma)/\Gamma$ in $(X/\Pi, h)$. For all $0 < k < 1$, and for all $K \geqslant 1$, we define

$$\mathrm{Ent}_{k,K}(h) := \liminf_{A\to\infty} \frac{\log\left(\#\{[\Gamma] \in \mathrm{QF} \mid \mathrm{QS}([\Gamma]) \leqslant K \text{ and } \mathrm{A}_{k,h}([\Gamma]) \leqslant A\}\right)}{A\log(A)}. \tag{5.2}$$

Note that this is an increasing function of $K$. We define

$$\mathrm{Ent}_k(h) := \lim_{K\to 1} \mathrm{Ent}_{k,K}(h). \tag{5.3}$$

**Lemma 5.1.2**

If $\mathrm{hyp}$ is a hyperbolic metric over $X$, then, for all $1 < K \leqslant \infty$, and for all $0 < k < 1$,

$$\mathrm{Ent}_{k,K}(\mathrm{hyp}) = \frac{(1-k)}{2\pi}. \tag{5.4}$$

**Proof:** Indeed, by the Gauss-Bonnet theorem, for all $[\Gamma] \in \mathrm{QF}$,

$$\mathrm{A}_{k,\mathrm{hyp}}([\Gamma]) = \frac{2\pi|\chi(\Gamma)|}{(1-k)} = \frac{4\pi(\mathrm{g}([\Gamma])-1)}{(1-k)}.$$

Thus, by (1.5),

$$\mathrm{Ent}_{k,K}(\mathrm{hyp}) = \liminf_{g\to\infty} \frac{(1-k)\log\left(\#\{[\Gamma] \in \mathrm{QF} \mid \mathrm{QS}([\Gamma]) \leqslant K \text{ and } \mathrm{g}([\Gamma]) \leqslant g\}\right)}{4\pi g\log(g)} = \frac{(1-k)}{2\pi},$$

as desired. $\square$

Theorem 1.1.1 will be a consequence of the following stronger result.

**Theorem 5.1.3**

For all $0 < k < 1$,
$$\mathrm{Ent}_k(h) \geqslant \frac{(1-k)}{2\pi}, \qquad (5.5)$$

with equality holding if and only if h has constant sectional curvature equal to $-1$.

**Remark 5.1.1.** Theorem 1.1.1 follows immediately from Theorem 5.1.3 by monotonicity of $\mathrm{Ent}_{k,K}$ in $K$.

Theorem 5.1.3 will follow from the following technical result.

**Lemma 5.1.4**

For any sequence $(c_m)_{m\in\mathbb{N}}$ in $\mathrm{QC}_*^+$ satisying the conclusions of Theorem 5.1.1,
$$\limsup_{m\to\infty} \frac{\mathrm{Area}(\mathrm{D}_{k,h}(c_m)/\Gamma(c_m))}{|\chi(\Gamma(c_m))|} \leqslant \frac{2\pi}{(1-k)}, \qquad (5.6)$$

with equality holding if and only if h has constant sectional curvature equal to $-1$.

**Proof:** Indeed, for all $m$, denote $\Gamma_m := \Gamma(c_m)$, $D_m := \mathrm{D}_{k,h}(c_m)$, and $A_m := \mathrm{Area}(D_m/\Gamma_m)$. For all $m$, we view $D_m$ both as a fibre of $\mathrm{MKD}_{k,h}$ and as a submanifold of $X$. For all $m$, let $h_m$ and $\mathrm{dA}_m$ denote respectively the metric and area form that it inherits from $(X,h)$, and let $\kappa_m$ denote the curvature of $h_m$. We denote by $\sigma: SX \to \mathbb{R}$ the function whose value at the point $\xi_x \in S_xX$ is the sectional curvature of the orthogonal complement of $\xi_x$ in $T_xX$. We also denote by $\sigma: \mathrm{MKD}_{k,h} \to \mathbb{R}$ the function whose value at the point $(D,p)$ is the sectional curvature of $T_pD$ in $X$.

By the Gauss-Bonnet theorem, for all $m$,

$$\begin{aligned} A_m &= \frac{1}{(1-k)} \int_{D_m/\Gamma_m} (1-k)\mathrm{dA}_m \\ &= \frac{1}{(1-k)} \int_{D_m/\Gamma_m} (1-|\sigma|)\mathrm{dA}_m + \frac{1}{(1-k)} \int_{D_m/\Gamma_m} |\kappa_m|\,\mathrm{dA}_m \\ &= \frac{1}{(1-k)} \int_{D_m/\Gamma_m} (1-|\sigma|)\mathrm{dA}_m + \frac{2\pi|\chi(\Gamma_m)|}{(1-k)}. \end{aligned}$$

Since the first term is negative, (5.6) follows. Moreover, equality holds if and only if, up to extraction of a subsequence,
$$\lim_{m\to\infty} \frac{1}{|\chi(\Gamma_m)|} \int_{D_m/\Gamma_m} (|\sigma|-1)\mathrm{dA}_m = 0. \qquad (5.7)$$

Suppose now that equality holds in (5.6), so that (5.7) holds. Let $\pi: \mathrm{MKD}_{k,h} \to X$ denote the canonical projection. Let $\Delta$ be a generic fundamental domain of $\Pi$. Viewing $\sigma$ now as a function over $\mathrm{MKD}_{k,h}$, we claim that
$$\lim_{m\to\infty} \int_{\pi^{-1}(\Delta)} (|\sigma|-1) d\hat{\mu}_{k,h}(c_m) = 0. \qquad (5.8)$$

Indeed, choose $m \in \mathbb{N}$, and let $\tau: \Pi/\Gamma_m \to \Pi$ be a section. By Lemma 4.2.2,

$$\begin{aligned} \int_{\pi^{-1}(\Delta)} (|\sigma|-1)d\hat{\mu}_{k,h}(c_m) &\leqslant \frac{(C-k)}{2\pi|\chi(\Gamma_m)|} \sum_{\alpha\in\Pi/\Gamma_m} \int_{\pi^{-1}(\Delta)} (|\sigma|-1) d\delta(\mathrm{D}_{k,h}(\tau(\alpha)\cdot c_m)) \\ &= \frac{(C-k)}{2\pi|\chi(\Gamma_m)|} \sum_{\alpha\in\Pi/\Gamma_m} \int_{\pi^{-1}(\tau(\alpha)^{-1}(\Delta))} (|\sigma|-1) d\delta(\mathrm{D}_{k,h}(c_m)) \\ &= \frac{(C-k)}{2\pi|\chi(\Gamma_m)|} \int_{\pi^{-1}(\Delta_{m,\tau})} (|\sigma|-1) d\delta(D_m), \end{aligned}$$

where $\Delta_{m,\tau}$ is as in (4.30). In particular, since $\Delta_{m,\tau} \cap D_m$ is a fundamental domain of the action of $\Gamma_m$ over this disk, it follows that

$$\int_{\pi^{-1}(\Delta)} (|\sigma| - 1) d\hat{\mu}_{k,h}(c_m) \leqslant \frac{(C-k)}{2\pi |\chi(\Gamma_m)|} \int_{D_m/\Gamma_m} (|\sigma| - 1) dA_m,$$

and (5.8) follows by (5.7).

By Theorem 5.1.1, $(d\hat{\mu}_{k,h}(c_m))_{m \in \mathbb{N}}$ weakly subconverges to a measure of full support over $\mathrm{MKD}_{k,h}(1)$, so that $(|\sigma| - 1)$ vanishes over this subset. However, by (1.1) and Theorem 2.1.3, $\mathrm{MKD}_{k,h}(1)$ identifies with $SX$. It follows that $h$ has constant sectional curvature equal to $-1$, and this completes the proof. $\square$

Before proceeding with the proof of Theorem 5.1.3, we first note that Lemma 5.1.4 immediately yields the following simpler rigidity result for the marked area spectrum of compact $k$-surfaces in $(X, h)$ which we consider worth stating separately.

**Theorem 5.1.5**

Let $(X, h)$ be a Cartan–Hadamard manifold of sectional curvature bounded above by $-1$, acted on cocompactly by the group $\Pi$. The metric $h$ has constant sectional curvature equal to $-1$ if and only if, for all $[\Gamma] \in \mathrm{QF}(\Pi)$,

$$A_{k,h}([\Gamma]) = \frac{2\pi |\chi(\Gamma)|}{(1-k)}. \tag{5.9}$$

**Proof:** Indeed, if $h$ has constant sectional curvature equal to $-1$, then (5.9) holds for all $[\Gamma] \in \mathrm{QF}$ by the Gauss–Bonnet theorem. Conversely, if (5.9) holds for all $[\Gamma] \in \mathrm{QF}$, then equality holds in (5.6), and $h$ therefore has constant sectional curvature equal to $-1$, as desired. $\square$

We now complete the proof of Theorem 5.1.3.

**Proof of Theorem 5.1.3:** Indeed, as in Lemma 5.1.4, we verify that, for every $[\Gamma] \in \mathrm{QF}$,

$$A_{k,h}([\Gamma]) \leqslant A_{k,\mathrm{hyp}}([\Gamma]),$$

from which it readily follows that

$$\mathrm{Ent}_k(h) \geqslant \mathrm{Ent}_k(\mathrm{hyp}).$$

It only remains to prove rigidity. Suppose that

$$\mathrm{Ent}_k(h) = \mathrm{Ent}_k(\mathrm{hyp}) = \frac{(1-k)}{2\pi}.$$

We claim that, for every $(c_m)_{m \in \mathbb{N}} \in \mathrm{QC}_*^+$ which satisfies the conclusions of Theorem 5.1.1,

$$\limsup_{m \to \infty} \frac{\mathrm{Area}(D_{k,h}(c_m)/\Gamma(c_m))}{|\chi(\Gamma(c_m))|} = \frac{2\pi}{(1-k)}.$$

The result will then follow by Lemma 5.1.4. To this end, suppose the contrary. There exists $\varepsilon > 0$, and a sequence $(c_m)_{m \in \mathbb{N}}$, satisfying the conclusions of Theorem 5.1.1, such that

$$\limsup_{m \to \infty} \frac{\mathrm{Area}(D_{k,h}(c_m)/\Gamma(c_m))}{|\chi(\Gamma(c_m))|} \leqslant \frac{2\pi(1-2\varepsilon)}{(1-k)}.$$

In particular, for all sufficiently large $m$,

$$\mathrm{Area}(D_{k,h}(c_m)/\Gamma(c_m)) \leqslant \frac{2\pi(1-\varepsilon) |\chi(\Gamma(c_m))|}{(1-k)}.$$

Now let $\text{QF}'$ denote the set of all $[\Gamma] \in \text{QF}$ such that $\Gamma$ preserves $\gamma \cdot c_m$, for some $m$, and for some $\gamma \in \Pi$.

$$\begin{aligned}
\text{Ent}_k(h) &= \lim_{K \to 1} \liminf_{A \to \infty} \frac{\log\left(\#\{[\Gamma] \in \text{QF} \mid \text{QS}([\Gamma]) \leqslant K \text{ and } \text{A}([\Gamma]) \leqslant A\}\right)}{A \log(A)} \\
&\geqslant \lim_{K \to 1} \liminf_{A \to \infty} \frac{\log\left(\#\{[\Gamma] \in \text{QF}' \mid \text{QS}([\Gamma]) \leqslant K \text{ and } \text{A}([\Gamma]) \leqslant A\}\right)}{A \log(A)} \\
&\geqslant \lim_{K \to 1} \liminf_{A \to \infty} \frac{\log\left(\#\{[\Gamma] \in \text{QF}' \mid \text{QS}([\Gamma]) \leqslant K \text{ and } \text{g}([\Gamma]) \leqslant (1-k)A/4\pi(1-\varepsilon)\}\right)}{A \log(A)}
\end{aligned}$$

However, upon applying Stirling's approximation to the Müller–Puchta formula (Formula (2) of [28]), we find that, for all $K$, there exists $B_1 > 0$ such that, for all sufficiently large $g$,

$$\#\{[\Gamma] \in \text{QF}' \mid \text{QS}([\Gamma]) \leqslant K \text{ and } \text{g}([\Gamma]) \leqslant g\} \geqslant (B_1 g)^{2g},$$

from which it readily follows that

$$\text{Ent}_k(h) \geqslant \frac{(1-k)}{2\pi(1-\varepsilon)}.$$

This is absurd, the assertion thus follows, and this completes the proof. $\square$

## A - Notation.

FUNCTIONS

| Function | Description |
| --- | --- |
| $\text{A}_{k,h}([\Gamma])$ | The area of the surface $(\text{D}_{k,h}\partial_\infty \Gamma)/\Gamma$ in $(X/\Pi,h)$ |
| $\text{D}_{k,h}(c)$ | The $k$-disk spanned by the Jordan curve $c$ |
| $\partial \alpha$ | The boundary of the image of the univalent map $\alpha$ |
| $\partial_\infty \alpha$ | The ideal boundary of the image of the conformal parametrization $\alpha$ |
| $\partial D$ | The boundary of the topological disk $D$ |
| $\partial_\infty D$ | The ideal boundary of the $k$-disk $D$ |
| $\partial(D,z)$ | The boundary of the marked quasidisk $(D,z)$ |
| $\partial_\infty(D,p)$ | The ideal boundary of the marked $k$-disk $(D,p)$ |
| $\hat{\partial} K$ | The boundary of the convex set $K$ |
| $\partial_{\text{fin}} K$ | The finite component of the boundary of the convex set $K$ |
| $\partial_\infty K$ | The ideal component of the boundary of the convex set $K$ |
| $\text{Ext}(c)$ | The exterior of the oriented Jordan curve $c$ |
| $\overline{\text{Ext}}(c)$ | The closure of the exterior of the oriented Jordan curve $c$ |
| $\text{g}([\Gamma])$ | The genus of the surface subgroup $\Gamma$ |
| $\text{Int}(c)$ | The interior of the oriented Jordan curve $c$ |
| $\overline{\text{Int}}(c)$ | The closure of the interior of the oriented Jordan curve $c$ |
| $\text{QS}(c)$ | The quasisymmetry constant of the quasicircle $c$ |
| $\text{QS}([\Gamma])$ | The quasisymmetry constant of $\partial_\infty \Gamma$ |
| $\Gamma(c)$ | The stabilizer subgroup of the oriented quasicircle $c$ in $\Pi$ |
| $\Phi_{k,h,\alpha}(D,p)$ | The first parametrization of $SX$ by MD |
| $\chi([\Gamma])$ | The Euler characteristic of the surface subgroup $\Gamma$ |
| $\Psi_{k,h,\alpha}(D,p)$ | The second parametrization of $SX$ by MD |

## TOPOLOGICAL SPACES

| Space | Topology | Description |
|---|---|---|
| $C^+$ | Hausdorff | Oriented round circles in $\hat{\mathbb{C}}$ |
| CC | Hausdorff | Compact convex subsets of $\overline{X}$ |
| $CC_{k,h}(\Omega)$ | Hausdorff | Elements $K \subseteq CC$ such that $(\partial_{\text{fin}} K) \cap \Omega$ is $C^\infty$ of constant curvature $k$ |
| $CC_{k,h}$ | Hausdorff | Elements $K \subseteq CC$ such that $\partial_{\text{fin}} K$ is $C^\infty$ of constant curvature $k$ |
| $\partial CC_{k,h}(\Omega)$ | Hausdorff | The topological boundary of $CC_{k,h}(\Omega)$ in CC |
| $\partial CC_{k,h}$ | Hausdorff | The topological boundary of $CC_{k,h}$ in CC |
| $JC^+$ | Hausdorff | Oriented Jordan curves in $\hat{\mathbb{C}}$ |
| $KD_{k,h}(K)$ | $C^\infty_{\text{loc}}$ | $k$-disks in $(X,h)$ spanning $K$-quasicircles |
| $KD_{k,h}$ | Colimit | $k$-disks in $(X,h)$ spanning quasicircles |
| MD | Hausdorff$\times \mathbb{S}^2$ | Marked round disks in $\hat{\mathbb{C}}$ |
| $MKD_{k,h}$ | $KD_{k,h} \times X$ | Marked $k$-disks in $(X,h)$ spanning quasicircles |
| MQD | $QD \times \hat{\mathbb{C}}$ | Marked quasidisks in $\hat{\mathbb{C}}$ |
| $QC^+(K)$ | Hausdorff | Oriented $K$-quasicircles in $\hat{\mathbb{C}}$ |
| $QC^+$ | Colimit | Oriented quasicircles in $\hat{\mathbb{C}}$ |
| $QC^+_*$ | $QC^+$ | Elements $c \in QC^+$ with $\Gamma(c)$ a compact surface subgroup |
| $QC^+_g$ | $QC^+$ | Elements $c \in QC^+_*$ with $\Gamma(c)$ of genus at most $g$ |
| $QD(K)$ | Hausdorff | $K$-quasidisks in $\hat{\mathbb{C}}$ |
| QD | Colimit | Quasidisks in $\hat{\mathbb{C}}$ |
| $QH^+(K)$ | Compact-open | $K$-quasiconformal homeomorphisms of $\hat{\mathbb{C}}$ |
| $QH^+$ | Colimit | Quasiconformal homeomorphisms of $\hat{\mathbb{C}}$ |
| $QH^+(\hat{\mathbb{R}})$ | $QH^+$ | Elements $\alpha \in QC^+$ preserving the oriented, extended real line $\hat{\mathbb{R}}$ |
| $UKD_{k,h}(K)$ | Compact-open | Uniformizations of $k$-disks in $(X,h)$ spanning $K$-quasicircles |
| $UKD_{k,h}$ | Colimit | Uniformizations of $k$-disks in $(X,h)$ spanning quasicircles |
| $UQD(K)$ | Compact-open | Uniformizations of $K$-quasidisks in $\hat{\mathbb{C}}$ |
| UQD | Colimit | Uniformizations of quasidisks in $\hat{\mathbb{C}}$ |

## CLASSICAL SURFACE THEORY

| Object | Formula | Description |
|---|---|---|
| $(S,e)$ | | Immersed surface |
| $\nu_e$ | | Unit normal vector field |
| $I_e$ | $I_e(\xi,\mu) := g(De\cdot\xi, De\cdot\mu)$ | First fundamental form |
| $II_e$ | $II_e(\xi,\mu) := g(\nabla_\xi \nu_e, De\cdot\mu)$ | Second fundamental form |
| $A_e$ | $I_e(A_e\cdot\xi,\mu) := II_e(\xi,\mu)$ | Shape operator |
| $K_e$ | $K_e := \text{Det}(A_e)$ | Extrinsic curvature |

**Bibliography.**

[1] Ahlfors L., An extension of Schwarz's lemma, *Trans. AMS.*, **43**, 359–364, (1938)

[2] Aronszajn N., A unique continuation theorem for elliptic differential equations or inequalities of the second order, *J. Math. Pures Appl.*, **36**, 235–249, (1957)

[3] Alvarez S., Smith G., Prescription de courbure des feuilles des laminations: retour sur un théorème de Candel, *Ann. Inst. Fourier*, **71**, no. 6, 2549–2593, (2021)

[4] Benoist Y., Foulon P., Labourie F., Flots d'Anosov à distributions de Liapounov différentiables I., *Ann. l'I.H.P. Physique Théorique*, **53**, no. 4, 395–412, (1990)

[5] Besson G., Courtois G., Gallot S., Entropies et rigidités des espaces localement symétriques de courbure strictement négative, *G.A.F.A.*, **5**, no. 5, 731–799, (1995)

[6] Bonahon F., Geodesic laminations with transverse Hölder distributions, *Ann. Sci. ENS.*, **30**, no. 2, 205–240, (1997)

[7] Calabi E., An extension of E. Hopf's maximum principle with an application to Riemannian geometry, *Duke Math. J.*, **25**, no. 1, 45-56, (1958)

[8] Calegari D., Marques F. C., Neves A, Counting minimal surfaces in negatively curved 3-manifolds, *Duke Math. J.*, **171**, no. 8., 1615–1648, (2022)

[9] Colding T. H., Minicozzi W. P. II, *A course in minimal surfaces*, Graduate Studies in Mathematics, **121**, Amer. Math. Soc., (2011)

[10] Erlandsson V., Souto J., *Mirzakhani's Curve Counting and Geodesic Currents*, Progress in Mathematics, **345**, Springer–Verlag, (2022)

[11] Evans L. C., Gariepy R. F., *Measure theory and fine properties of functions*, Routledge, (2015)

[12] Gromov M., Foliated Plateau problem, part I: Minimal varieties, *G.A.F.A.*, **1**, 14–79, (1991)

[13] Gromov M., Foliated Plateau problem, part II: Harmonic maps of foliations, *G.A.F.A.*, **1**, 253–320, (1991)

[14] Gromov M., Sign and geometric meaning of curvature, *Rendiconti del Seminario Matematico e Fisico di Milano*, **61**, 9–123, (1991)

[15] Gromov M., Three remarks on geodesic dynamics and fundamental group, *Enseign. Math.*, **46**, 391–402, (2000)

[16] Guillemin V., Pollack A., *Differential Topology*, American Mathematical Society, (1974)

[17] Hamenstädt U., Incompressible surfaces in rank one locally symmetric spaces, *G.A.F.A*, **25**, 815–859, (2015)

[18] Kahn J., Marković V., Immersing almost geodesic surfaces in a closed hyperbolic three-manifold, *Ann. Math.*, **175**, 1127–1190, (2012)

[19] Kahn J., Marković V., Counting essential surfaces in a closed hyperbolic three-manifold, *Geom. Topo.*, **16**, 601–624, (2012)

[20] Kassel F., Groupes de surface dans les réseaux des groupes de Lie semi-simples (d'après J. Kahn, V. Marković, U. Hamenstädt, F. Labourie et S. Mozes), *Séminaire Bourbaki*, **1181**, Astérisque, **438**, 1–72, (2022)

[21] Labourie F., Problèmes de Monge-Ampère, courbes pseudo-holomorphes et laminations, *G.A.F.A.*, **7**, 496–534, (1997)

[22] Labourie F., Un lemme de Morse pour les surfaces convexes, *Inv. Math.*, **141**, no. 2, 239–297, (2000)

[23] Labourie F., Random $k$-surfaces, *Ann. Math.*, **161**, no. 1, 105–140, (2005)

[24] Labourie F., Asymptotic counting of minimal surfaces in hyperbolic 3-manifolds, *Séminaire Bourbaki*, 05/2021, $73^o$ année, no. 1179, 2020–2021

[25] Lehto O., Virtanen K. I., *Quasiconformal mappings in the plane*, Springer-Verlag, (1973)

[26] Lowe B., Deformations of totally geodesic foliations and minimal surfaces in negatively curved 3-manifolds, *G.A.F.A.*, **31**, 895–929, (2021)

[27] Lowe B., Neves A., Minimal surface entropy and average area ratio, arXiv:2110.09451

[28] Müller T. W., Puchta J-C., Character theory of symmetric groups and subgroup growth of surface groups, *J. London Math. Soc.*, **66**, 623–640, (2002)

[29] McDuff D., Salamon D., *J-holomorphic curves and symplectic topology*, Colloquium Publications, **52**, Amer. Math. Soc., (2004)

[30] Ratner M., Raghunathan's topological conjecture and distributions of unipotent flows, *Duke Math. J.*, **63**, 235–280, (1991)

[31] Rudin W., *Real and complex analysis*, McGraw–Hill, (1987)

[32] Seppi A., Minimal discs in hyperbolic space bounded by a quasicircle at infinity, *Comment. Math. Helv.*, **91**, 807–839, (2016)

[33] Simon L., *Geometric measure theory*, Proceedings of the Centre for Mathematical Analysis, **3**, Australian National University, (1984)

[34] Smith G., On the asymptotic Plateau problem in Cartan–Hadamard manifolds, arXiv:2107.14670
[35] Smith G., Quaternions, Monge-Ampère structures and $k$-surfaces, to appear in *Surveys in Geometry, Vol. II*, (Papadopoulos A. ed.), Springer-Verlag, (2023)
[36] Smith G., *Global singularity theory for the Gauss curvature equation*, Ensaios Matemáticos, **28**, 1–114, (2015)
[37] Yau S. T., A general Schwarz lemma for Kähler manifolds, *Amer. J. Math.*, **100**, no. 1, 197–203 (1978)